\documentclass[smallextended]{svjour3}
\pdfoutput=1
\usepackage{epsfig,amsfonts,color}
\usepackage{amsmath}
\usepackage{blkarray}
\usepackage{bm}
\usepackage{algorithm}
\usepackage{algorithmic}
\usepackage{multirow}
\usepackage{graphicx}
\usepackage{caption}
\usepackage{subcaption}
\usepackage{placeins}
\usepackage{marvosym}

\newcommand{\mods}[1]{\left|  #1 \right|}

\usepackage{amssymb, palatino, geometry}
\usepackage[hyphens]{url}
\usepackage[colorlinks=true,linkcolor=blue,citecolor=blue,urlcolor=blue]{hyperref}
\usepackage[titletoc,title]{appendix}

\usepackage{multirow}
\usepackage{mathtools}

\newcommand{\be}{\begin{equation}}
\newcommand{\ee}{\end{equation}}
\newcommand{\bea}{\begin{eqnarray}}
\newcommand{\eea}{\end{eqnarray}}

\newcommand{\bvec}{\left(\begin{array}{c}}
\newcommand{\evec}{\end{array}\right)}
\newcommand{\bsub}{\begin{subequations}}
\newcommand{\esub}{\end{subequations}}

\usepackage{listings}
\DeclareCaptionFont{white}{\color{white}}
\DeclareCaptionFormat{listing}{{\parbox{\dimexpr\textwidth-2\fboxsep\relax}{#1#2#3}}}
\usepackage[dvipsnames]{xcolor}

\lstdefinelanguage{julia}
{
keywordsprefix=\@,
morekeywords={
exit,whos,edit,load,is,isa,isequal,typeof,tuple,ntuple,uid,hash,finalizer,convert,promote,
subtype,typemin,typemax,realmin,realmax,sizeof,eps,promote_type,method_exists,applicable,
invoke,dlopen,dlsym,system,error,throw,assert,new,Inf,Nan,pi,im,begin,while,for,in,return,
break,continue,macro,quote,let,if,elseif,else,try,catch,end,bitstype,ccall,do,using,module,
import,export,importall,baremodule,immutable,local,global,const,Bool,Int,Int8,Int16,Int32,
Int64,Uint,Uint8,Uint16,Uint32,Uint64,Float32,Float64,Complex64,Complex128,Any,Nothing,None,
function,type,typealias,abstract,get_node,add_edge,create_estimation_model,set_solution,
solve, get_solution, solve_ss_problem, create_estimation_problem, addnode
},
morekeywords = [2]{},
sensitive=true,
morecomment=[l]{\#},
morestring=[b]',
morestring=[b]"
}

\begin{document}
\raggedbottom

\title{A Graph-Based Modeling Abstraction for Optimization: \\ Concepts and Implementation in Plasmo.jl}
\titlerunning{Graph-Based Modeling for Optimization}

\date{}
\author{Jordan Jalving \and Sungho Shin \and Victor M. Zavala}
\authorrunning{Jordan Jalving, Sungho Shin, Victor M. Zavala}

\institute{\large \Letter \\
\normalsize
Jordan Jalving \\
jalving@wisc.edu\\
\\
Sungho Shin \\
sungho.shin@wisc.edu\\
\\
Victor M. Zavala \\
victor.zavala@wisc.edu\\
\at  Department of Chemical and Biological Engineering \\ University of Wisconsin-Madison, 1415 Engineering Dr, Madison, WI 53706, USA}

\maketitle

\abstract {We present a general graph-based modeling abstraction for optimization that we call an {\em OptiGraph}. Under this abstraction,
{\em any} optimization problem is treated as a hierarchical hypergraph in which nodes represent optimization subproblems and edges represent connectivity between such subproblems.
The abstraction enables the modular construction of highly complex models in an intuitive manner, facilitates the use of graph analysis tools (to perform partitioning, aggregation,
and visualization tasks), and facilitates communication of structures to decomposition algorithms. We provide an open-source implementation of the abstraction in the {\tt Julia}-based
package {\tt Plasmo.jl}.  We provide tutorial examples and large application case studies to illustrate the capabilities.
}
\\


\section{Motivation and Background}\label{sec:introduction}
Advances in decomposition {\em algorithms} and {\em computing architectures} are continuously expanding the complexity and scale of optimization models that
can be tackled \cite{Kang2014}. Application examples include
financial planning \cite{Gondzio2006}, supply chain scheduling \cite{Maravelias2012}, enterprise-wide management \cite{Grossmann2012},
infrastructure optimization \cite{kang2015}, and network control \cite{rawlings2009}.  Decomposition algorithms seek to address computational bottlenecks
(in terms of memory, robustness, and speed) by exploiting structures present in a model  \cite{Scattolini2009,Zavala2016,Brunaud2017,Conejo2006,Cao2019}.
Well-known algorithms to exploit structures at the problem level include Lagrangian decomposition \cite{Fisher1985},
Benders decomposition \cite{Sahinidis1991}, Dantzig-Wolfe decomposition \cite{Bergner2015}, progressive hedging \cite{Watson2012}, Gauss-Seidel \cite{shin2018multi},
and the alternating direction method of multipliers \cite{Boyd2011}. Fundamentally, these algorithms seek to solve the original problem by solving subproblems defined
over {\em subproblem partitions} and by coordinating subproblem solutions via communication of primal-dual information. Well-known paradigms to exploit structures at the
linear-algebra level (inside optimization solvers) include block elimination and preconditioning (e.g., Schur and Riccati) \cite{Gondzio2003,Kang2014,Rodriguez2020,rao1998application}.
In these schemes, the original problem is solved by using a general algorithm (e.g., interior-point, sequential quadratic programming, or augmented Lagrangian) and
decomposition occurs during the computation of the search step.

The efficient use of decomposition algorithms (and of computing architectures under which they are executed) relies on the ability to communicate model structures in a {\em flexible} manner.
Surprisingly, the development of {\em modeling environments} that support the development of decomposition algorithms has remained rather limited. As a result,
decomposition algorithms have been mostly used to tackle models that have rather obvious structures such as stochastic optimization
\cite{Linderoth2003,Kim2016,Steinbach2003,Lubin2012,Hubner2020},  network optimization  \cite{Shin2019}, dynamic optimization \cite{biel2019efficient}, and hierarchical
optimization \cite{grossman2013}.  Some examples of environments that support structured modeling include the structure-conveying modeling
language ({\em SML}) \cite{Colombo2009}, which is an extension of {\tt AMPL} \cite{ampl} that conveys structures in the form of
variable and constraint blocks.  SML was one of the earliest attempts to automate structure identification and exploitation in modeling environments and was motivated by the availability of structure-exploiting interior point solvers.
{\tt Pyomo} \cite{hart2017pyomo} is a {\em Python}-based modeling package that enables the expression of structures in the form of variable and constraint blocks.
{\tt Pyomo} also provides a structured modeling template called {\tt PySP} \cite{Watson2012}  that greatly facilitates the expression of multi-stage stochastic programs and their solution using progressive hedging and
Benders decomposition. {\tt StuctJuMP} \cite{Huchette2014}
and {\tt StochasticPrograms.jl} \cite{biel2019efficient} are extensions of the {\tt Julia}-based package
{\tt JuMP.jl} \cite{DunningHuchetteLubin2017} to express stochastic optimization structures.

\subsection{Graph-Based Model Representations}
It is important to recall that, at a fundamental level, an optimization model is a collection of algebraic functions (constraints and objectives) that are connected via variables.
Connectivity between functions and variables can be represented as {\em graphs}.
This concept is by no means new; graph representations of optimization models are routinely used in modeling environments to perform automatic differentiation and model
processing tasks \cite{ampl}.  The underlying graph structure of the model is also implicitly communicated to solvers in the form of sparsity patterns of constraint, objective,
and derivative matrices. However, these graph representations operate at a level of granularity that might not be particularly useful for decomposition.
Specifically, the benefits of decomposition tend to become apparent when operating over {\em blocks} of functions and variables.
An issue that arises here is that identifying blocks that are suitable for a decomposition algorithm and computing architecture is not an easy task.
For instance, we might want to identify blocks that have sparse external connectivity (to reduce the amount of inter-block coupling) \cite{Rodriguez2020}.
Moreover, if decomposition is to be executed on a parallel computer, one might want to ensure that the blocks
are of similar size (to avoid load imbalance issues) \cite{kim2019asynchronous}. Another issue that arises here is that blocks might be {\em degenerate}
(e.g., they might have more constraints than variables). Existing modeling environments do not provide capabilities to handle such issues.

\begin{figure}[]
\centering
\begin{subfigure}[t]{0.32\textwidth}
    \centering
    \includegraphics[width=4.9cm, height=8cm]{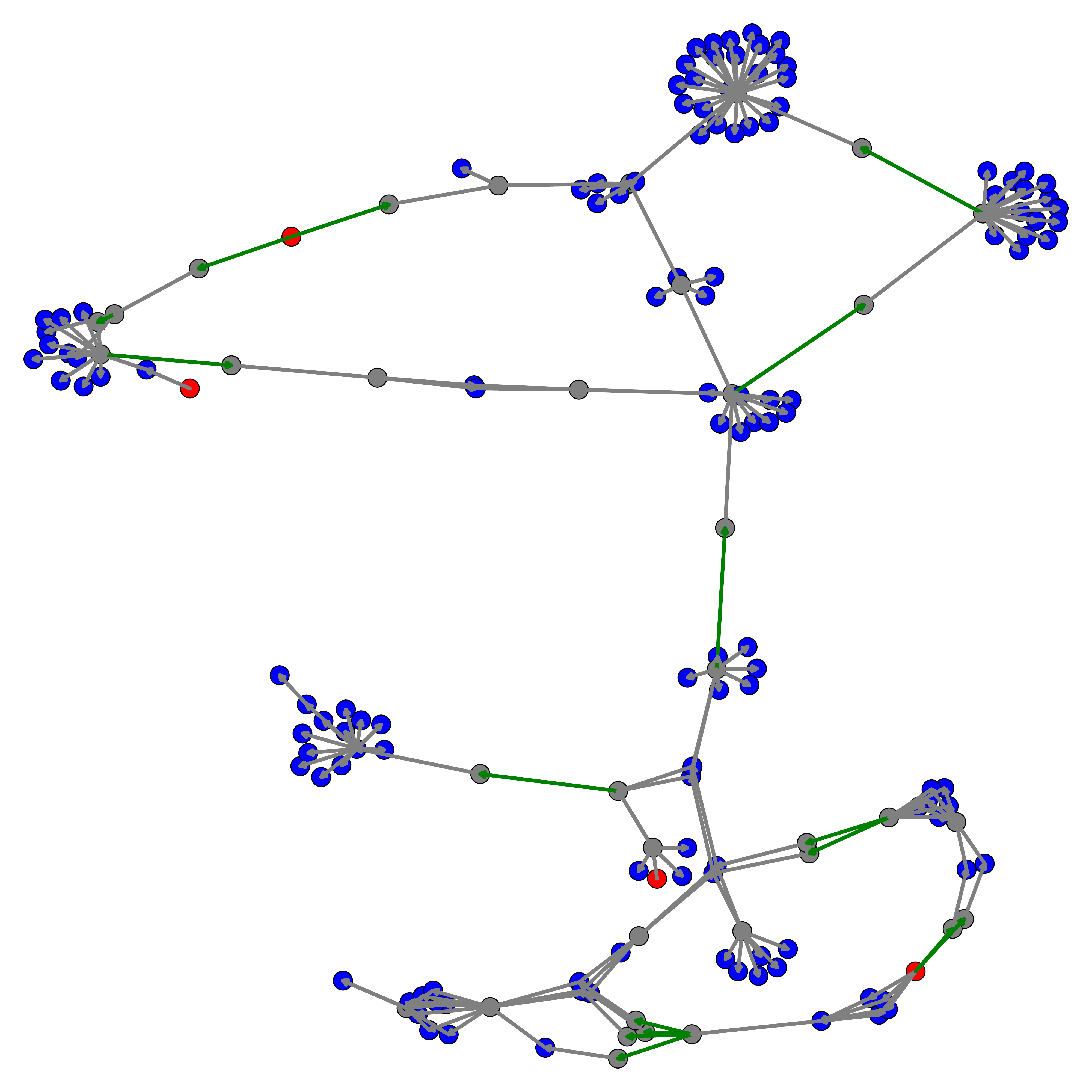}
    \caption{Network Topology}
    \label{fig:gas_network_layout}
\end{subfigure}
\begin{subfigure}[t]{0.32\textwidth}
    \centering
    \includegraphics[width=4.9cm, height=8cm]{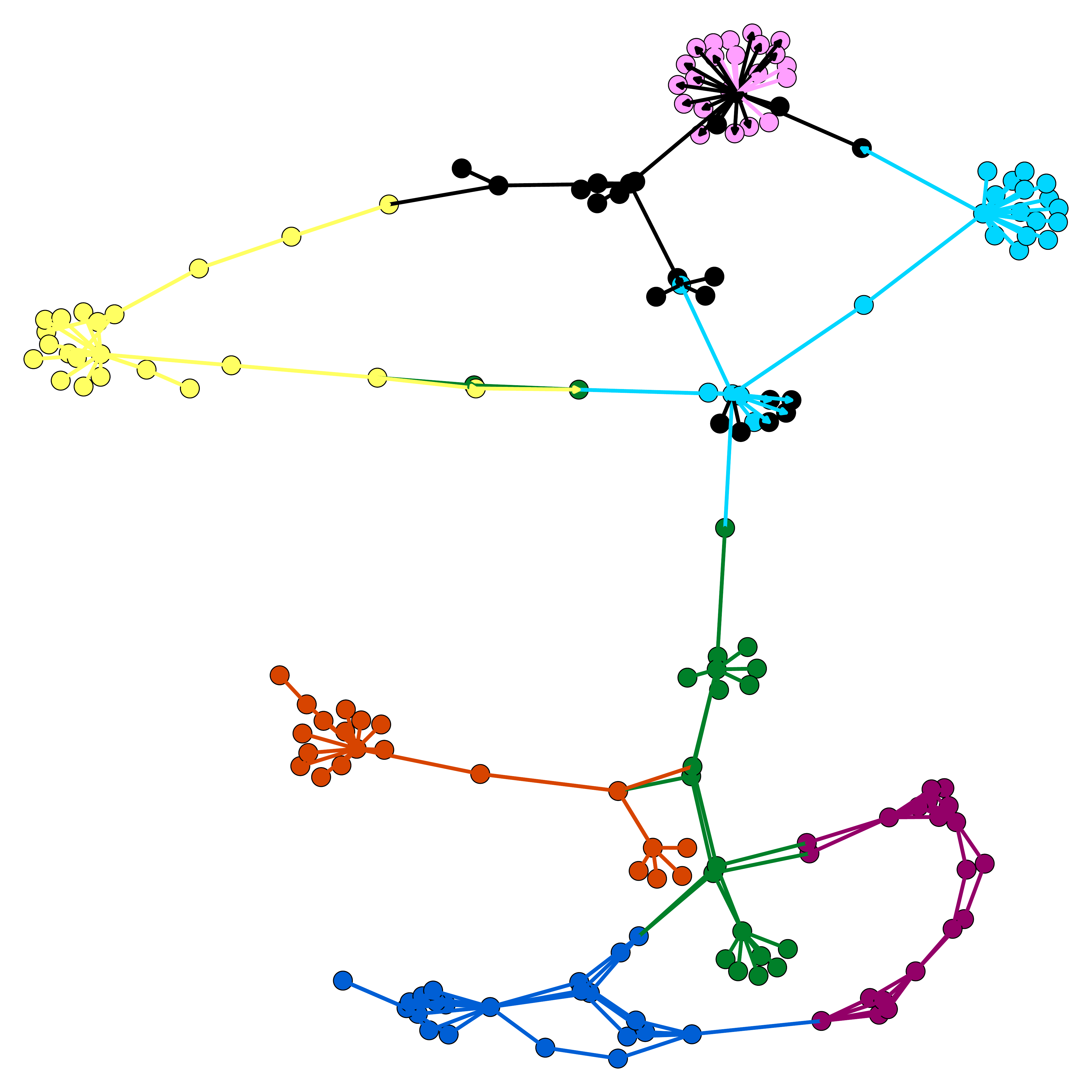}
    \caption{Network Partition}
    \label{fig:gas_network_partition}
\end{subfigure}
\begin{subfigure}[t]{0.32\textwidth}
    \centering
    \includegraphics[width=4.9cm, height=8cm]{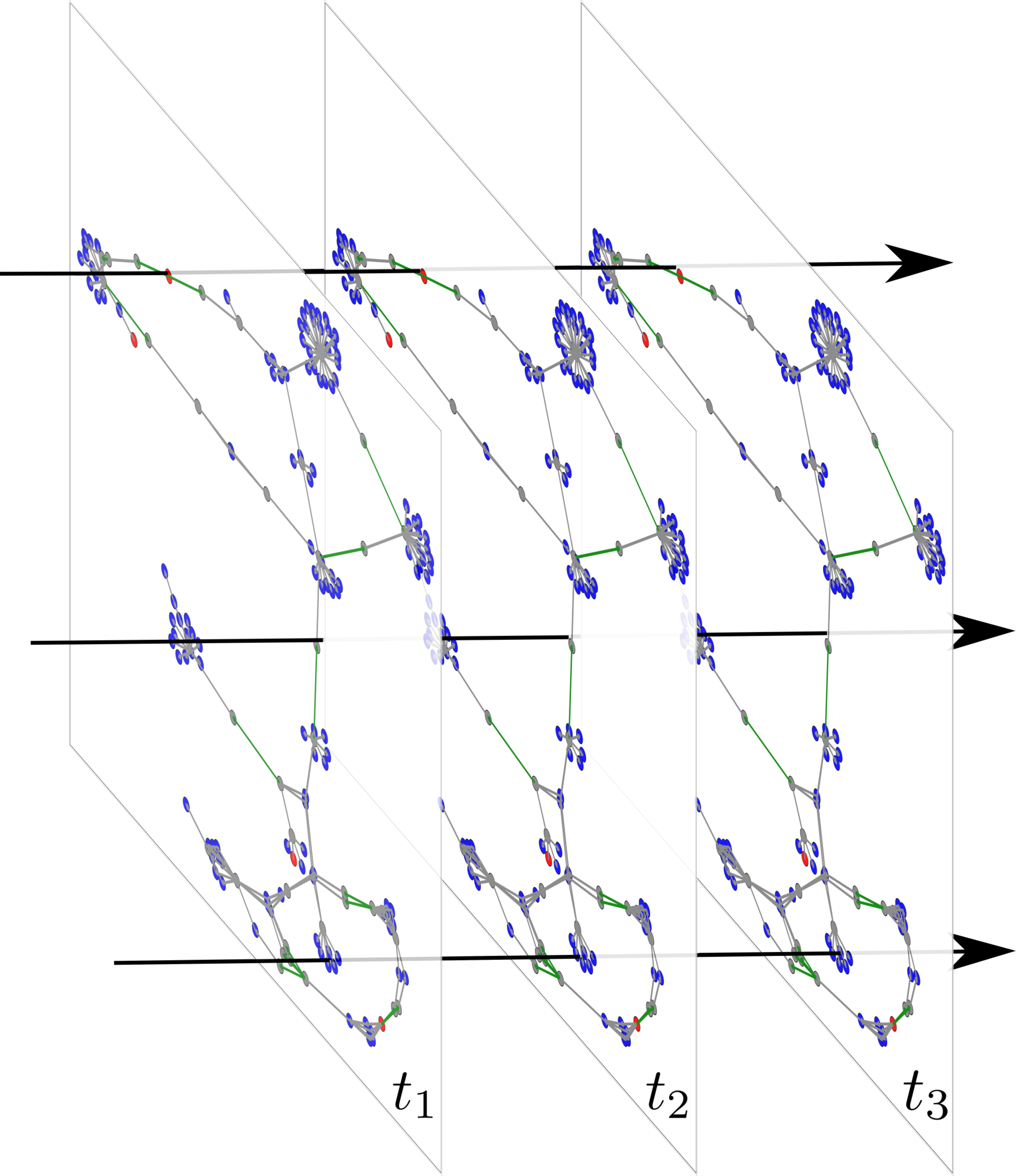}
    \caption{Time Partition}
    \label{fig:gas_time_partition}
\end{subfigure}
\caption{Natural gas transmission network and possible partitioning strategies.
The network layout of the system (left)), the system split into eight network partitions (middle), and the system
represented by three time partitions (right).}
\label{fig:gasnetwork}
\end{figure}

\begin{figure}[]
\centering
\begin{subfigure}[t]{0.48\textwidth}
    \centering
    \includegraphics[width=7.0cm, height=8cm]{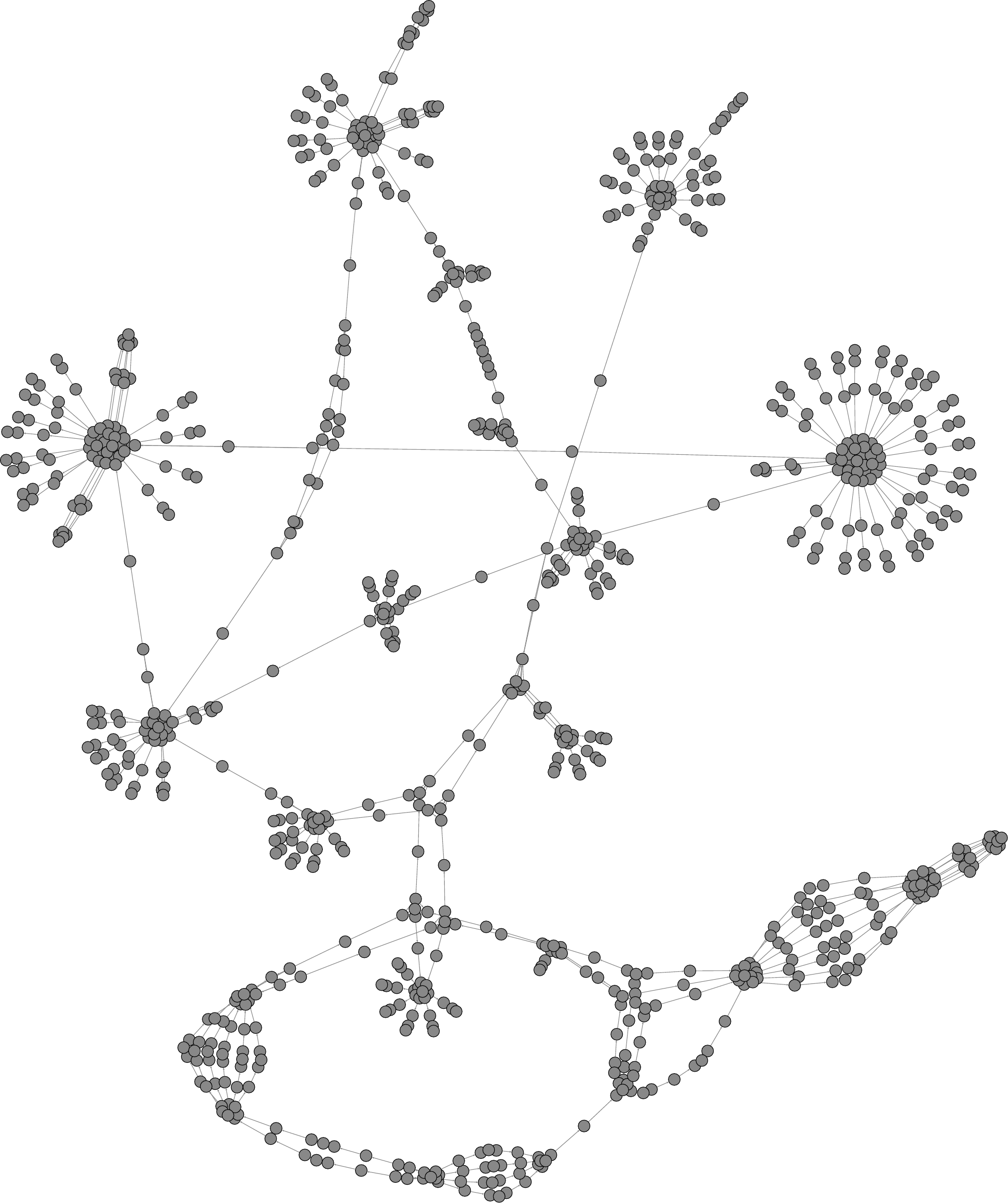}
    \caption{Space Unfolding}
    \label{fig:gas_ocp_space}
\end{subfigure}
\begin{subfigure}[t]{0.48\textwidth}
    \centering
    \includegraphics[width=7.0cm, height=8cm]{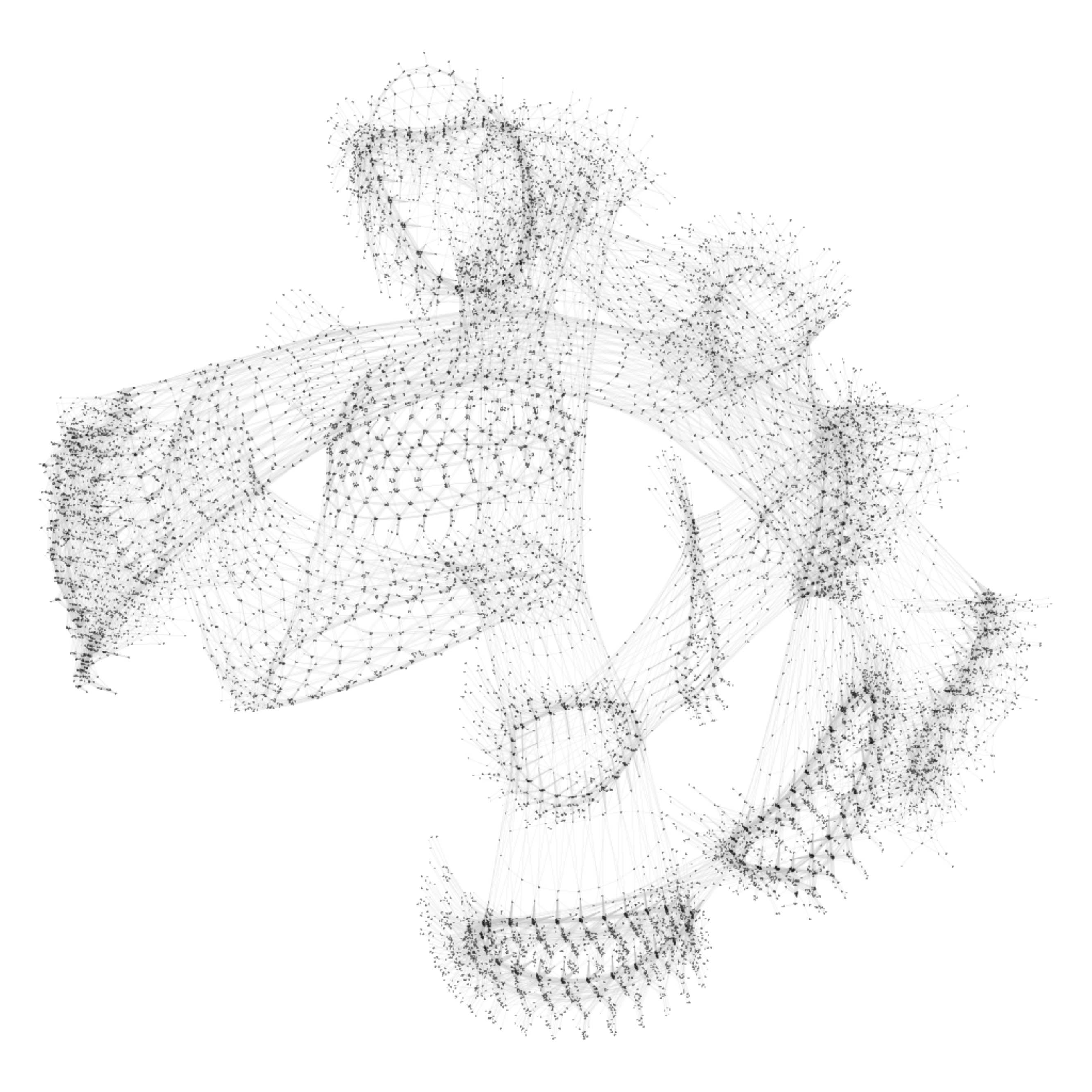}
    \caption{Space-Time Unfolding}
    \label{fig:gas_ocp_space_time}
\end{subfigure}
\caption{Space  and space-time unfolding of natural gas network (uncovering components).}
\label{fig:gasnetwork}
\end{figure}

To highlight some of the challenges that arise in the modeling structured problems, consider the natural gas network depicted in the left panel of Figure \ref{fig:gas_network_layout}.
This is a regional gas transmission system that contains 215 pipelines, 16 compressors, 172 junction points \cite{Chiang2016}.
The problem has an obvious structure induced by the network topology. This representation naturally conveys connectivity between components (pipelines, compressors, and junctions).
Hidden in this representation, however, is the internal connectivity present in such components; this inner connectivity is usually complex and induced by variables and constraints that
capture physical coupling (e.g., conservation of mass, energy, and momentum). The issue here is that there is a severe imbalance in the complexity of the components
(e.g, pipelines might contain partial differential equations to capture flow while compressors contain much simpler algebraic equations). Moreover, hidden in the representation
is the connectivity of the components across {\em time}, which gives rise to complex space-time coupling.  Intuitively, one could seek to decompose the model {\em spatially} by exploiting
the network topology, as shown in the middle panel of Figure \ref{fig:gas_network_partition}. Unfortunately, this approach does not account for inner component complexity and results in
partitions with a disparate number of variables and constraints. This will give rise to load balancing issues and limit the benefits of parallel computation. To deal with this issue,
one could seek to decompose the problem \emph{temporally}, as shown in the right panel of Figure \ref{fig:gas_time_partition}.
This approach leads to well-balanced partitions (every partition is a copy of the network), but leads to significantly more coupling between partitions.
A high degree of coupling will also limit the scalability of algorithms (increased coupling tends to increase communication and slow down convergence).
One would intuitively expect that there exist partitions that can balance coupling and load balancing (by traversing space-time and exploiting inner block complexity).
Identifying such partitions, however, requires {\em unfolding} of the model structure.  Figure \ref{fig:gas_ocp_space} provides a visual rendering of the spatial representation of
the natural gas system that unfolds inner component connectivity. Here, the size of the clusters gives us an initial indication of component complexity.
Figure \ref{fig:gas_ocp_space_time} is a visual rendering of the space-time representation that further unfolds temporal connectivity between components and this clearly
reveals interesting (but complex) structures. Determining effective partitions from such structures requires advanced graph partitioning techniques.

So the question is: {\em If an optimization model is a graph after all, why is it that we do not use graph concepts to build optimization models?} Most optimization modelers
follow an algebraic modeling paradigm in which the model is assembled by adding functions and variables
(e.g., using packages such as {\tt GAMS}, {\tt AMPL}, {\tt Pyomo}).  There is, however, another modeling paradigm known as object-oriented modeling,
in which the model is assembled by adding blocks of functions and variables. This paradigm is widely used in  engineering communities to perform simulation
(e.g., using packages such as {\tt Modelica}, {\tt AspenPlus}, and {\tt gProms} \cite{mattsson1998physical,dowling2015framework}). An important observation is that
object-oriented modeling more naturally lends itself to the expression and exploitation of problem structures because the underlying graph is progressively built by the user.
For instance, packages such as {\tt AspenPlus} use the underlying graph structure induced by blocks to partition the model and communicate it to a decomposition algorithm.
One could thus think of object-oriented modeling as a form of graph-based modeling. The purely object-oriented paradigm however, may suffer from
imblance of component complexity or from a high degree of coupling when performing decomposition.

\begin{figure}[]
\centering
\begin{subfigure}[t]{0.32\textwidth}
    \centering
    \includegraphics[width=4.8cm, height=8cm]{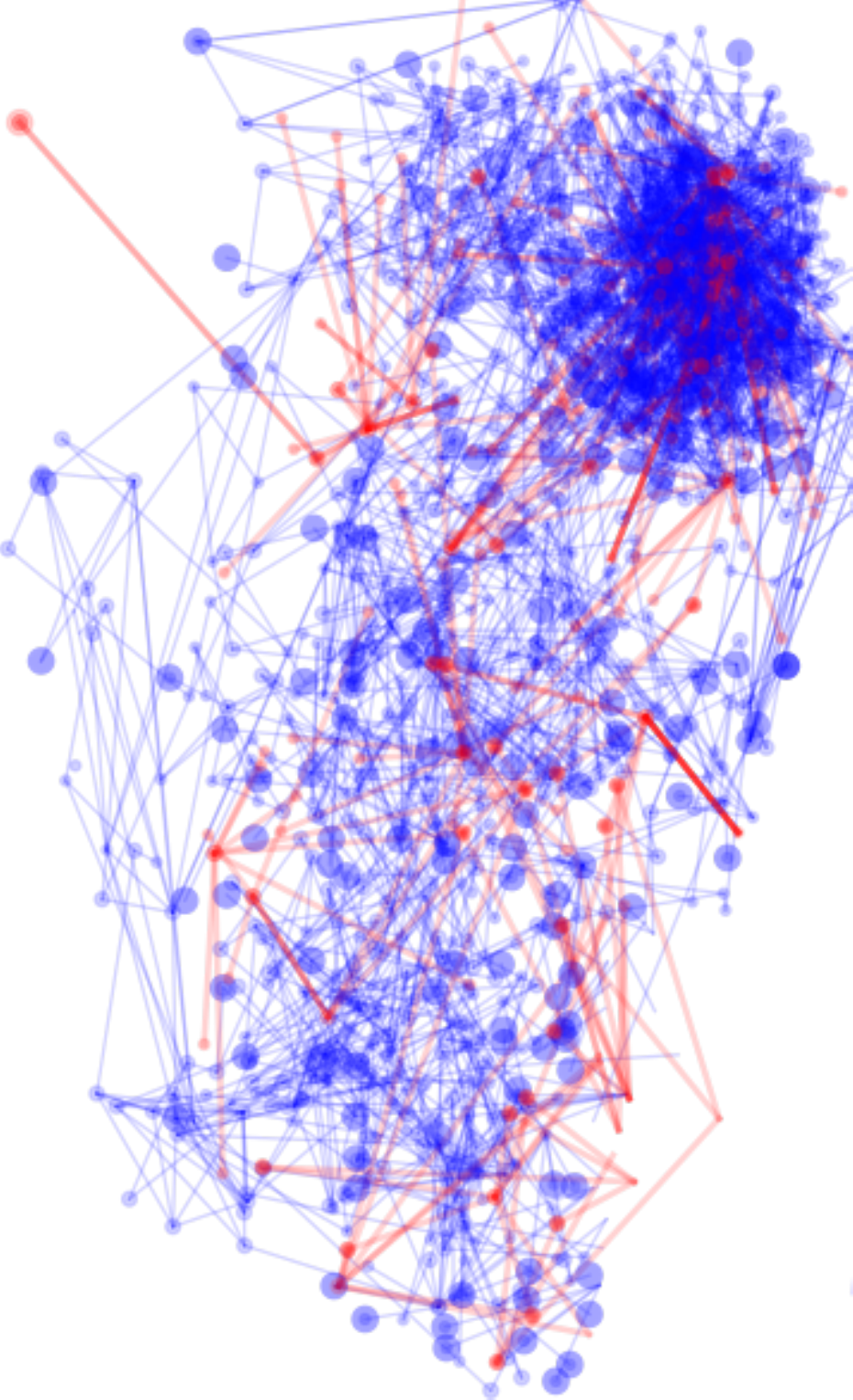}
    \caption{Coupled Networks}
    \label{fig:coupled_gas_electric}
\end{subfigure}
\begin{subfigure}[t]{0.32\textwidth}
    \centering
    \includegraphics[width=4.8cm, height=8cm]{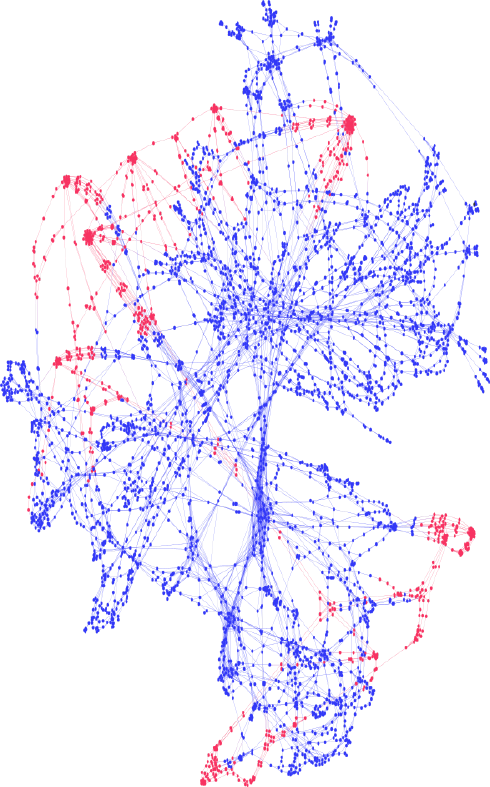}
    \caption{Space Unfolding}
    \label{fig:coupled_space_optimization_structure}
\end{subfigure}
\begin{subfigure}[t]{0.32\textwidth}
    \centering
    \includegraphics[width=4.8cm, height=8cm]{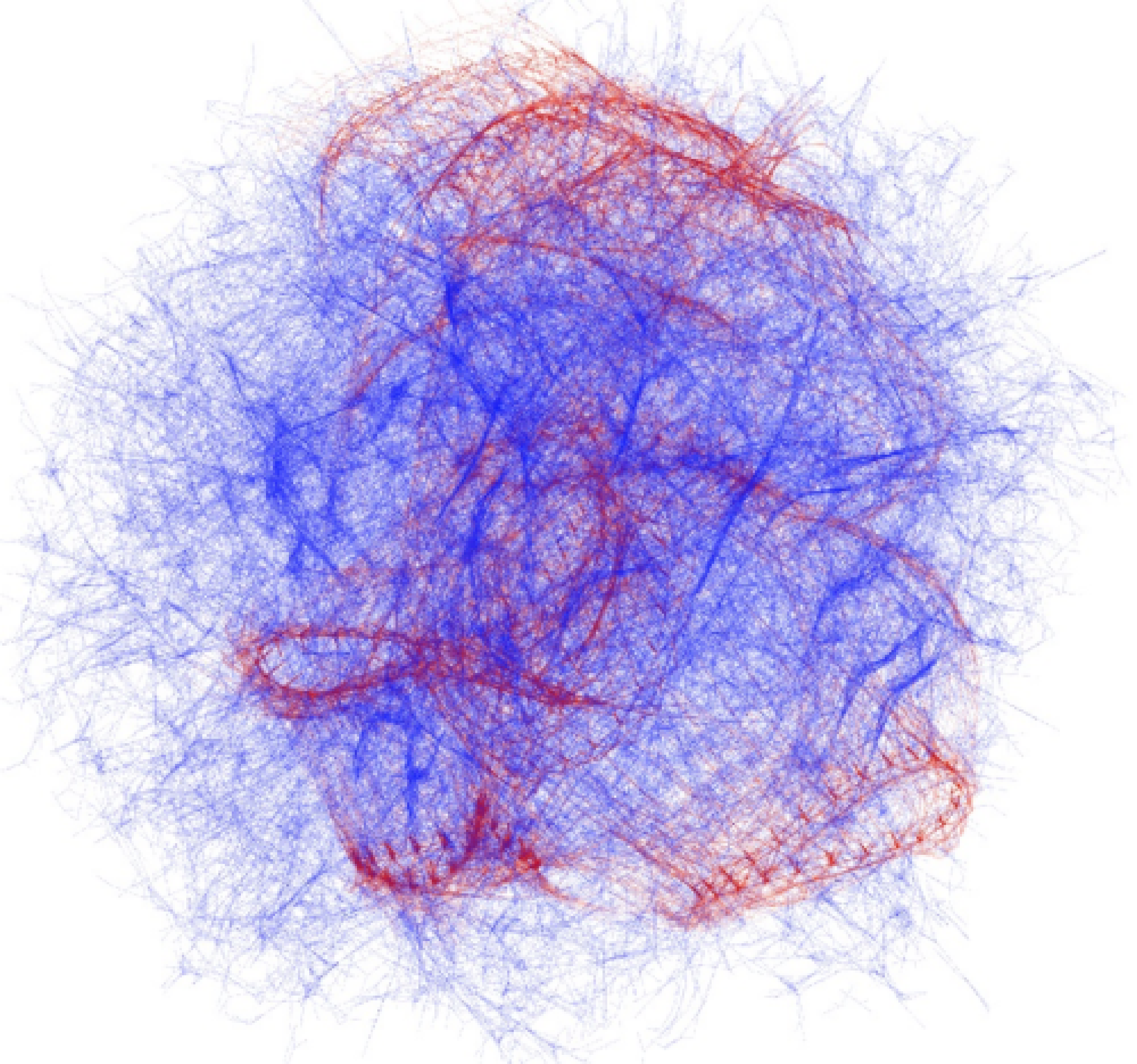}
    \caption{Space-Time Unfolding}
    \label{fig:coupled_space_time_optimization_structure}
\end{subfigure}
\caption{Network topology of coupled natural gas (in blue) and electric (in red) networks (left), space unfolding of components (middle),
and space-time unfolding of components (right).}
\label{fig:coupled_system}
\end{figure}

Recently, we proposed a graph-based abstraction to build optimization models \cite{Jalving2017,Jalving2019}. In this abstraction,
any optimization model is treated as a {\em hierarchical graph}. At a given level, a graph comprises a set of nodes and edges;
each node contains an optimization model (with variables, objectives,  constraints, and data) and each edge captures connectivity between node models.
Importantly, the optimization model in each node can be expressed algebraically (as in a standard algebraic modeling language) or as a graph
(thus enabling inheritance and creation of hierarchical graphs). This provides flexibility to capture both algebraic and object-oriented modeling paradigms under the same framework.
The abstraction naturally exposes problem structure to algorithms and provides a modular approach to construct models.
Modularization enables collaborative model building, independent processing of model components (e.g., automatic differentiation), and data management.
The abstraction naturally captures a wide range of problem classes such as stochastic optimization (graph is a tree), dynamic optimization (graph is a line),
network (graph is the network itself), PDE optimization (graph is a mesh), and multiscale optimization (graph is a meshed tree).
Moreover, graphs can be naturally built by combining graphs from different classes (e.g., stochastic PDE optimization).
This approach is thus more general than other graph-based abstractions proposed for specific problem classes such as network optimization
and control \cite{Heo2014,Jogwar2015,Moharir2017,Tang2018,Hallac2015}.  This modeling abstraction also generalizes those used in simulation
packages such as {\tt Modelica}, {\tt AspenPlus}, {\tt gProms}, which are tailored to specific physical systems.
The hierarchical graph structure can be communicated to decomposition solvers or it can be collapsed into an optimization model that can be solved with off-the-shelf solvers
(e.g., {\tt Gurobi} or {\tt Ipopt}).  To illustrate the types of capabilities enabled by hierarchical graph modeling, we take the natural-gas network in
Figure \ref{fig:gas_network_layout} and couple it to an electrical network, as shown in Figure \ref{fig:coupled_gas_electric}.  Here, the natural gas graph and the electrical
power graph can be built independently and then coupled to construct a higher level graph. The graph structure can be communicated to a graph visualization tool
and this allows us to analyze the graph using different representations.  In the left panel we see a traditional representation of infrastructure
coupling (that hides temporal and component coupling), while in the middle and right panels we unfold spatial and spatio-temporal coupling. The space-time graph contains
over 100,000 nodes and 300,000 edges and reveals that there exist a large number of non-obvious structures that might be beneficial from a computational standpoint.

Graph-based model representations can take advantage of powerful graph {\em analysis} tools.  For instance, graph partitioning tools such as
{\tt Metis} \cite{Karypis1998} and {\tt Scotch} \cite{pellegrini} provide efficient algorithms to automatically analyze problem structure and to identify suitable
partitions to be exploited by decomposition algorithms. Graph partitioning seeks to decompose a domain described by a graph such that communication between subdomains is minimized subject to
load balancing (i.e., the domains are about the same size).  The most ubiquitous graph partitioning applications target
parallelizing scientific simulations \cite{schloegel2003}, performing sparse-matrix operations \cite{Gupta1997}, and preconditioning systems of PDEs \cite{Schulz15coursenotes}.
Hypergraph partitioning generalizes graph partitioning concepts to effectively
decompose non-symmetric domains \cite{Devine2006} and has been applied to physical network design and sparse-matrix-vector
multiplication \cite{Papa2007}. Popular hypergraph partitioning tools include {\tt hMetis} \cite{KarypishMetis} and {\tt PaToH} \cite{Çatalyürek2011}.
Community detection approaches have also been used to find partitions
that maximize modularity \cite{Newman2006,Allman2018}. Graph partitioning approaches have been used to decompose optimization problems in various contexts.
Graph approaches have been used to partition network optimization problems
using graph coloring techniques \cite{Zenios1988} and
block partitioning schemes \cite{Zenios1992}. Bipartite graph representations
have been used to permute linear programs \cite{Ferris1998} into block diagonal form to enable parallel solution. Hypergraph partitioning has been used to
decompose mixed-integer programs to formulate Dantzig-Wolfe decompositions \cite{Wang2013,Bergner2015}
using {\tt hMetis} and {\tt PaToH}. Community detection approaches have been used to automate structure identification in general optimization problems and with this
enable higher efficiency of decomposition algorithms
\cite{Tang2017,Allman2018}. We highlight that these approaches start with a given algebraic optimization model that is then transformed into a graph
representation; consequently, these are not graph modeling abstractions {\em per se}.

Graph-based representations have been adopted in {\em modeling environments} for specific problem classes. {\tt SnapVX} \cite{Hallac2015} uses a graph topology
abstraction to formulate network optimization problems that are solved using
ADMM \cite{Boyd2011}. {\tt GRAVITY} \cite{Hijazi2018} is an algebraic modeling framework that incorporates graph-aware modeling syntax but this is not the paradigm driving the environment.
{\tt DISROPT}  is a {\tt Python}-based environment for distributed network optimization \cite{Farina}.  {\tt DeCODe} \cite{decode} is a package written in {\tt Matlab} and
{\tt Python} that automatically creates a bipartite representation of optimization models created with {\tt Pyomo} and uses community detection to determine partitions for use
with Benders or Lagrangian decomposition. In recent work we have provided a {\tt Julia}-based implementation of our hierarchical graph abstraction,
that we called {\tt Plasmo.jl} \cite{Jalving2017,Jalving2019}. This package provides modeling syntax to systematically create hierarchical graph models,
to integrate algebraic {\tt JuMP.jl} \cite{DunningHuchetteLubin2017} models within nodes, and to facilitate the expression of connectivity between node models.
The first version of {\tt Plasmo.jl} used an abstraction that was targeted to handle coupled infrastructure networks.
This abstraction sought to generalize the network modeling capabilities of {\tt DMNetwork} \cite{abhyankar2018petsc} (used for simulation) to handle optimization
problems \cite{Jalving2017}. Subsequent development provided data management capabilities
to facilitate model reduction and re-use (for real-time optimization applications) \cite{Jalving2018} and provided an interface
to communicate graph structures to the parallel interior-point solver {\tt PIPS-NLP} \cite{cao2017scalable}. The abstraction was later used to
create a {\tt computational graph} abstraction wherein nodes are computational tasks and edges communicate data between tasks.
This abstraction enables hierarchical modeling of cyber-physical systems, workflows, and algorithms \cite{Jalving2019}.

\subsection{Contributions of This Work}

In this work, we present a new and general graph modeling abstraction for optimization that we call an {\tt OptiGraph}. As in our previously proposed abstraction,
we represent any optimization model as a hierarchical graph wherein nodes contain optimization models with corresponding objectives, variables, constraints, and data.
The novel feature of our abstraction lies on how {\em edge connectivity} between models is represented.
Specifically, we express connectivity in the form of {\em hyperedges} (edges that can connect multiple nodes), thus generalizing the concept of an edge
(which connects two edges). This hypergraph representation enables a more intuitive expression of structures at the modeling level and enables the use of
efficient hypergraph analysis tools. The hypergraph representation can also be transformed into a standard graph representation (via lifting) and to other graph
representations (e.g., bipartite) and this enables the use of different graph analysis techniques and tools and enables flexible expression of structures for decomposition.
We provide an efficient implementation of the {\tt OptiGraph} abstraction in {\tt Plasmo.jl} that can handle hundreds of thousands to millions of nodes and edges.
We implement an interface to {\tt LightGraphs.jl} and {\tt Gephi} to demonstrate that the implementation facilitates interfacing to graph analysis and
visualization tools \cite{Bromberger17,gephi}. We illustrate how to use these tools and topological information obtained from the {\tt OptiGraph} to automatically
obtain partitions that trade-off coupling and load balancing and to compute useful descriptors of the graph (e.g., neighborhoods, degree distribution, minimum spanning tree, modularity).
Moreover, we show how to manipulate the graph to perform aggregation functions (merge nodes).  These features provide flexibility to visualize, decompose,
and solve large models using different techniques.  We implement an interface to the parallel solver {\tt PIPS-NLP} to analyze computational efficiency
obtained with different partitioning strategies. We also demonstrate that the abstraction facilitates the implementation of decomposition algorithms;
to do so, we implement an overlapping Schwarz scheme that is tailored to graph-structured problems \cite{Shin2020}.

\subsection{Paper Structure}
The rest of the paper is organized as follows: Section \ref{sec:software_framework} describes the {\tt OptiGraph} abstraction and
its implementation in {\tt Plasmo.jl}.  In Section \ref{sec:algorithm_interfaces} we
show how the {\tt OptiGraph} structure can be exploited by different optimization algorithms. Section \ref{sec:case_study} provides case studies that demonstrate capabilities on a large natural gas network problem and on a direct current optimal power flow (DC OPF) problem. Section \ref{sec:conclusion} concludes the manuscript and discusses future directions for {\tt Plasmo.jl}.

\section{Hierarchical Graph Abstraction and Implementation}\label{sec:software_framework}

This section introduces the {\tt OptiGraph} abstraction alongside its implementation in {\tt Plasmo.jl}.
In a nutshell, {\tt Plasmo.jl} is a graph-structured modeling language that offers concise syntax and access to graph tools that enable model analysis and manipulation
(e.g., partitioning, aggregation, visualization).  {\tt Plasmo.jl} is an
official {\tt Julia} package hosted at \url{https://github.com/zavalab/Plasmo.jl}. All example scripts in this
section can be found at \url{https://github.com/zavalab/JuliaBox/tree/master/PlasmoExamples}.

\subsection{Graph Abstraction}

The {\tt OptiGraph} abstraction proposed is composed of a set of \emph{OptiNodes} $\mathcal{N}$ (each embedding an optimization model with its local
variables, constraints, objective function, and data) and a set of  \emph{OptiEdges} $\mathcal{E}$ (each embedding a set of \emph{linking  constraints})
that capture coupling between {\em OptiNodes}. The {\tt OptiGraph} is denoted as $\mathcal{G}(\mathcal{N},\mathcal{E})$ and contains the optimization model of interest. The {\tt OptiEdges} $\mathcal{E}$ are hyperedges that connect two or more {\tt OptiNodes}. The {\tt OptiGraph} is an undirected hypergraph.  Whenever clear from context, we simply refer to the {\tt OptiGraph} as graph, to the {\tt OptiNodes} as nodes, and to {\em OptiEdges} as edges. We denote the set of nodes that belong to $\mathcal{G}$ as $\mathcal{N}(\mathcal{G})$ and the set of edges as
$\mathcal{E}(\mathcal{G})$. The topology of the hypergraph is encoded in the
incidence matrix $A \in \mathcal{R}^{|\mathcal{N}|\times |\mathcal{E}|}$; here, the notation $|\mathcal{S}|$ denotes cardinality of set $\mathcal{S}$. The neighborhood of node $n$ is denoted as $\mathcal{N}(n)$ and this is the set of nodes connected to $n$.  The set of nodes that support an edge $e$ are denoted as $\mathcal{N}(e)$ and the set of edges that are incident to node $n$ are denoted as $\mathcal{E}(n)$.

The optimization model associated with an {\tt OptiGraph} can be represented mathematically as:
\begin{subequations}\label{eq:model-graph-compact}
    \begin{align}
        \min_{{\{x_n}\}_{n \in \mathcal{N}(\mathcal{G})}} & \quad \sum_{n \in \mathcal{N(\mathcal{G})}} f_n(x_n) \quad & (\textrm{Objective})\label{eq:graph-objective} \\
        \textrm{s.t.} & \quad x_n \in \mathcal{X}_n,      \quad n \in \mathcal{N(\mathcal{G})}, \quad & (\textrm{Node Constraints})\label{eq:node-constraints} \\
        & \quad g_e(\{x_n\}_{n \in \mathcal{N}(e)}) = 0,  \quad e \in \mathcal{E(\mathcal{G})}. &(\textrm{Link Constraints})  \label{eq:link-constraints}
    \end{align}
\end{subequations}
Here, ${\{x_n}\}_{n \in \mathcal{N}(\mathcal{G})}$ is the collection of decision variables over the entire set of nodes $\mathcal{N}(\mathcal{G})$ and $x_n$ is the set of
variables on node $n$. Function
\eqref{eq:graph-objective} is the graph objective function which is given by the sum of objective functions $f_n(x_n)$,
\eqref{eq:node-constraints} represents the collection of constraints over all nodes $\mathcal{N}(\mathcal{G})$, and \eqref{eq:link-constraints} is the collection of
linking constraints over all edges $\mathcal{E}(\mathcal{G})$. Here, the constraints of a node $n$ are represented by the set $\mathcal{X}_n$ while the linking constraints induced by an
edge $e$ are represented by the vector function $g_e(\{x_n\}_{n \in \mathcal{N}(e)})$ (an edge can contain multiple linking constrains). The optimization problem representation captures
sufficient features to facilitate the discussion but differs from the actual implementation.  For example, here we assume that the graph objective is obtained via a linear combination of the
node objectives but other combinations are possible (e.g., to handle conflict resolution formulations wherein nodes represent different stakeholders). In addition, we assume that coupling
between nodes arises in the form of complicating constraints but definition of complicating variables is also possible (coupling variables can always be represented as coupling
constraints via lifting).

\subsection{Syntax and Usage}

\begin{figure}[]
    \centering
    \includegraphics[scale=0.18]{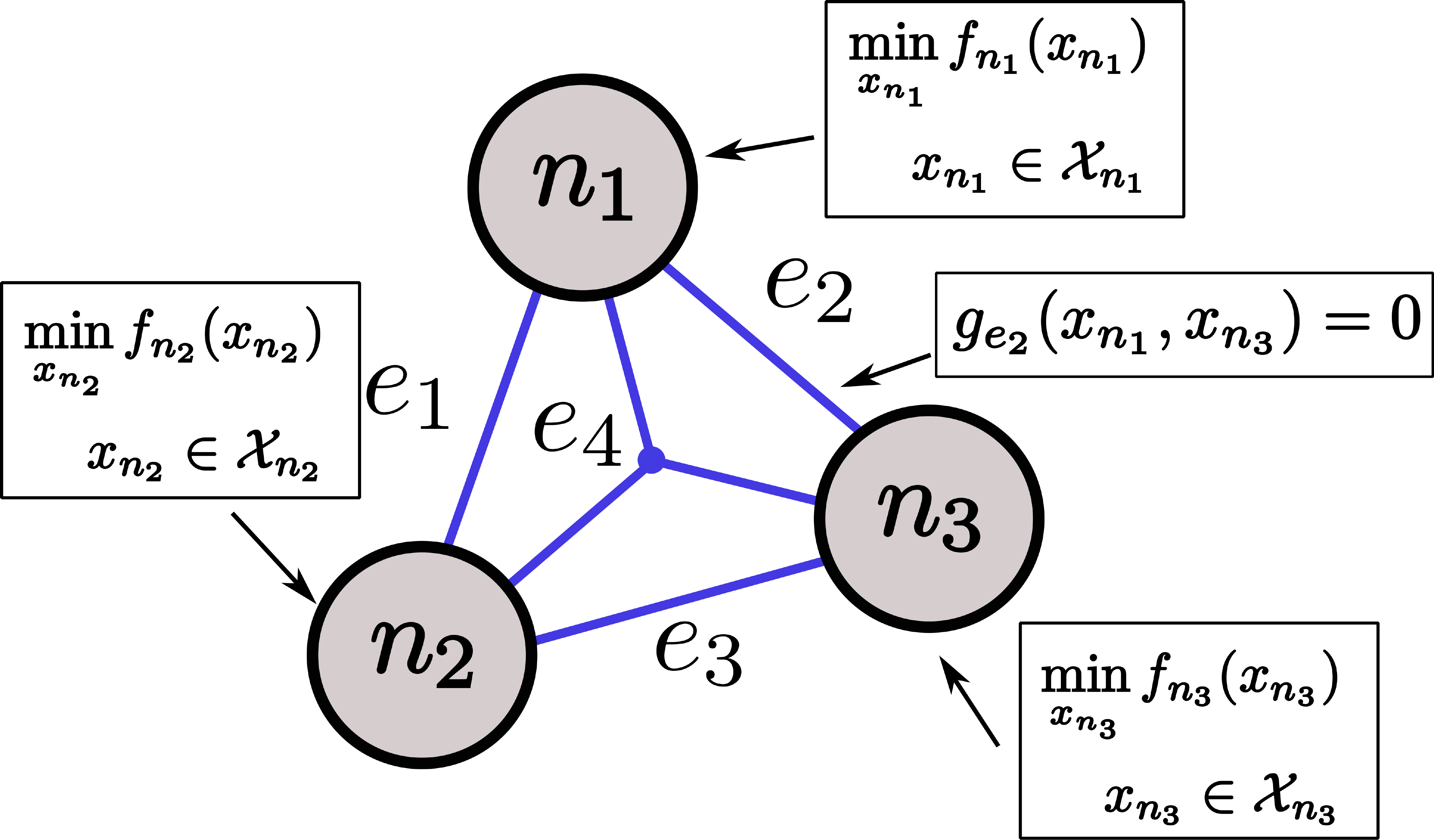}
    \caption{{\tt OptiGraph}  with three {\tt OptiNodes} connected by four {\tt OptiEdges}}.
    \label{fig:OptiGraph_representation}
\end{figure}

We now introduce the {\tt Plasmo.jl} implementation of the {\tt OptiGraph} and show how to use its syntax to create, analyze, and solve an optimization model. {\tt Plasmo.jl} implements the graph model described by \eqref{eq:model-graph-compact} and illustrated in Figure \ref{fig:OptiGraph_representation}.  In our implementation, an {\tt OptiNode} encapsulates a {\tt Model} object from the {\tt JuMP.jl}
modeling language. This harnesses the algebraic modeling syntax and processing functionality of {\tt JuMP.jl}.  An {\tt OptiEdge} object encapsulates the linking constraints that define coupling  between the nodes.

\subsubsection{Example 1 : Creating an OptiGraph}
We start with the simple example given by \eqref{eq:example1} to demonstrate {\tt OptiGraph} syntax:

\begin{subequations} \label{eq:example1}
    \begin{align}
         \min \quad & y_{n_1} + y_{n_2} + y_{n_3} \quad & \text{(Objective)}  \\
        \textrm{s.t.} \quad & x_{n_1} \ge 0, y_{n_1} \ge 2, x_{n_1} + y_{n_1} \ge 3 \quad & \text{(Node 1 Constraints)} \label{eq:example1_node1} \\
        & x_{n_2}  \ge 0,  x_{n_2} + y_{n_2} \ge 3  \quad & \text{(Node 2 Constraints)} \label{eq:example1_node2}\\
        & x_{n_3} \ge 0, x_{n_3} + y_{n_3} \ge 3  \quad & \text{(Node 3 Constraints)} \label{eq:example1_node3}\\
        & x_{n_1} + x_{n_2} + x_{n_3} = 3  \quad & \text{(Link Constraint)} \label{eq:example1_link}
    \end{align}
\end{subequations}
In this model, equations \eqref{eq:example1_node1}, \eqref{eq:example1_node2}, and \eqref{eq:example1_node3} represent
individual node constraints and \eqref{eq:example1_link} is a linking constraint that couples the three nodes. We formulate and solve this optimization
model as shown in Snippet \ref{code:example1_code_snippet}.
We import {\tt Plasmo.jl} into a {\tt Julia} session on Line \ref{line:importmg} as well as the off-the-shelf linear programming solver {\tt GLPK} \cite{glpk} to solve the problem.
We define {\tt graph1} (an {\tt OptiGraph}) on Line \ref{line:OptiGraph} and then create
three {\tt OptiNodes} on Lines \ref{line:model_create_start}-\ref{line:model_create_end} using the {\tt @optinode} macro.  We also
use the {\tt @variable}, {\tt @constraint}, and {\tt @objective} macros (extended from {\tt JuMP.jl}) to define node model attributes. Next, we use the {\tt @linkconstraint}
macro on Line \ref{line:link_constraint} to create a linking constraint between the three nodes. Importantly, this {\em automatically} creates an {\tt OptiEdge}.
This feature is key, as the user does not need to express the topology of the graph (this is automatically created behind the scenes as linking constraints are added).
In other words, the user does not need to provide an adjacency matrix (which can be highly complex). We solve the problem using the {\tt GLPK} optimizer
on Line \ref{line:solve_glpk}.  Since {\tt GLPK} does not exploit graph structure, {\tt Plasmo.jl} automatically transforms the graph into a standard LP format.
We query the solution for each variable using the {\tt value} function on Line \ref{line:query_solution} which
accepts the corresponding node and variable we wish to query. Another important feature is that the solution data retains the structure of the {\tt OptiGraph},
and this facilitates query and analysis. We also highlight that the syntax is similar to that of {\tt JuMP.jl} but operates at a higher level of abstraction.

\begin{figure}
\begin{minipage}{\textwidth}
\begin{minipage}[t]{\textwidth}
\begin{scriptsize}
\lstset{language=Julia,breaklines = true}
\begin{lstlisting}[label = {code:example1_code_snippet},caption = Creating and Solving an {\tt OptiGraph} in Plasmo.jl]
using Plasmo       |\label{line:importmg}|
using GLPK
using Plots

graph1 = OptiGraph()   |\label{line:OptiGraph}|

#Create three OptiNodes
@optinode(graph1,n1)                |\label{line:model_create_start}|
@variable(n1, y >= 2)
@variable(n1,x >= 0)
@constraint(n1,x + y >= 3)
@objective(n1, Min, y)

@optinode(graph1,n2)
@variable(n2, y)
@variable(n2,x >= 0)
@constraint(n2,x + y >= 3)
@objective(n2, Min, y)

@optinode(graph1,n3)
@variable(n3, y )
@variable(n3,x >= 0)
@constraint(n3,x + y >= 3)
@objective(n3, Min, y)          |\label{line:model_create_end}|

#Create link constraint between nodes (automatically creates an optiedge on graph1)
@linkconstraint(graph1, n1[:x] + n2[:x] + n3[:x] == 3)  |\label{line:link_constraint}|

#Optimize with GLPK
optimize!(graph1,GLPK.Optimizer)  |\label{line:solve_glpk}|

#Query Solution
value(n1,n1[:x])            |\label{line:query_solution}|
value(n2,n2[:x])
value(n3,n3[:x])
objective_value(graph1)

#Visualize graph topology
plt_graph1 = plt_graph1 = Plots.plot(graph1,node_labels = true,    |\label{line:plot_graph1}|
markersize = 60,labelsize = 30, linewidth = 4,
layout_options = Dict(:tol => 0.01,:iterations => 2));

#Visualize graph adjacency
plt_matrix1 = Plots.spy(graph1,node_labels = true,markersize = 30);    |\label{line:spy_graph1}|
\end{lstlisting}
\end{scriptsize}
\end{minipage}%
\\
\begin{minipage}{0.5\textwidth}
\includegraphics[width = 7.3cm,keepaspectratio]{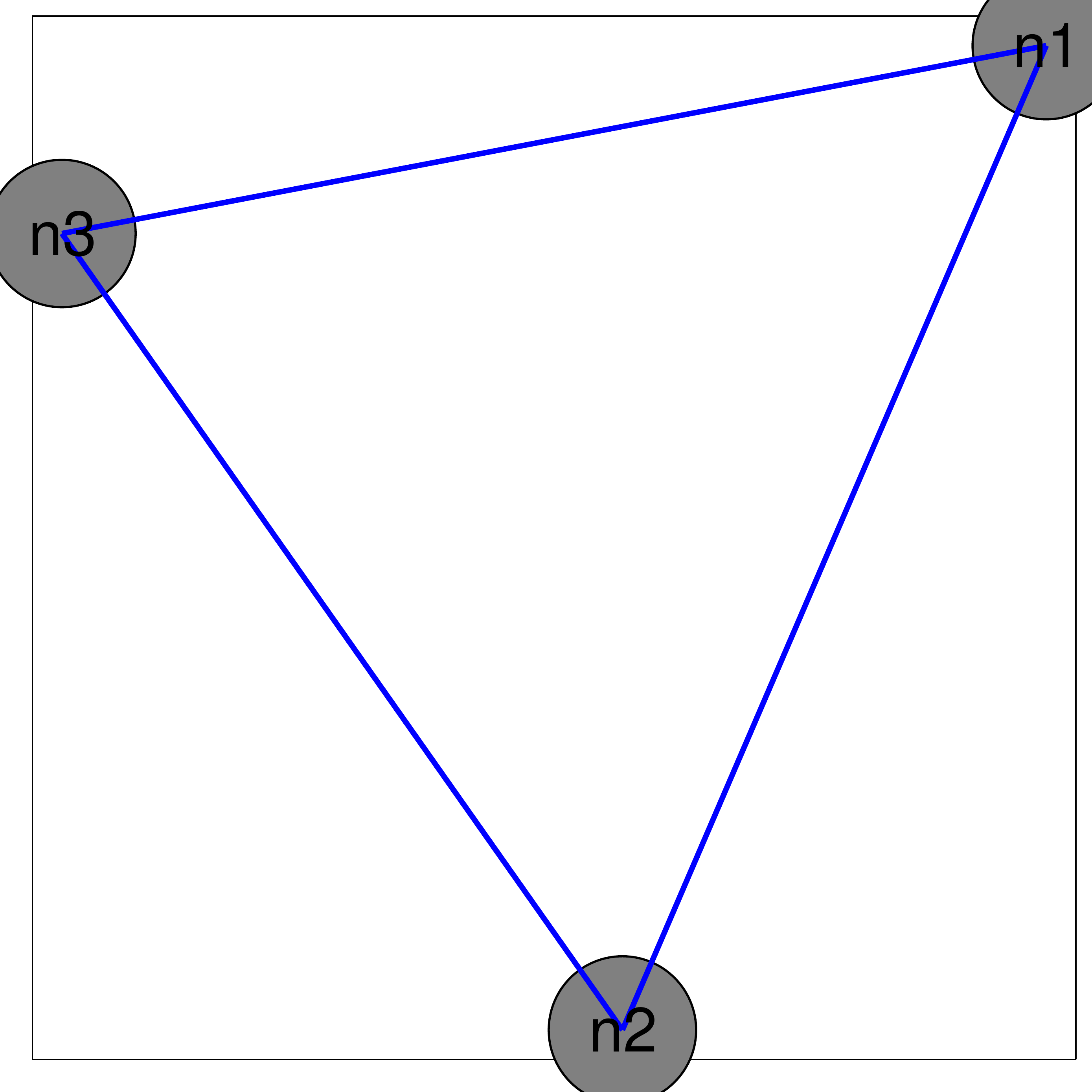}
\end{minipage}%
\begin{minipage}{0.5\textwidth}
\hspace{0.2cm}
\includegraphics[width = 7.3cm,keepaspectratio]{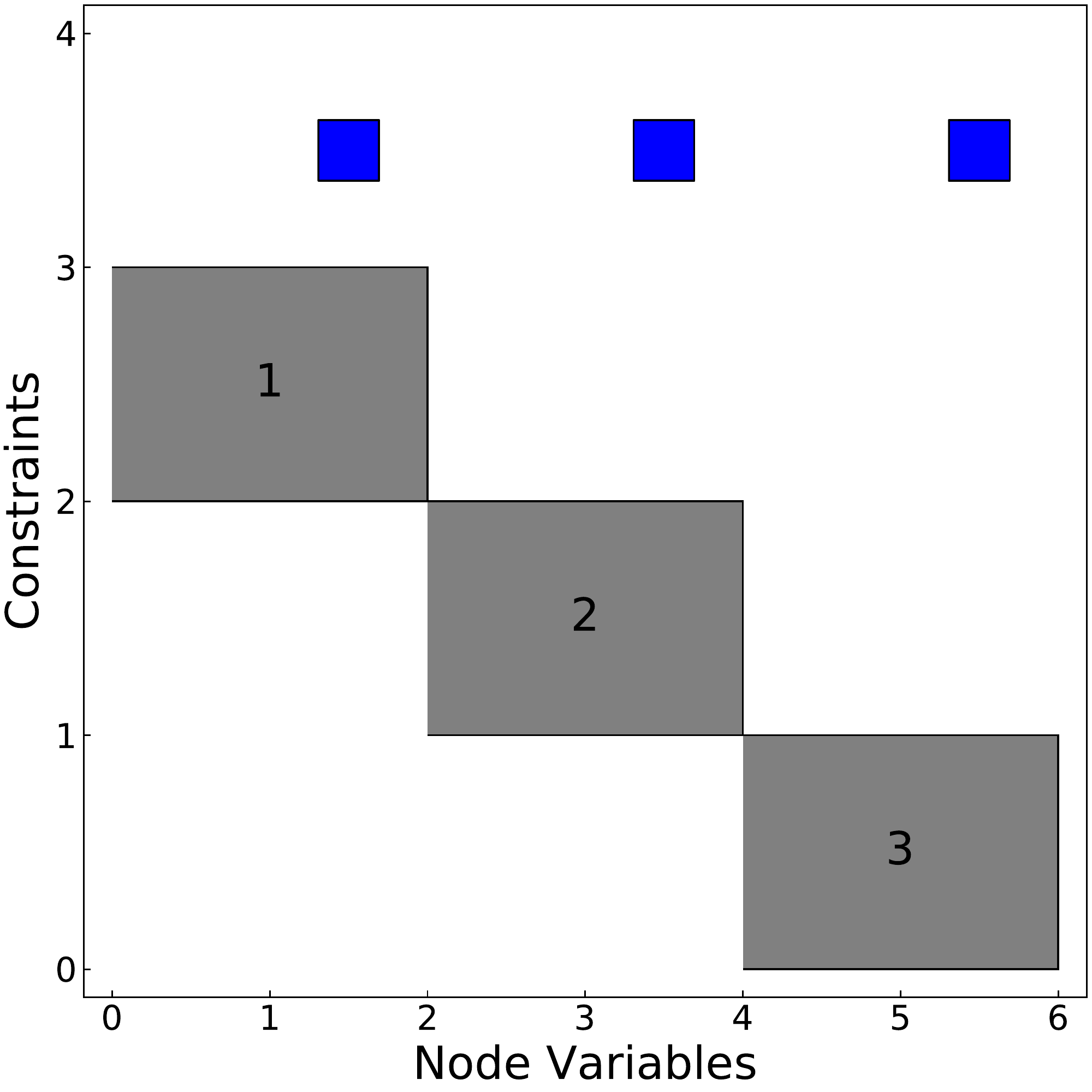}
\end{minipage}
\captionof{figure}{Output visuals for Code Snippet \ref{code:example1_code_snippet}.  Graph topology obtained with {\tt plot} \\ function (left) and graph matrix representation
obtained with {\tt spy} function (right).}
\label{fig:example1_result}
\end{minipage}
\end{figure}

We have implemented capabilities to visualize the structure of {\tt OptiGraphs} by extending plotting functions available in {\tt Plots.jl}.
We layout the {\tt OptiGraph} topology on Line \ref{line:plot_graph1} using the {\tt plot} function and we plot the underlying adjacency matrix
structure on Line \ref{line:spy_graph1} using the {\tt spy} function.  Both of these functions can accept keyword arguments to customize their layout or appearance.
The matrix visualization also encodes information on the number of variables and constraints in each node and edge. The results are depicted in Figure \ref{fig:example1_result};
the left figure shows a standard graph visualization where we draw an edge between each pair of nodes if they share an edge, and the right figure shows the
matrix representation where labeled blocks correspond to nodes and blue marks represent linking constraints that connect their variables.
The node layout helps visualize the overall connectivity of the graph while the matrix layout helps visualize the size of nodes and edges.

\subsection{Hierarchical Graphs}\label{sec:modeling_with_subgraphs}

A key novelty of the {\tt OptiGraph}  abstraction is that it can cleanly represent hierarchical structures. This feature enables expression of models with multiple embedded structures and
enables modular model building (e.g., by merging existing models).
To illustrate this, consider you have the {\tt Optigraphs} $\mathcal{G}_i,\mathcal{G}_j$ each with its own local nodes and edges. We assemble these low-level {\tt OptiGraphs} (subgraphs) to
build a high-level {\tt OptiGraph}
$\mathcal{G}(\{\mathcal{G}_i,\mathcal{G}_j\},\mathcal{N}_g,\mathcal{E}_g)$. We use the notation $\mathcal{G}(\mathcal{N},\mathcal{E})$ to indicate that the high-level
graph has nodes $\mathcal{N} = \mathcal{N}_g \cup \mathcal{N}(\mathcal{G}_i) \cup \mathcal{N}(\mathcal{G}_j)$ and edges
$\mathcal{E} = \mathcal{E}_g \cup \mathcal{E}(\mathcal{G}_i) \cup \mathcal{E}(\mathcal{G}_j)$. Here,
$\mathcal{E}_g$ are global edges in $\mathcal{G}$ (connect nodes across low-level graphs $\mathcal{G}_i,\mathcal{G}_j$ but are not elements of such subgraphs) and $\mathcal{N}_g$ are global
nodes in $\mathcal{G}$ (can be connected to nodes in low-level graphs but that are not elements of such subgraphs).
For every global edge $e\in\mathcal{E}_g$, we have that $\mathcal{N}(e) \in \mathcal{N}(\mathcal{G}_i)$, $\mathcal{N}(\mathcal{G}_j)$, and $\mathcal{N}_g$ where
$\mathcal{N}(\mathcal{G}_i) \cap \mathcal{N}(\mathcal{G}_j) \cap \mathcal{N}_g = \emptyset$.
In other words, the edges $\mathcal{E}_g$ only connect nodes across low-level graphs $\mathcal{G}_i$ and $\mathcal{G}_j$ or between global
nodes $\mathcal{N}_g$ in $\mathcal{G}$ and low-level graphs. Consequently, if we treat the elements of a node set $\mathcal{N}$ as a graph (i.e., we combine  subgraphs into a single {\tt OptiNode}),
we can represent  $\mathcal{G}(\{\mathcal{G}_i,\mathcal{G}_j\}, \mathcal{E}_g,\mathcal{N}_g)$ as $\mathcal{G}(\mathcal{N}, \mathcal{E})$. This nesting procedure can be carried over multiple levels to form a hierarchical graph.
\\

The formulation given by \eqref{eq:model-graph-compact} can be extended to describe a hierarchical graph with an arbitrary number of subgraphs. This is shown in \eqref{eq:model-graph-subgraphs}; here, ${\{x_n}\}_{n \in \mathcal{N}(\mathcal{G})}$
represents the collection of variables over the entire set of nodes in the graph. Equation \eqref{eq:subgraph_links} are the linking constraints over the edges for each subgraph. Here, we use the notation
$\mathcal{SG} \in AG(\mathcal{G})$ to indicate that the subgraph elements $\mathcal{SG}$ are part of the recursive set of subgraphs in $\mathcal{G}$ (i.e., set contains all subgraphs in the graph). Specifically, we define the mapping $AG : \mathcal{G} \rightarrow \{\mathcal{SG}_1,\mathcal{SG}_2,\mathcal{SG}_3,...,\mathcal{SG}_N$\} where
$N$ is the total number of subgraphs in $\mathcal{G}$.

\begin{subequations}\label{eq:model-graph-subgraphs}
    \begin{align}
        &\min_{{\{x_n}\}_{n \in \mathcal{N}}} \quad \sum_{\mathcal{SG} \in AG(\mathcal{G})}\sum_{n \in \mathcal{N}(\mathcal{SG})} f_n(x_n) \quad & (\textrm{Objective})\\
        \textrm{s.t.} & \quad x_n \in \mathcal{X}_n, \qquad\quad n \in \mathcal{N}(\mathcal{SG}), \mathcal{SG} \in AG(\mathcal{G}) \quad & (\textrm{Node Constraints})\\
        & \quad g_e(\{x_n\}_{n \in \mathcal{N}(e)}), e \in \mathcal{E}(\mathcal{SG}), \mathcal{SG} \in AG(\mathcal{G}) \quad & (\textrm{Subgraph Link Constraints}) \label{eq:subgraph_links}\\
        & \quad g_e(\{x_n\}_{n \in \mathcal{N}(e)}), e \in \mathcal{E}(\mathcal{G}). \quad & (\textrm{Graph Link Constraints} \label{eq:parent_graph_links})
    \end{align}
\end{subequations}

\begin{figure}[!htp]
    \centering
    \includegraphics[scale=0.14]{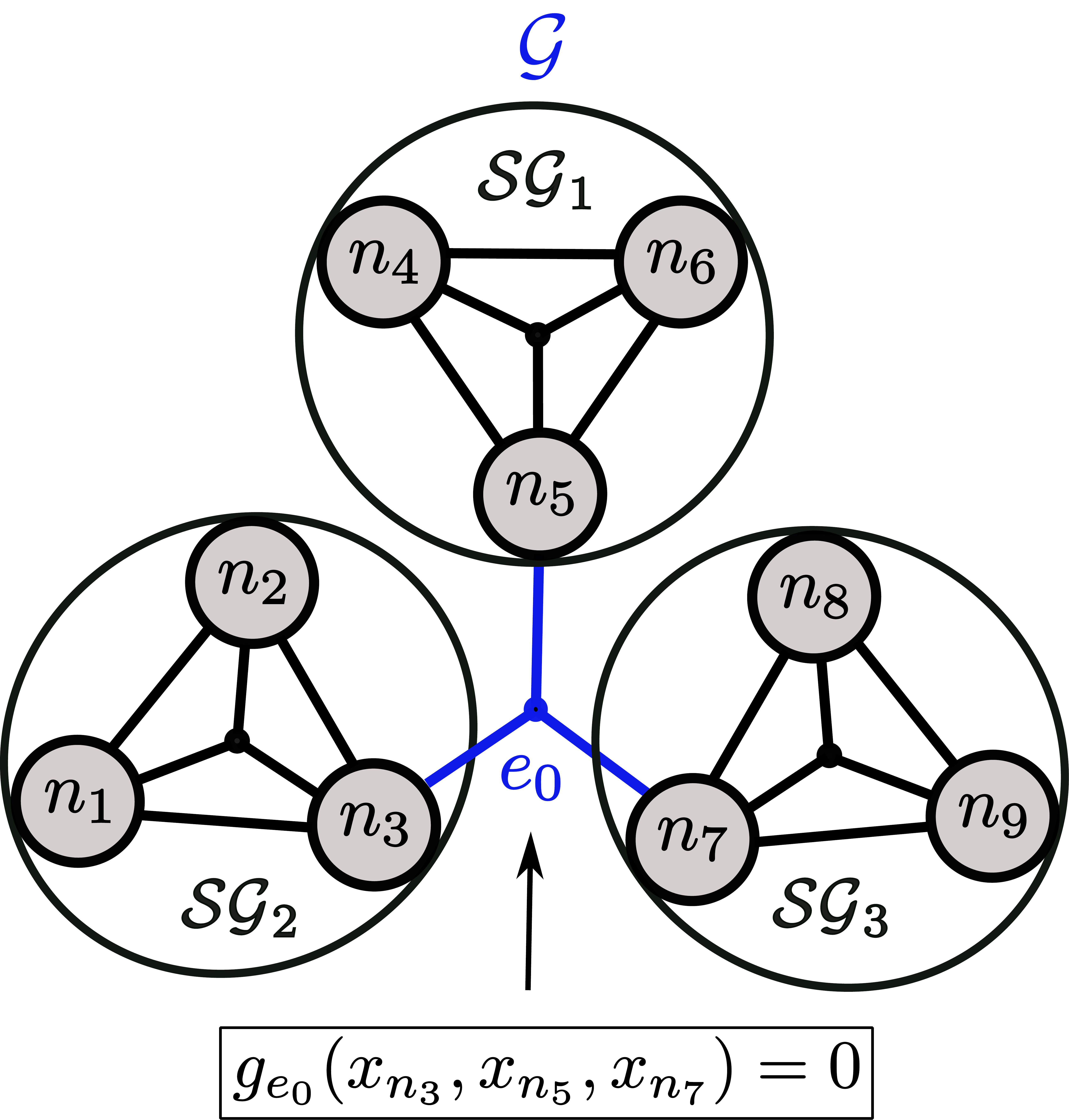}
    \caption{An {\tt OptiGraph} that contains three subgraphs.  The subgraphs are coupled through the global edge $e_0$ in {\tt OptiGraph} $\mathcal{G}$.}
    \label{fig:OptiGraph_subgraphs}
\end{figure}

\begin{figure}[!htp]
    \centering
    \includegraphics[scale=0.14]{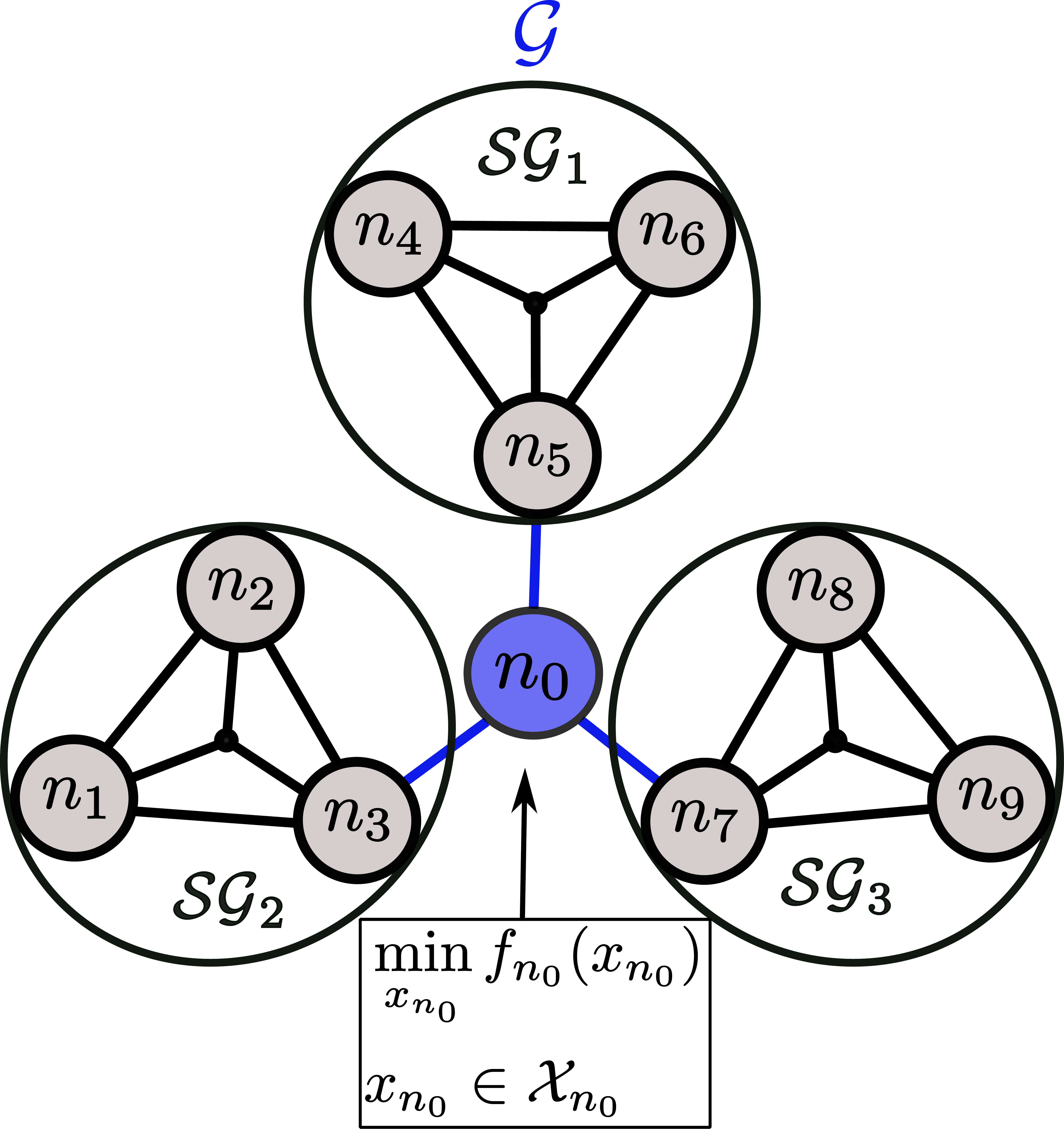}
    \caption{An {\tt OptiGraph} that contains three subgraphs.  The subgraphs are coupled to the global node $n_0$ in {\tt OptiGraph} $\mathcal{G}$.}
    \label{fig:OptiGraph_subgraphs_master}
\end{figure}

Figures \ref{fig:OptiGraph_subgraphs} and \ref{fig:OptiGraph_subgraphs_master} depict
examples of graphs with subgraph nodes. In Figure \ref{fig:OptiGraph_subgraphs} we have a graph $\mathcal{G}$ that contains three subgraphs
$\mathcal{SG}_1,\mathcal{SG}_2,\mathcal{SG}_3$ with a total of nine nodes (three in each subgraph).  The graph $\mathcal{G}$  contains a global hyperedge that connects
to local nodes in the subgraphs. Figure \ref{fig:OptiGraph_subgraphs_master} shows a hierarchical  graph $\mathcal{G}$ with three subgraphs $\mathcal{SG}_1,\mathcal{SG}_2,\mathcal{SG}_3$
and nine total nodes.  The graph $\mathcal{G}$  contains a global node $n_0$ that is connected to nodes in the subgraphs.  This type of structure arises when there is a parent (master)
optimization problem that is connected to children subproblems.

\subsubsection{Example 2: Hierarchical Graph with Global Edge}

\begin{figure}
\begin{minipage}{\textwidth}
\begin{minipage}{\textwidth}
\begin{scriptsize}
\lstset{language=Julia,breaklines = true}
\begin{lstlisting}[caption = Hierarchical Connectivity using Global Edge, label = {code:example2}]
    #Create low-level graph2
    graph2 = OptiGraph()           |\label{line:create_graph2}|

    @optinode(graph2,n4)                |\label{line:graph2_start}|
    @variable(n4, x >= 0);  @variable(n4, y >= 2)
    @constraint(n4,x + y >= 5); @objective(n4, Min, y)

    @optinode(graph2,n5)
    @variable(n5, x >= 0);  @variable(n5, y >= 2)
    @constraint(n5,x + y >= 5); @objective(n5, Min, y)

    @optinode(graph2,n6)
    @variable(n6, x >= 0); @variable(n6, y >= 2 )
    @constraint(n6,x + y >= 5); @objective(n6, Min, y)

    #Create graph2 linking constraint
    @linkconstraint(graph2, n4[:x] + n5[:x] + n6[:x] == 5)    |\label{line:graph2_end}|

    #Create low-level graph 3
    graph3 = OptiGraph()                   |\label{line:graph3_start}|

    @optinode(graph3,n7)
    @variable(n7, x >= 0); @variable(n7, y >= 2)
    @constraint(n7,x + y >= 7); @objective(n7, Min, y)

    @optinode(graph3,n8)
    @variable(n8, x >= 0);  @variable(n8, y >= 2)
    @constraint(n8,x + y >= 7); @objective(n8, Min, y)

    @optinode(graph3,n9)
    @variable(n9,x >= 0); @variable(n9, y >= 2)
    @constraint(n9,x + y >= 7); @objective(n9, Min, y)

    #Create graph3 linking constraint
    @linkconstraint(graph3, n7[:x] + n8[:x] + n9[:x] == 7)      |\label{line:graph3_end}|

    #Create high-level graph0
    graph0 = OptiGraph()                |\label{line:ex2_graph0_create}|

    #Add subgraphs to graph0
    add_subgraph!(graph0,graph1)          |\label{line:ex2_graph0_subgraph_start}|
    add_subgraph!(graph0,graph2)
    add_subgraph!(graph0,graph3)          |\label{line:ex2_graph0_subgraph_end}|

    #Add linking constraint to graph0 connecting its subgraphs
    @linkconstraint(graph0,n3[:x] + n5[:x] + n7[:x] == 10)      |\label{line:ex2_graph0_linkconstraint}|

    #Optimize with GLPK
    optimize!(graph0,GLPK.Optimizer)             |\label{line:ex2_solve}|
\end{lstlisting}
\end{scriptsize}
\end{minipage}
\\
\begin{minipage}{0.5\textwidth}
\includegraphics[width=7.3cm,keepaspectratio]{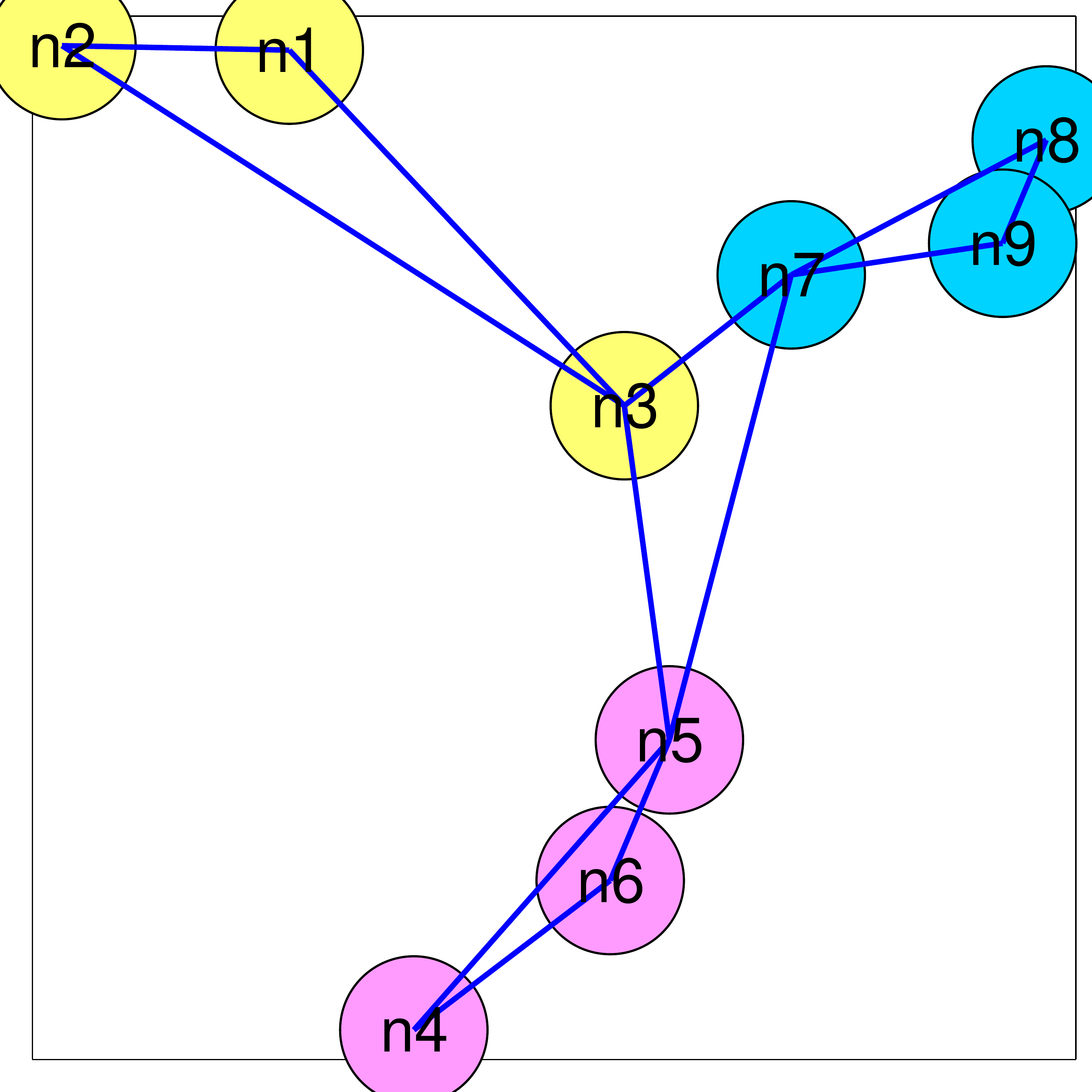}
\end{minipage}%
\begin{minipage}{0.5\textwidth}
\hspace{0.2cm}
\includegraphics[width=7.3cm,keepaspectratio]{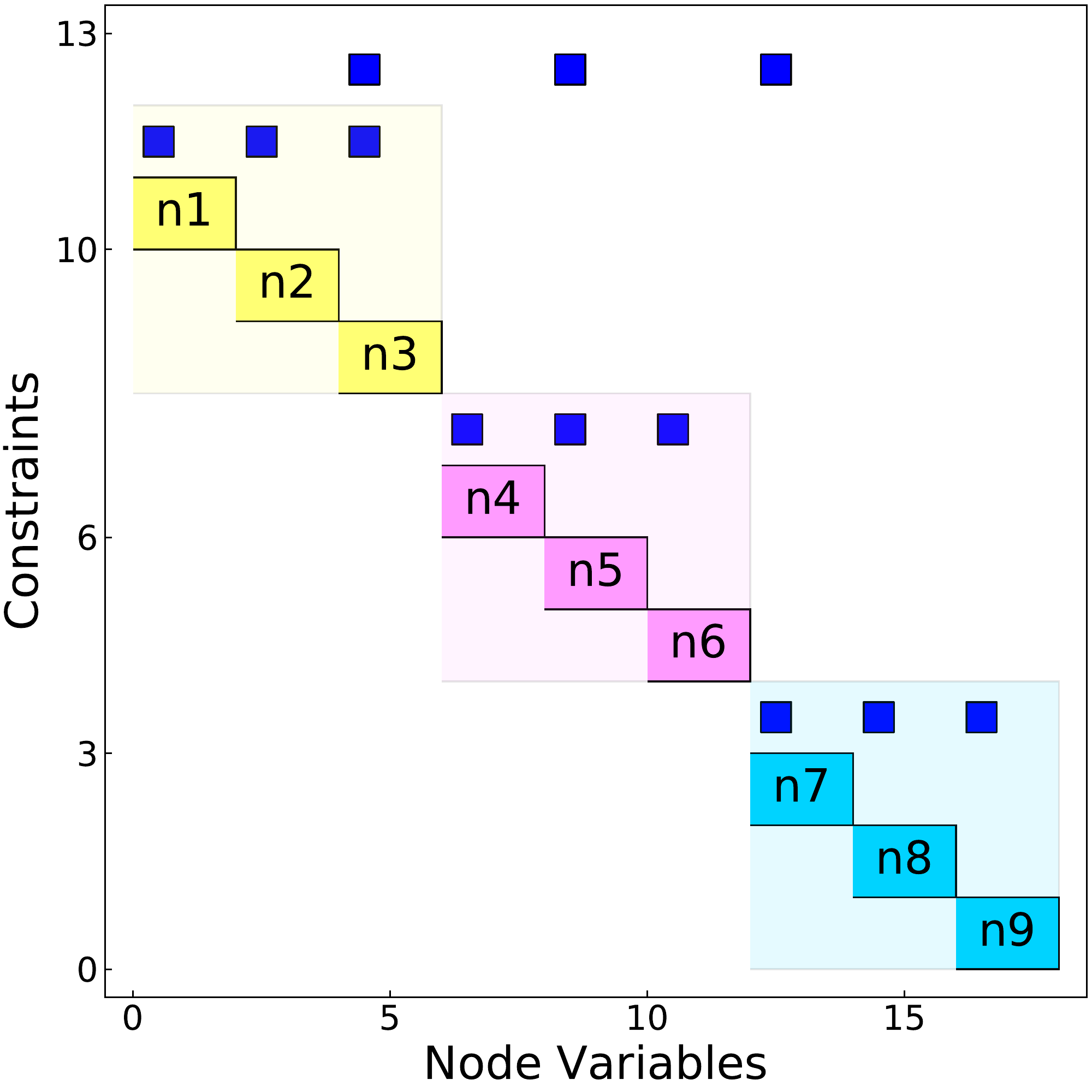}
\end{minipage}
\captionof{figure}{Output visuals for Code Snippet \ref{code:example2} showing hierarchical structure of the {\tt OptiGraph} \\ with three subgraphs.}
\label{fig:example2_result}
\end{minipage}
\end{figure}

An {\tt OptiGraph} object manages its own nodes and edge in a {\em self-contained} manner (without requiring references to other higher-level graphs).
Consequently, we can define subgraphs independently (in a modular manner) and these can be coupled together by using global edges or nodes defined on a high-level graph.
This feature is key because each node can have its own syntax (syntax does not need to be consistent or non-redundant over the entire model).
This is a fundamental difference with algebraic modeling languages (syntax has to be consistent and non-redundant over the entire model). The formulation in \eqref{eq:example2}
illustrates how to express hierarchical connectivity using a global edge.  Equation \eqref{eq:example2_objective} is the summation of
every node objective function in the graph; \eqref{eq:example2_node1}, \eqref{eq:example2_node2} and \eqref{eq:example2_node3}
describe node constraints; \eqref{eq:example2_subgraph1}, \eqref{eq:example2_subgraph2}, and \eqref{eq:example2_subgraph3}
represent link constraints within each subgraph; and \eqref{eq:example2_parent} defines a link constraint at the higher level graph
(that links nodes from each individual subgraph). Formulation \eqref{eq:example2} can be expressed as a hierarchical
{\tt OptiGraph} using the {\tt add\_subgraph!} function. This functionality is shown in Code Snippet \ref{code:example2}, where we extend
 {\tt graph1} from Code Snippet \ref{code:example1_code_snippet}.
We create a new graph called {\tt graph2} on Line \ref{line:create_graph2} and setup nodes and link them together
on Lines \ref{line:graph2_start} through \ref{line:graph2_end}.  We also construct {\tt graph3} in the same fashion
on Lines \ref{line:graph3_start} through \ref{line:graph3_end}.  Next, we create {\tt graph0}
on Line \ref{line:ex2_graph0_create} and add graphs {\tt graph1}, {\tt graph2}, and {\tt graph3} as subgraphs to {\tt graph0} on
Lines \ref{line:ex2_graph0_subgraph_start} through \ref{line:ex2_graph0_subgraph_end}.
We add a linking constraint to {\tt graph0} that couples nodes on each subgraph
on Line \ref{line:ex2_graph0_linkconstraint} and solve the graph on Line \ref{line:ex2_solve}.  We present the graph visualization
in Figure \ref{fig:example2_result}. Here we can see the hierarchical structure of the {\tt OptiGraph} and the local and global coupling constraints.
This structure is compatible with that shown in Figure \ref{fig:OptiGraph_subgraphs}.

\begin{subequations}\label{eq:example2}
    \begin{align}
        & \min \sum_{i = 1:9} y_{n_i}  \quad & \text{(Node Objectives)} \label{eq:example2_objective}\\
        \textrm{s.t.}\; & x_{n_i} \ge 0, y_{n_i} \ge 2, x_{n_i} + y_{n_i} \ge 3, i \in \{1,2,3\} \quad & \text{(Subgraph 1 Constraints)} \label{eq:example2_node1} \\
        & x_{n_i}  \ge 0,  y_{n_i}  \ge 2, x_{n_i} + y_{n_i} \ge 5, i \in \{4,5,6\}  \quad & \text{(Subgraph 2 Constraints)}         \label{eq:example2_node2} \\
        & x_{n_i} \ge 0, y_{n_i}  \ge 2, x_{n_i} + y_{n_i} \ge 7, i \in \{7,8,9\}    \quad & \text{(Subgraph 3  Constraints)}   \label{eq:example2_node3} \\
        & x_{n_1} + x_{n_2} + x_{n_3} = 3  \quad & \text{(Subgraph 1 Link Constraint)} \label{eq:example2_subgraph1} \\
        & x_{n_4} + x_{n_5} + x_{n_6} = 5  \quad & \text{(Subgraph 2 Link Constraint)}  \label{eq:example2_subgraph2} \\
        & x_{n_7} + x_{n_8} + x_{n_9} = 7  \quad & \text{(Subgraph 3 Link Constraint)}  \label{eq:example2_subgraph3}\\
        & x_{n_3} + x_{n_5} + x_{n_7} = 10 \quad & \text{(Global Link Constraint)}  \label{eq:example2_parent}
    \end{align}
\end{subequations}

\FloatBarrier

\subsubsection{Example 3: Hierarchical Graph with Global Node}

\begin{figure}
\noindent\begin{minipage}{\textwidth}
\begin{minipage}{\textwidth}
\begin{scriptsize}
\lstset{language=Julia,breaklines = true}
\begin{lstlisting}[caption= Hierarchical Connectivity using Global Node, label = {code:example3}]
#Create graph0
graph0 = OptiGraph()     |\label{line:ex3_graph0_start}|

#Create a node on graph0
@optinode(graph0,n0)
@variable(n0,x)
@constraint(n0,x >= 0)      |\label{line:ex3_graph0_end}|

#Add subgraphs to graph0
add_subgraph!(graph0,graph1)        |\label{line:ex3_subgraph_start}|
add_subgraph!(graph0,graph2)
add_subgraph!(graph0,graph3)        |\label{line:ex3_subgraph_end}|

#Create linking constraints on graph0 connecting it to its subgraphs
@linkconstraint(graph0,n0[:x] + n3[:x] == 3)    |\label{line:ex3_linkconstraint_start}|
@linkconstraint(graph0,n0[:x] + n5[:x] == 5)
@linkconstraint(graph0,n0[:x] + n7[:x] == 7)     |\label{line:ex3_linkconstraint_end}|

#Optimize with GLPK
optimize!(graph0,GLPK.Optimizer)                 |\label{line:ex3_solve}|
\end{lstlisting}
\end{scriptsize}
\end{minipage}\\
\\
\begin{minipage}{0.5\textwidth}
\includegraphics[width=7.3cm,keepaspectratio]{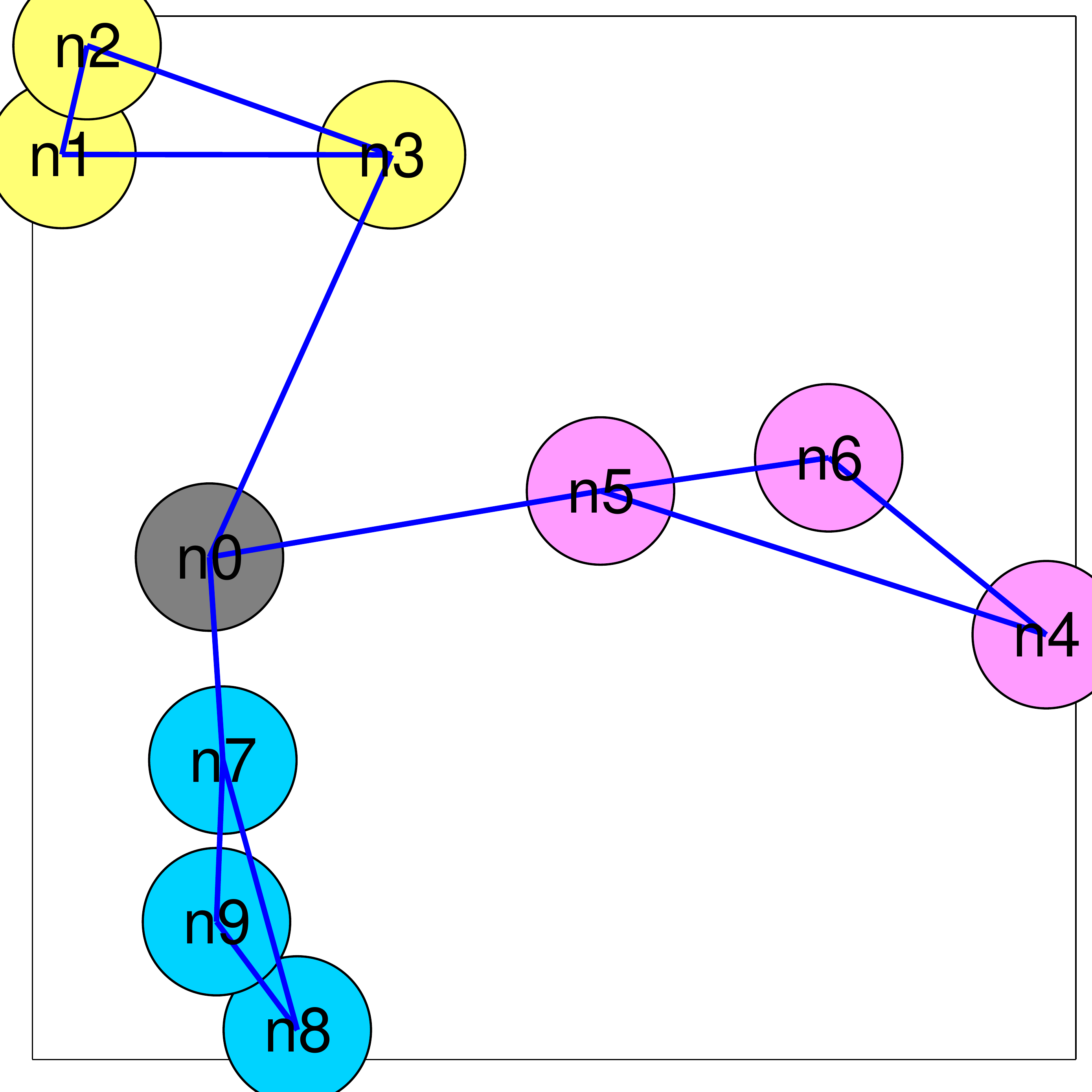}
\end{minipage}%
\begin{minipage}{0.5\textwidth}
\hspace{0.2cm}
\includegraphics[width=7.3cm,keepaspectratio]{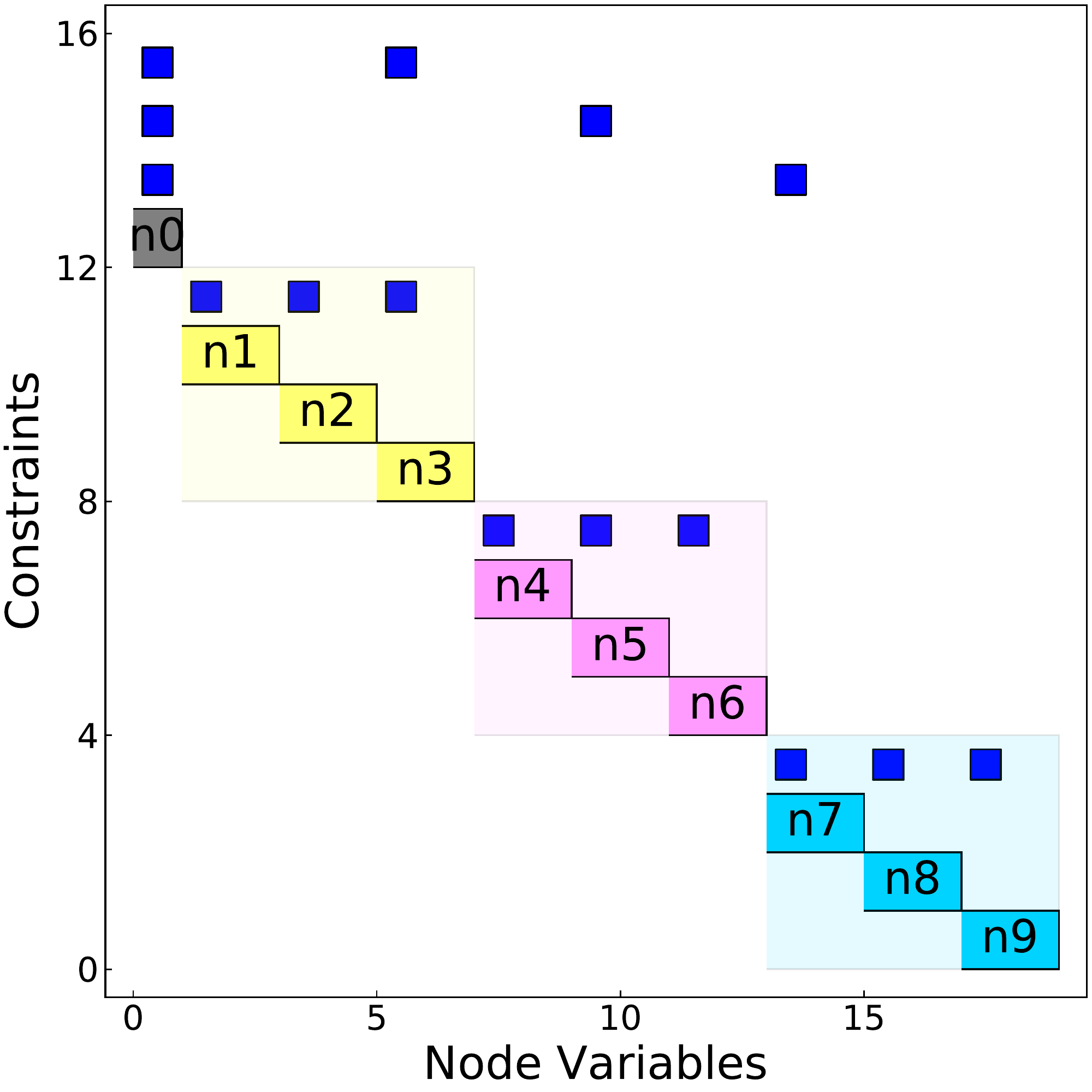}
\end{minipage}
\captionof{figure}{Output visuals for Code Snippet \ref{code:example3} showing hierarchical structure of the {\tt OptiGraph} \\ with three subgraphs.}
\label{fig:example3_result}
\end{minipage}
\end{figure}

We can express hierarchical connectivity within a {\tt OptiGraph} by defining a global node that is connected with subgraph nodes.
Formulation \eqref{eq:example3} illustrates this idea; this is analogous to \eqref{eq:example2} where we have removed
the high level linking constraint \eqref{eq:example2_parent} and have replaced it with a high level node \eqref{eq:example3_node0}
and three linking constraints that couple the graph to its subgraphs \eqref{eq:example3_parent}.  The implementation of the formulation
in \eqref{eq:example3} is shown in Code Snippet \ref{code:example3}. Here, we assume that we already have {\tt graph1}, {\tt graph2}, and {\tt graph3}
defined from Snippets \ref{code:example1_code_snippet} and \ref{code:example2}. We recreate {\tt graph0} and setup the node {\tt n0} on Lines \ref{line:ex3_graph0_start}
through \ref{line:ex3_graph0_end}. We add subgraphs {\tt graph1}, {\tt graph2}, and {\tt graph3} on Lines \ref{line:ex3_subgraph_start}
through \ref{line:ex3_subgraph_end} like in the previous snippet, and add linking constraints that connect
node {\tt n0} to nodes in each subgraph on Lines \ref{line:ex3_linkconstraint_start} through \ref{line:ex3_linkconstraint_end}.
We  solve the newly created {\tt graph0} on Line \ref{line:ex3_solve} and  present the visualization in Figure \ref{fig:example3_result}. This structure is compatible with
that shown in Figure \ref{fig:OptiGraph_subgraphs_master}.

\begin{subequations}\label{eq:example3}
    \begin{align}
       \min\; & \sum_{i = 1:9} y_{n_i}  \quad & \text{(Objective)} \label{eq:example3_objective}\\
        \textrm{s.t.}\; & x_{n_i} \ge 0, y_{n_i} \ge 2, x_{n_i} + y_{n_i} \ge 3, i \in \{1,2,3\} \quad & \text{(Subgraph 1 Constraints)} \label{eq:example3_node1} \\
        & x_{n_i}  \ge 0,  y_{n_i}  \ge 2, x_{n_i} + y_{n_i} \ge 5, i \in \{4,5,6\}  \quad & \text{(Subgraph 2 Constraints)}         \label{eq:example3_node2} \\
        & x_{n_i} \ge 0, y_{n_i}  \ge 2, x_{n_i} + y_{n_i} \ge 7, i \in \{7,8,9\}    \quad & \text{(Subgraph 3 Constraints)}   \label{eq:example3_node3} \\
        & x_{n_1} + x_{n_2} + x_{n_3} = 3  \quad & \text{(Subgraph 1 Link Constraint)} \label{eq:example3_subgraph1} \\
        & x_{n_4} + x_{n_5} + x_{n_6} = 5  \quad & \text{(Subgraph 2 Link Constraint)}  \label{eq:example3_subgraph2} \\
        & x_{n_7} + x_{n_8} + x_{n_9} = 7  \quad & \text{(Subgraph 3 Link Constraint)}  \label{eq:example3_subgraph3}\\
        & x_{n_0}  \ge 0 \quad & \text{(Graph Constraint)} \label{eq:example3_node0}\\
        & x_{n_0} + x_{n_3} = 3, x_{n_0} + x_{n_5} = 5,x_{n_0} + x_{n_7} = 7 \quad & \text{(Global Link Constraints)}  \label{eq:example3_parent}
    \end{align}
\end{subequations}

\subsection{Processing and Manipulating OptiGraphs}

In addition to the graph construction functions presented in the previous examples ({\tt @optinode}, {\tt @linkconstraint}, {\tt add\_subgraph!}),
it is also possible to query an {\tt OptiGraph} object to retrieve its attributes.  Table \ref{table:primary_functions} summarizes
the main {\tt Plasmo.jl} functions used to create and inspect an {\tt OptiGraph}.  We inspect the nodes, edges, linking constraints, and subgraphs
using {\tt getoptinodes}, {\tt getoptiedges}, {\tt getlinkconstraints}, and {\tt getsubgraphs} functions, and we can
collect \emph{all} of the corresponding graph attributes using recursive versions of these functions
({\tt all\_nodes}, {\tt all\_edges}, {\tt all\_linkconstraints} and {\tt all\_subgraphs}).

\begin{table}
\scriptsize
\begin{center}
\begin{tabular}{|l|l|}
\hline
Function & Description \\
\hline
{\tt \textbf{OptiGraph}()} & Create a new OptiGraph object.\\
{\tt \textcolor{Blue}{\textbf{@optinode}}(graph::OptiGraph,expr::Expr)} & Create OptiNodes using {\tt Julia} expression \\
{\tt \textcolor{Blue}{\textbf{@linkconstraint}}(graph::OptiGraph,expr::Expr)} & Create linking constraint between nodes in {\tt graph} using {\tt expr}\\
{\tt \textbf{add\_subgraph!}(graph::OptiGraph,sg::OptiGraph)} & Add subgraph {\tt sg} to {\tt graph}\\
\hline
{\tt \textbf{getoptinodes}(graph::OptiGraph)} & Retrieve local OptiNodes in {\tt graph}\\
{\tt \textbf{getoptiedges}(graph::OptiGraph)} & Retrieve local OptiEdges in {\tt graph}\\
{\tt \textbf{getlinkconstraints}(graph::OptiGraph)} & Retrieve linking constraints in {\tt graph}\\
{\tt \textbf{getsubgraphs}(graph::OptiGraph)} & Retrieve subgraphs in {\tt graph}\\
\hline
{\tt \textbf{all\_optinodes}(graph::OptiGraph)} & Retrieve all OptiNodes in {\tt graph} (including subgraphs)\\
{\tt \textbf{all\_optiedges}(graph::OptiGraph)} & Retrieve OptiEdges in {\tt graph} (including subgraphs)\\
{\tt \textbf{all\_linkconstraints}(graph::OptiGraph)} & Retrieve all linking constraints in {\tt graph} (including subgraphs)\\
{\tt \textbf{all\_subgraphs}(graph::OptiGraph)} & Retrieve all subgraphs in {\tt graph} (including subgraphs)\\
\hline
\end{tabular}
\caption{Overview of {\tt OptiGraph} construction and query functions in {\tt Plasmo.jl}.}\label{table:primary_functions}
\end{center}
\end{table}

Thus far we have discussed how to enable model construction using a \emph{bottom-up} approach.  Specifically, we showed how to
construct an {\tt OptiGraph} by adding subgraphs. We now show how to create an {\tt OptiGraph} using a \emph{top-down} approach. Specifically, we show how to partition an {\tt OptiGraph},  construct subgraphs from the partitions, and use these to derive an alternative representation of the {\tt OptiGraph}. As part of this, we will discuss how {\tt Plasmo.jl} interfaces to standard graph partitioning tools.

\subsubsection{Hypergraph Partitioning}\label{sec:hypergraph_partitioning}

The {\tt OptiGraph} is, by default, a hypergraph; as such, we can naturally exploit hypergraph partitioning capabilities.  Here, we focus on hypergraph partitioning concepts but we note that our
discussion also applies to standard graph partitioning frameworks (a hypergraph can be projected to various standard graph representations). Hypergraph partitioning has received significant interest
in the last few years because it naturally
represents complex {\em non-pairwise} relationships and more accurately captures coupling in such structures compared to traditional graphs. Popular hypergraph
partitioning tools include the well-known {\tt hMetis} and {\tt PaToH} packages,
as well as the comprehensive {\tt Zoltan} \cite{Devine2006} software suite which provides
hypergraph partitioning algorithms for {\tt C}, {\tt C++}, and {\tt Fortran} applications.
More recent frameworks have made advances to create large-scale hypergraphs \cite{Jiang2018},
perform fast hypergraph partitioning \cite{Mayer2018HYPEMH}, and create high-quality hypergraph partitions \cite{shhmss2016alenex}.

To provide an overview of hypergraph partitioning techniques, we use notation that is similar to that of an {\tt OptiGraph}.
A hypergraph contains a set of hypernodes $\mathcal{N}(\mathcal{H})$ and hyperedges $\mathcal{E}(\mathcal{H})$ where we denote
the hypergraph containing hypernodes and hyperedges as $\mathcal{H}(\mathcal{N},\mathcal{E})$.
In hypergraph partitioning, one seeks to partition the set of nodes $\mathcal{N}(\mathcal{H})$ into a collection $\mathcal{P}$ of at most $k$
disjoint subsets such that $\mathcal{P} = \{P_1,P_2,...,P_k\}$ while minimizing an
objective function over the edges such as \eqref{eq:edge-cut} or \eqref{eq:connectivity} subject to a
balance constraint \eqref{eq:hypergraph_balance} (such that partitions are roughly the same size).

\begin{subequations}\label{eq:hypergraph-partitioning}
    \begin{align}
        & \Phi_{cut}(\mathcal{P}) = \sum_{e \in \mathcal{E}_{cut}(\mathcal{P})} w(e) \label{eq:edge-cut}\\
        & \Phi_{con}(\mathcal{P}) = \sum_{e \in \mathcal{E}_{cut}(\mathcal{P})} w(e) (\lambda(e) - 1)  \label{eq:connectivity}\\
        & \frac{1}{k}\sum_{n \in \mathcal{N}(\mathcal{H})} s(n) - \epsilon_{max}  \le \sum_{n \in \mathcal{P}_i} s(n) \le \frac{1}{k}\sum_{n \in \mathcal{N}(\mathcal{H})} s(n) +
            \epsilon_{max}, \quad \forall i \in \{1,...,k\}  \label{eq:hypergraph_balance}
    \end{align}
\end{subequations}

Here, $\Phi_{cut}$ and $\Phi_{con}$ are the most prominent hypergraph partitioning objectives (called the minimum edge-cut and
minimum connectivity), where $\mathcal{E}_{cut}(\mathcal{P})$ is the set of \emph{cut edges} of the partitions in $\mathcal{P}$
(i.e., all edges that \emph{cross} partitions defined by $\mathcal{P}$). The formulation in \eqref{eq:hypergraph-partitioning} introduces edge weights $w(e) \rightarrow \mathbb{R}$ for each
edge $e \in \mathcal{E}(\mathcal{H})$ and node sizes $s(n) \rightarrow \mathbb{R}$ for each node $n \in \mathcal{N}(\mathcal{H})$ which can be used to express specific attributes to the partitioner.
For instance, large edge weights typically express tight coupling or high communication volume
and nodes sizes often represent computational load. The objective \eqref{eq:connectivity} includes the connectivity
value $\lambda(e)$ which denotes the number of partitions connected by a hyperedge $e$.
We also define the parameter $\epsilon_{max}> 0$ in \eqref{eq:hypergraph_balance}
which specifies the maximum imbalance tolerance of the partitions.  Lower values of $\epsilon_{max}$ enforce equal-sized partitions and higher values
allow for solutions with disparate partition sizes. The lower bound constraint in \eqref{eq:hypergraph_balance} is not always included in
partitioning algorithm implementations.

\begin{figure}[]
\centering
\begin{subfigure}[t]{0.32\textwidth}
    \centering
    \includegraphics[scale=0.40]{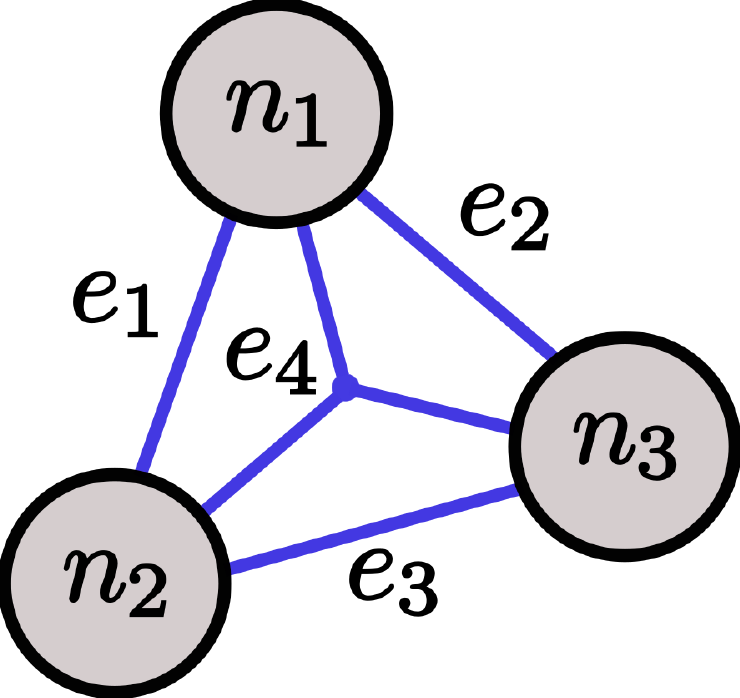}
    \caption{Hypergraph}
    \label{fig:row_net_hypergraph}
\end{subfigure}
\begin{subfigure}[t]{0.32\textwidth}
    \centering
    \includegraphics[scale=0.40]{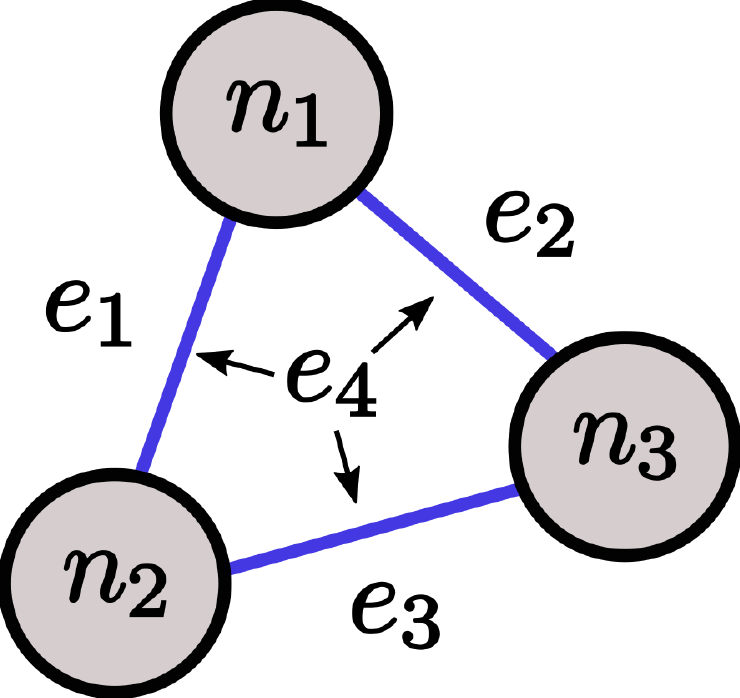}
    \caption{Standard Graph}
    \label{fig:clique_expansion}
\end{subfigure}
\begin{subfigure}[t]{0.32\textwidth}
    \centering
    \includegraphics[scale=0.40]{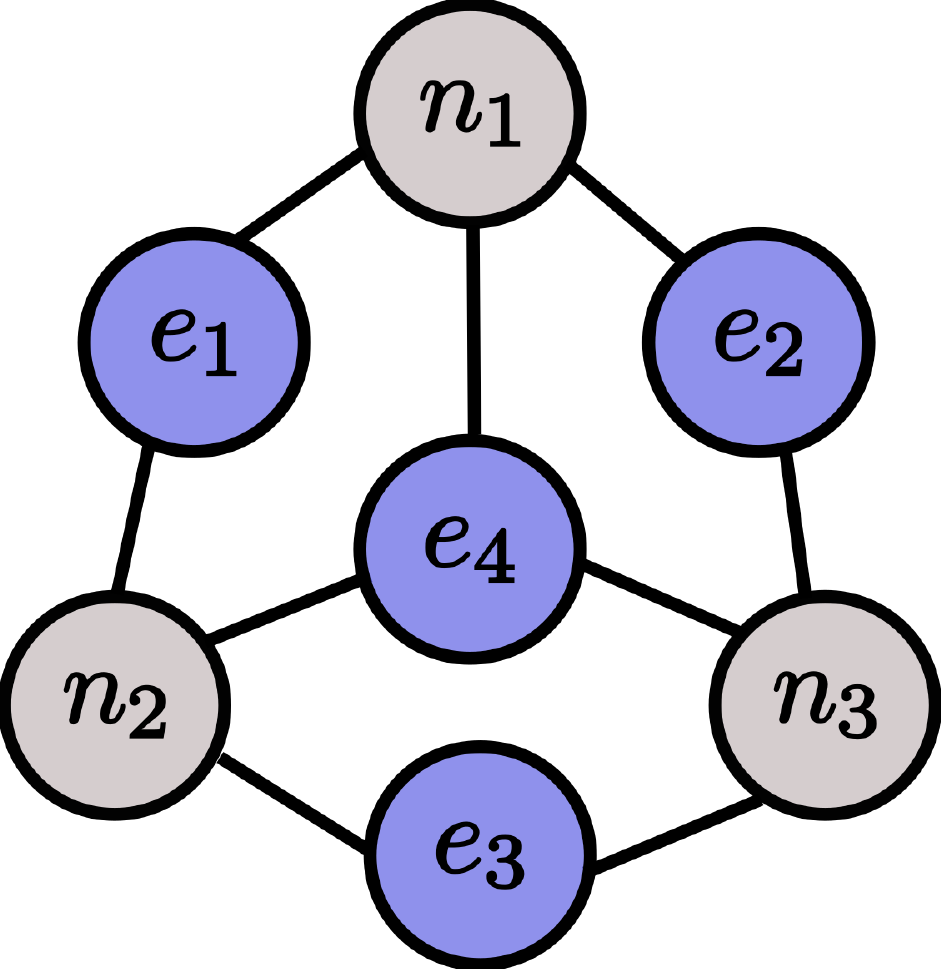}
    \caption{Bipartite Graph}
    \label{fig:star_expansion}
\end{subfigure}
\caption{Typical graph representations used in partitioning tools.  A hypergraph can be
projected to a standard graph or a bipartite graph.}
\label{fig:graph_representations}
\end{figure}

The hypergraph is an intuitive representation for an {\tt OptiGraph} but other representations are also possible.
The hypergraph in Figure \ref{fig:row_net_hypergraph} can be projected to a standard graph representation (shown in Figure \ref{fig:clique_expansion}) or to a
bipartite representation (shown in Figure \ref{fig:star_expansion}).
Standard graph representations of optimization problems can utilize mature partitioning tools such as {\tt Metis} or community detection strategies.
This provides flexibility to experiment with different techniques. In the remainder of the discussion we utilize the hypergraph representation for
partitioning but we highlight that a broader range of partitioning strategies can be captured in the proposed framework.

\subsubsection{OptiGraph Partitioning in Plasmo.jl}\label{sec:modelgraph_part_and_topology}

\begin{figure}[]
\centering
\begin{subfigure}[t]{0.32\textwidth}
    \centering
    \includegraphics[width = 4.5cm,keepaspectratio]{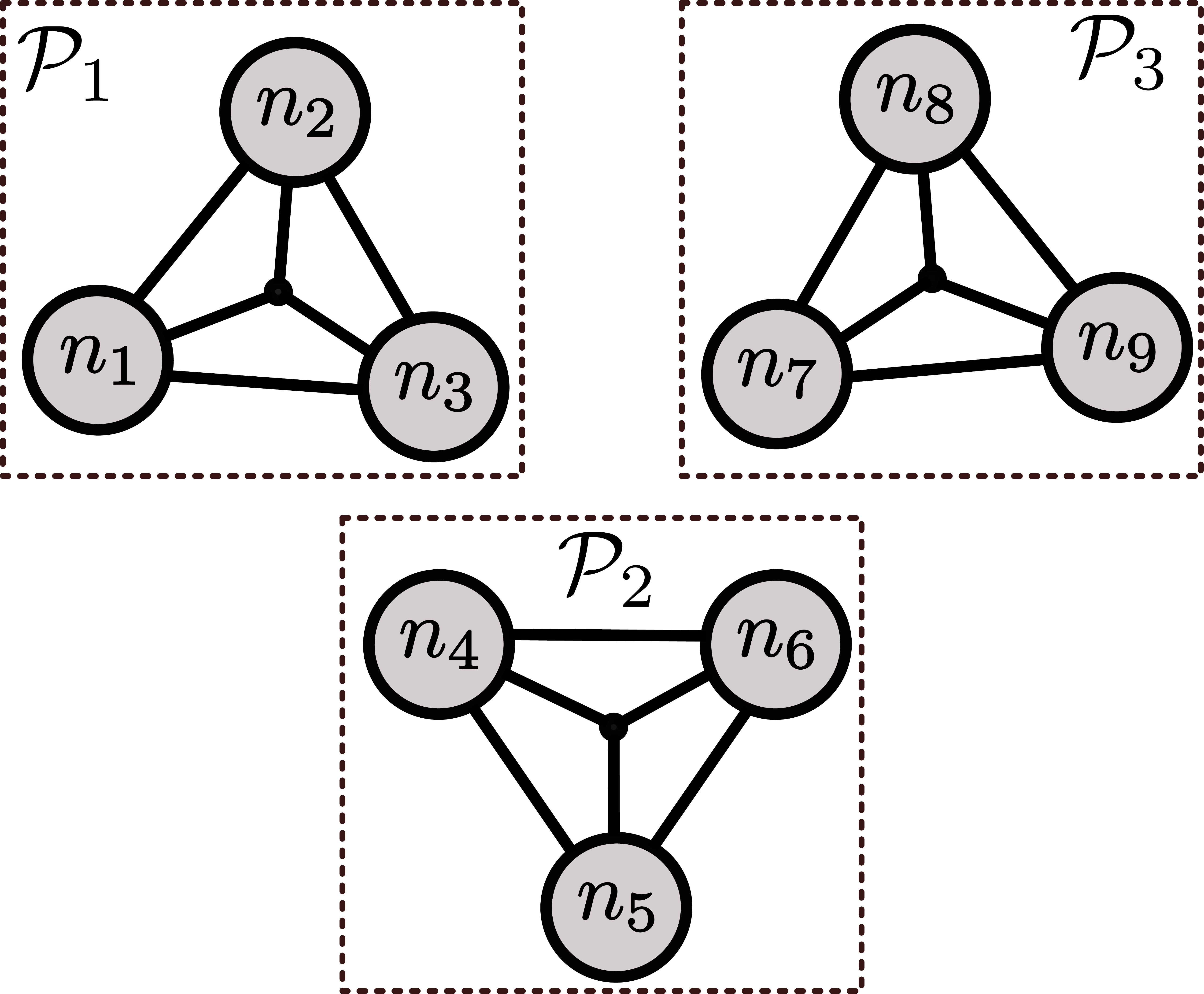}
    \caption{{\tt Partition}}
    \label{fig:modelgraph_partition}
\end{subfigure}
\begin{subfigure}[t]{0.32\textwidth}
    \centering
    \includegraphics[width = 4.5cm,keepaspectratio]{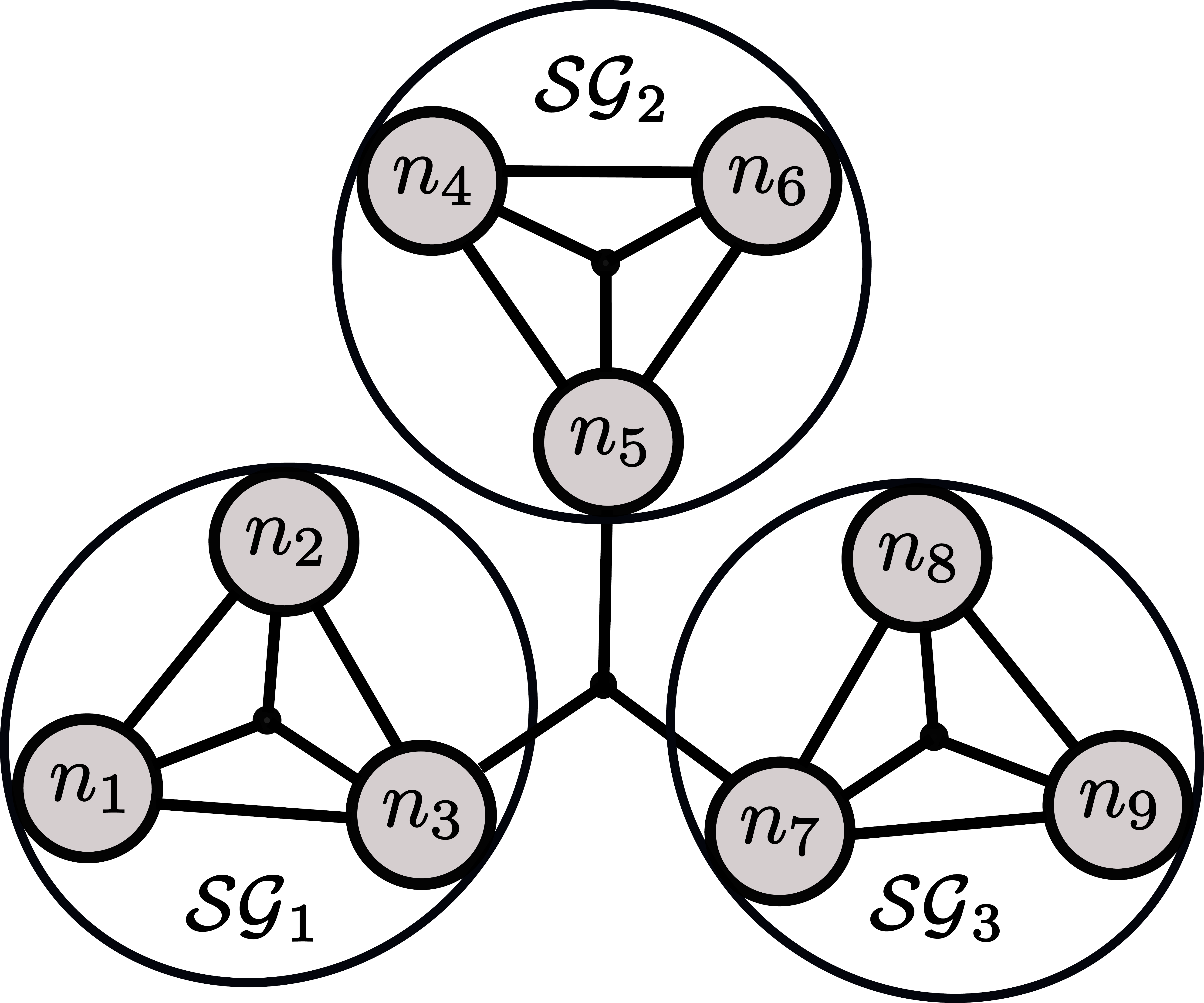}
    \caption{{\tt make\_subgraphs!}}
    \label{fig:modelgraph_subgraphs}
\end{subfigure}
\begin{subfigure}[t]{0.32\textwidth}
    \centering
    \includegraphics[width = 4.5cm,keepaspectratio]{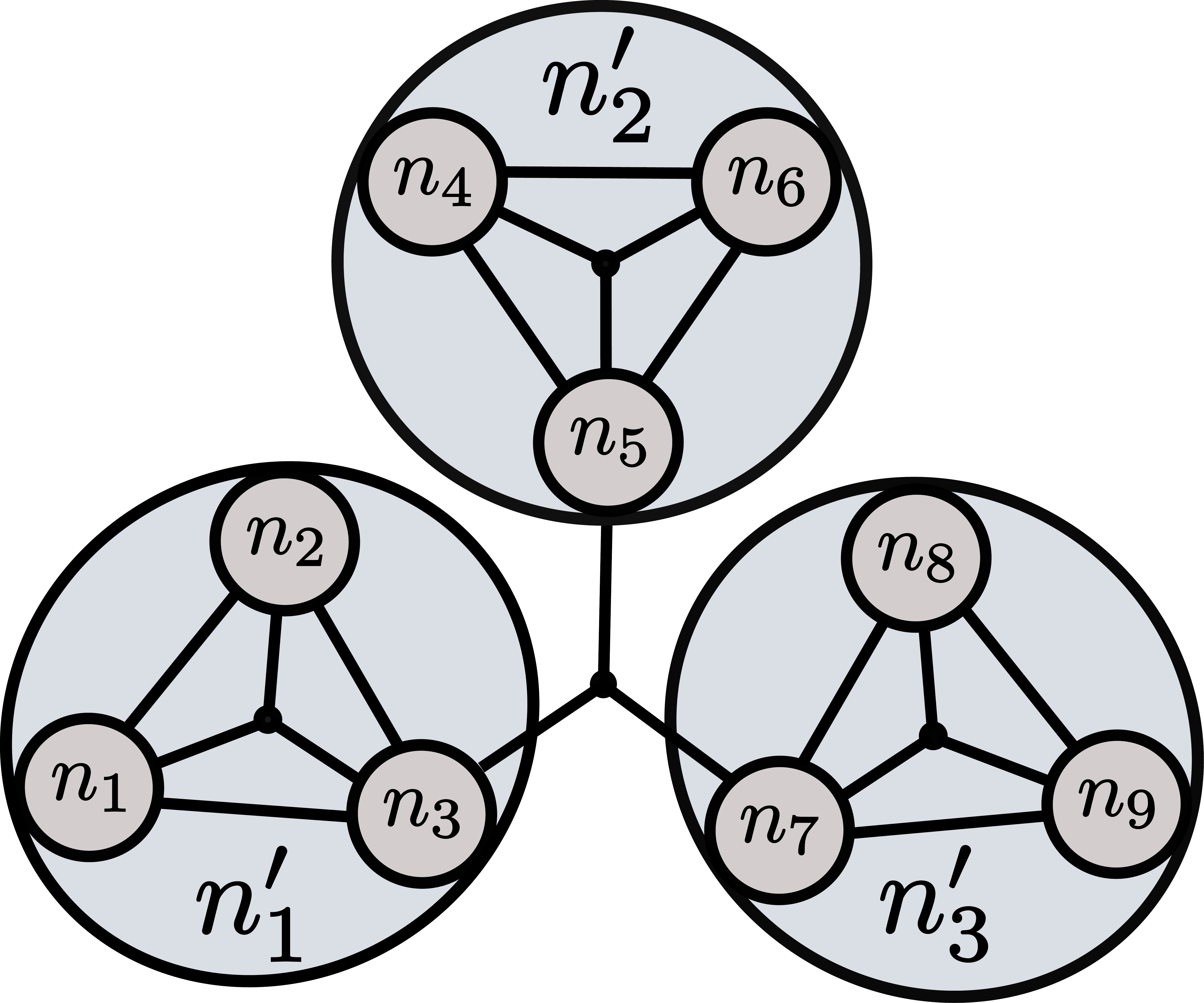}
    \caption{{\tt aggregate}}
    \label{fig:modelgraph_combine}
\end{subfigure}
\caption{Core partitioning functionality in {\tt Plasmo.jl}. (left) A {\tt Partition} with nine {\tt OptiNodes}, (middle)
corresponding {\tt OptiGraph}
containing three subgraphs, and (right) subgraphs aggregated into {\tt OptiNodes}.}
\label{fig:partitioning_modelgraphs}
\end{figure}

Here we discuss partitioning capabilities implemented in {\tt Plasmo.jl}. Figure \ref{fig:partitioning_modelgraphs} depicts the basic aspects of the partitioning framework.
We show how to create {\tt OptiGraph} partitions (with a {\tt Partition} object), how to formulate subgraphs (with the {\tt make\_subgraphs!} function), and how to \emph{combine} (aggregate)
subgraphs into stand-alone {\tt OptiNodes} (with the {\tt aggregate} function). For example,
Figure \ref{fig:modelgraph_partition} contains three partitions $P_1$,$P_2$, and $P_3$, Figure \ref{fig:modelgraph_subgraphs} shows the
corresponding subgraphs $\mathcal{SG}_1$,$\mathcal{SG}_2$, and $\mathcal{SG}_3$ created from the partitions, and Figure \ref{fig:modelgraph_combine} depicts
the resulting {\tt OptiNodes} $n_1^\prime$,$n_2^\prime$, and $n_3^\prime$ which represent optimization \emph{subproblems}. To perform hypergraph partitioning we use
the {\tt KaHyPar} \cite{shhmss2016alenex} hypergraph partitioner through the {\tt KaHyPar.jl} interface.
{\tt KaHyPar} targets the creation of high quality partitions and offers a straightforward {\tt C} library interface
which facilitates its connection with {\tt Plasmo.jl}. Throughout our examples we use the default {\tt KaHyPar} configuration which uses a direct multilevel k-way
algorithm with community detection initialization.

\begin{figure}[]
\centering
\begin{subfigure}[t]{0.32\textwidth}
    \centering
    \includegraphics[width = 4.5cm,keepaspectratio]{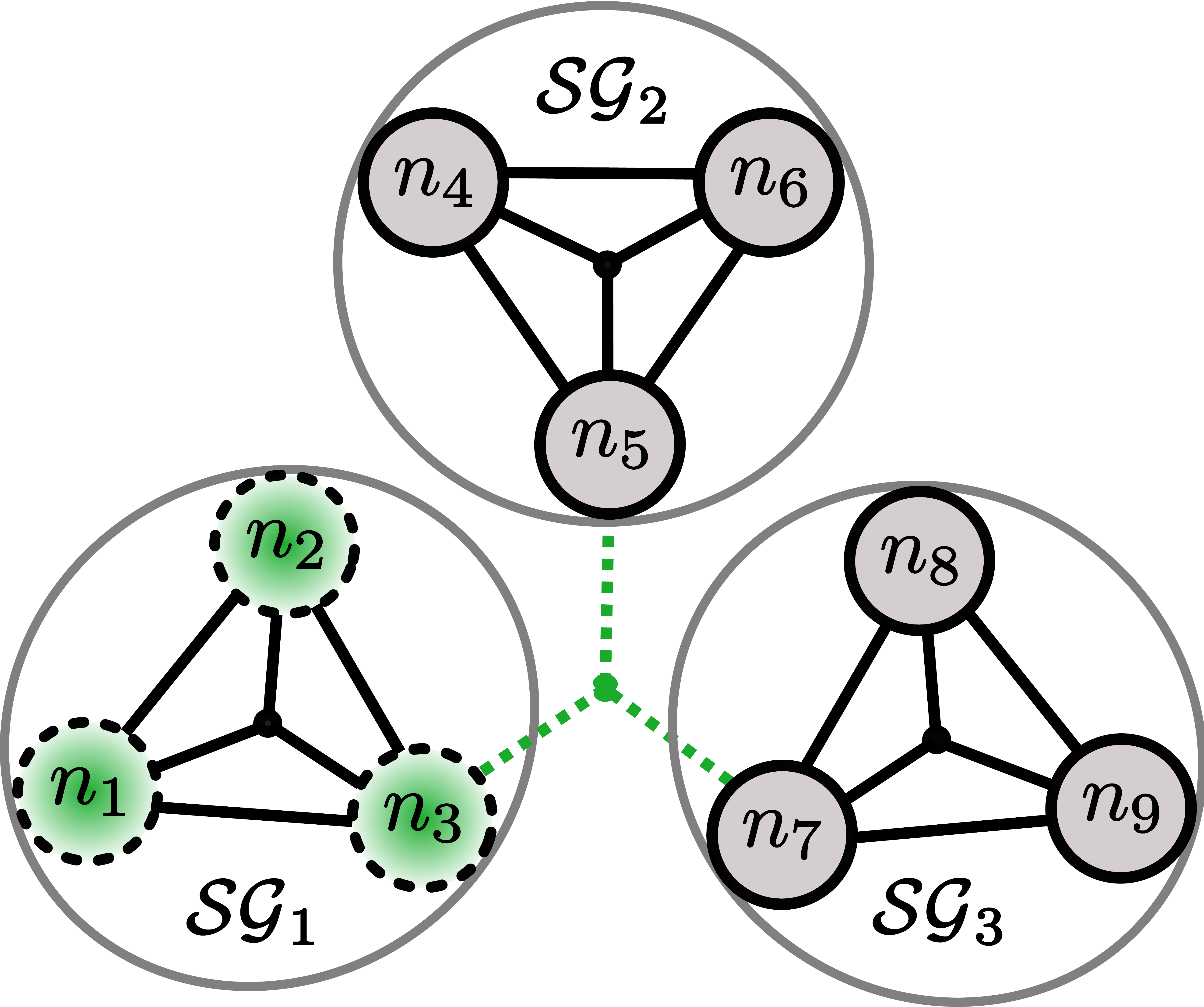}
    \caption{{\tt incident\_edges}}
    \label{fig:modelgraph_incident_edges}
\end{subfigure}
\begin{subfigure}[t]{0.32\textwidth}
    \centering
    \includegraphics[width = 4.5cm,keepaspectratio]{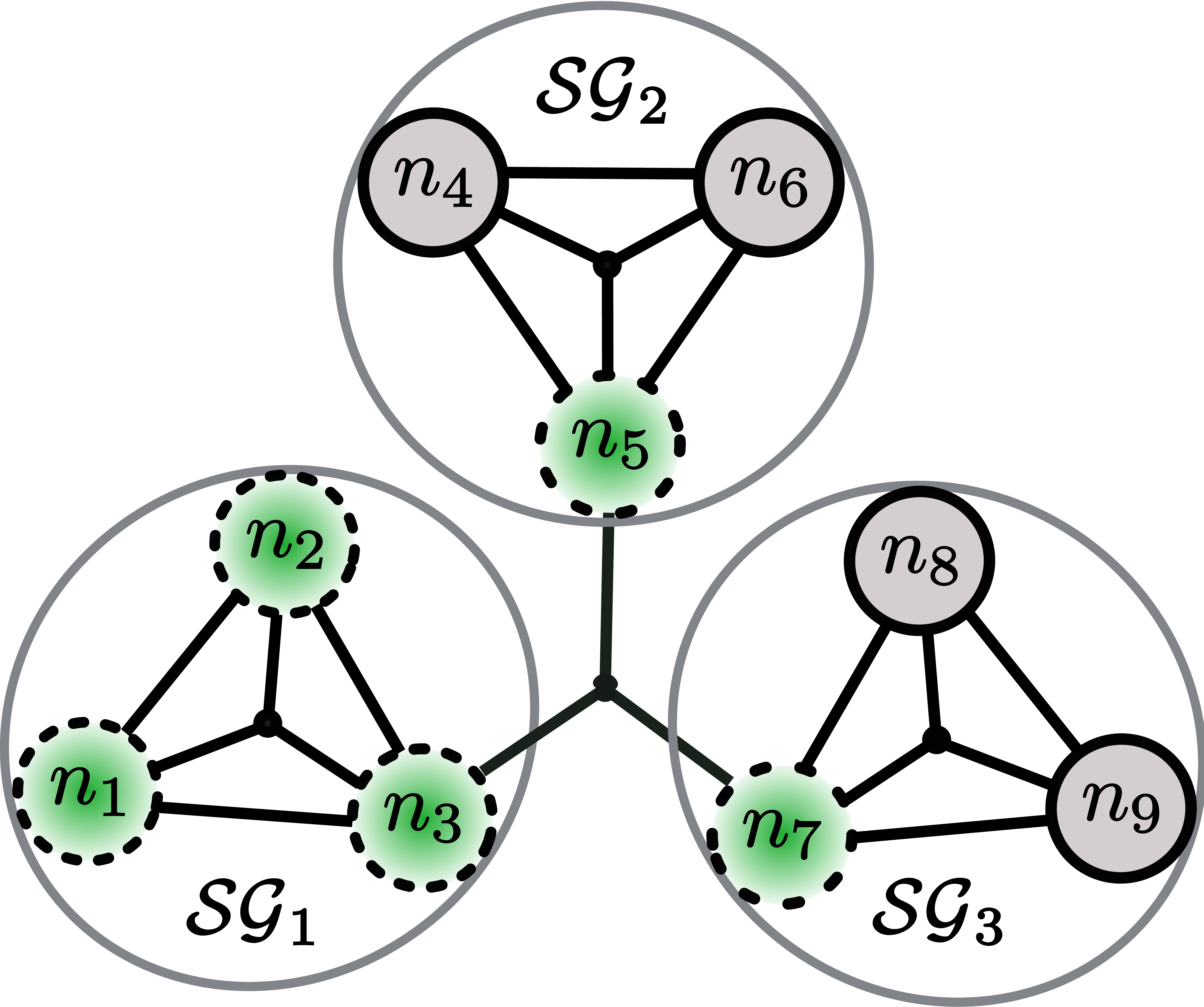}
    \caption{{\tt neighborhood}}
    \label{fig:modelgraph_neighborhood}
\end{subfigure}
\begin{subfigure}[t]{0.32\textwidth}
    \centering
    \includegraphics[width = 4.5cm,keepaspectratio]{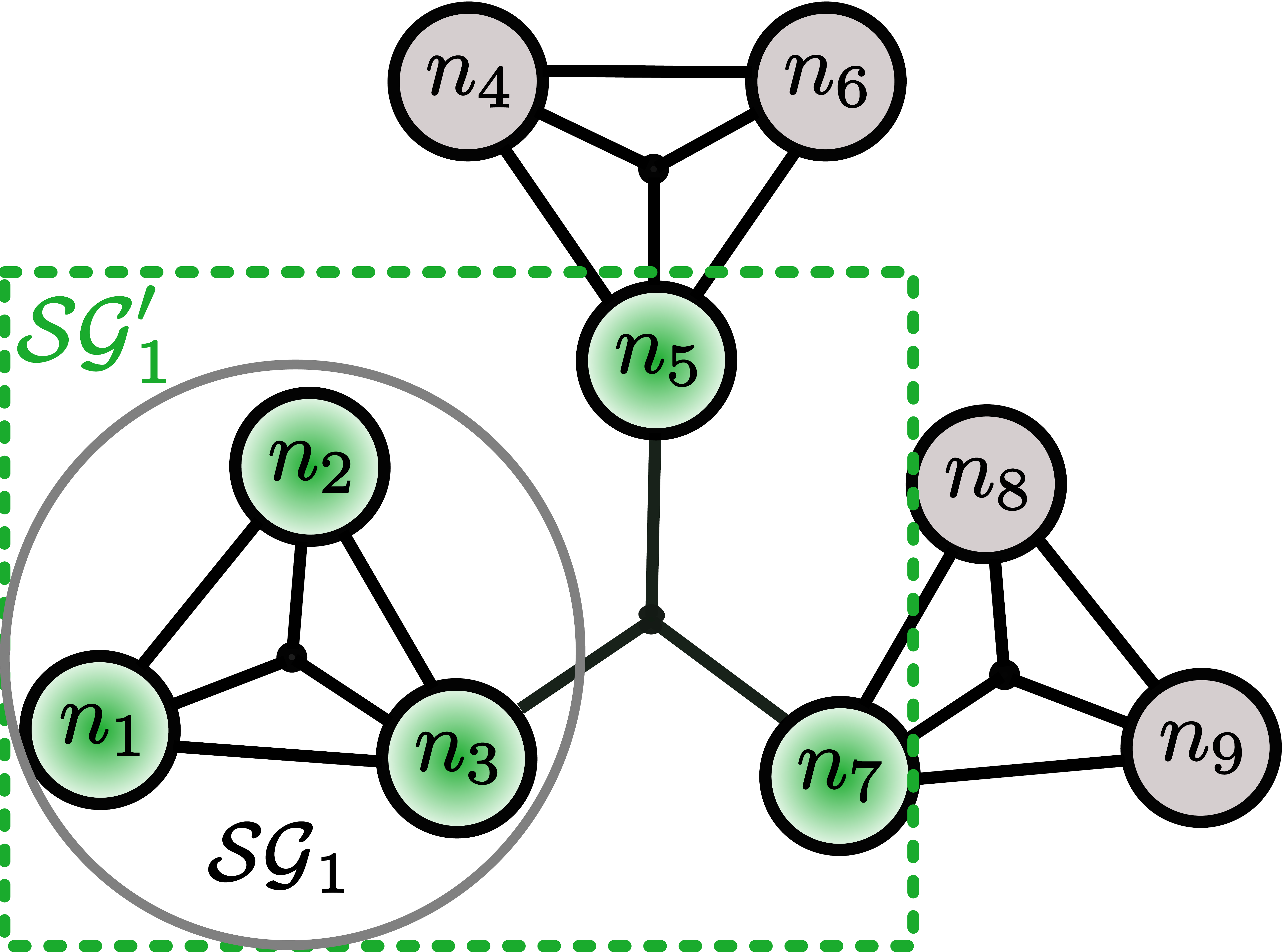}
    \caption{{\tt expand}}
    \label{fig:modelgraph_expand}
\end{subfigure}
\caption{Core graph topology functions in {\tt Plasmo.jl}. (left) querying incident edges to a subgraph, (middle) querying
a subgraph neighborhood, and (right) expanding a subgraph.}
\label{fig:modelgraph_topology_functions}
\end{figure}

The {\tt OptiGraph} object offers topology manipulation functionality that can be used, for instance, to modify subgraph structures  and formulate subproblems for algorithms. Figure \ref{fig:modelgraph_topology_functions} depicts core topology functions
we commonly use in the framework. We can query the  incident {\tt OptiEdges} (Figure \ref{fig:modelgraph_incident_edges}) to a set of {\tt OptiNodes} (or a subgraph) to inspect
coupling (i.e. inspect incident linking constraints).  We can also obtain the {\tt neighborhood} (Figure \ref{fig:modelgraph_neighborhood})
around a set of {\tt OptiNodes} to inspect an expanded problem domain, and we can {\tt expand}
(Figure \ref{fig:modelgraph_expand}) a subgraph into a larger domain and generate the corresponding subproblem.

Table \ref{table:partition_functions} summarizes the core graph partitioning and manipulation functions in {\tt Plasmo.jl}. The {\tt gethypergraph} function
returns a hypergraph object (extends a {\tt LightGraphs.jl} object) and a {\tt reference\_map} which maps the hypergraph elements back to the {\tt OptiGraph} (i.e. hypergraph node indices
are mapped back to {\tt OptiNodes}). We also introduce a {\tt Partition} object that describes how to formulate subgraphs within a graph.  As we will show, the {\tt Partition}
object is an intermediate interface to formulate subgraphs in a general way.

\begin{table}
\scriptsize
\begin{center}
\begin{tabular}{|l|}
\hline
Functions and Descriptions\\
\hline
Create a hypergraph representation of {\tt graph}.\\
{\tt hypergraph,ref = \textbf{gethypergraph}(graph::OptiGraph)}\\
\hline
Create a {\tt Partition} given an {\tt OptiGraph}, a {\tt vector} of integers and a mapping {\tt ref\_map}. \\
{\tt partition = \textbf{Partition}(graph::OptiGraph,vector::Vector\{Int\},mapping::Dict\{Int,OptiNode\})}\\
\hline
Reform {\tt graph} into subgraphs using {\tt partition}.\\
{\tt \textbf{make\_subgraphs!}(graph::OptiGraph,partition::Partition)}\\
\hline
Combine subgraphs in {\tt graph} such that {\tt n\_levels} of subgraphs remain.\\
{\tt \textbf{aggregate}(graph::OptiGraph,n\_levels::Int)}\\
\hline
Retrieve incident {\tt OptiEdges} of {\tt OptiNodes} in {\tt graph}.\\
{\tt \textbf{incident\_optiedges}(graph::OptiGraph,nodes::Vector\{OptiNode\})}\\
\hline
Retrieve neighborhood within {\tt distance} of {\tt nodes}.\\
{\tt \textbf{neighborhood}(graph::OptiGraph,nodes::Vector\{OptiNode\},distance::Int)}\\
\hline
Retrieve a subgraph from {\tt graph} including neighborhood nodes within {\tt distance} of {\tt sg}.\\
{\tt \textbf{expand}(graph::OptiGraph,sg::OptiGraph,distance::Int)}\\
\hline
\end{tabular}
\caption{Overview of core partitioning and topology functions in {\tt Plasmo.jl}.}
\label{table:partition_functions}
\end{center}
\end{table}

\subsubsection{Example 4: Partitioning a Dynamic Optimization Problem}
\begin{figure}
\begin{minipage}{\textwidth}
\begin{minipage}{\textwidth}
\begin{scriptsize}
\lstset{language=Julia,breaklines = true}
\begin{lstlisting}[caption = Construction of dynamic optimization problem \eqref{eq:simple_example},label = {code:ex4_construction}]
using Plasmo  |\label{line:ex4_load_packages_start}|
using Plots
using Ipopt

T = 100         #number of time points
d = sin.(1:T)   #disturbance   |\label{line:ex4_load_packages_end}|

#Create an OptiGraph
graph = OptiGraph()                 |\label{line:ex4_create}|

#Add nodes for states and controls
@optinode(graph,state[1:T])              |\label{line:ex4_state}|
@optinode(graph,control[1:T-1])          |\label{line:ex4_control}|

#Add state variables
for (i,node) in enumerate(state)   |\label{line:ex4_variables_start}|
    @variable(node,x)
    @constraint(node, x >= 0)
    @objective(node,Min,x^2)
end

#Add control variables
for node in control
    @variable(node,u)
    @constraint(node, u >= -1000)
    @objective(node,Min,u^2)
end

#Initial condition
n1 = state[1]
@constraint(n1,n1[:x] == 0)       |\label{line:ex4_variables_end}|

#Dynamic coupling
@linkconstraint(graph,[t = 1:T-1],state[t+1][:x] == state[t][:x] +   |\label{line:ex4_dynamics}|
                                    control[t][:u] + d[t])

#Optimize with Ipopt
ipopt = Ipopt.Optimizer
optimize!(graph,ipopt)                      |\label{line:ex4_solve}|

#Plot result structure
plt_graph4 = plot(graph,                     |\label{line:ex4_plt}|
layout_options = Dict(:tol => 0.1,:iterations => 500),
linealpha = 0.2,markersize = 6)

plt_matrix4 = spy(graph)                      |\label{line:ex4_spy}|
\end{lstlisting}
\end{scriptsize} %
\end{minipage}\\
\\
\begin{minipage}{0.5\textwidth}
\includegraphics[width=7.3cm,keepaspectratio]{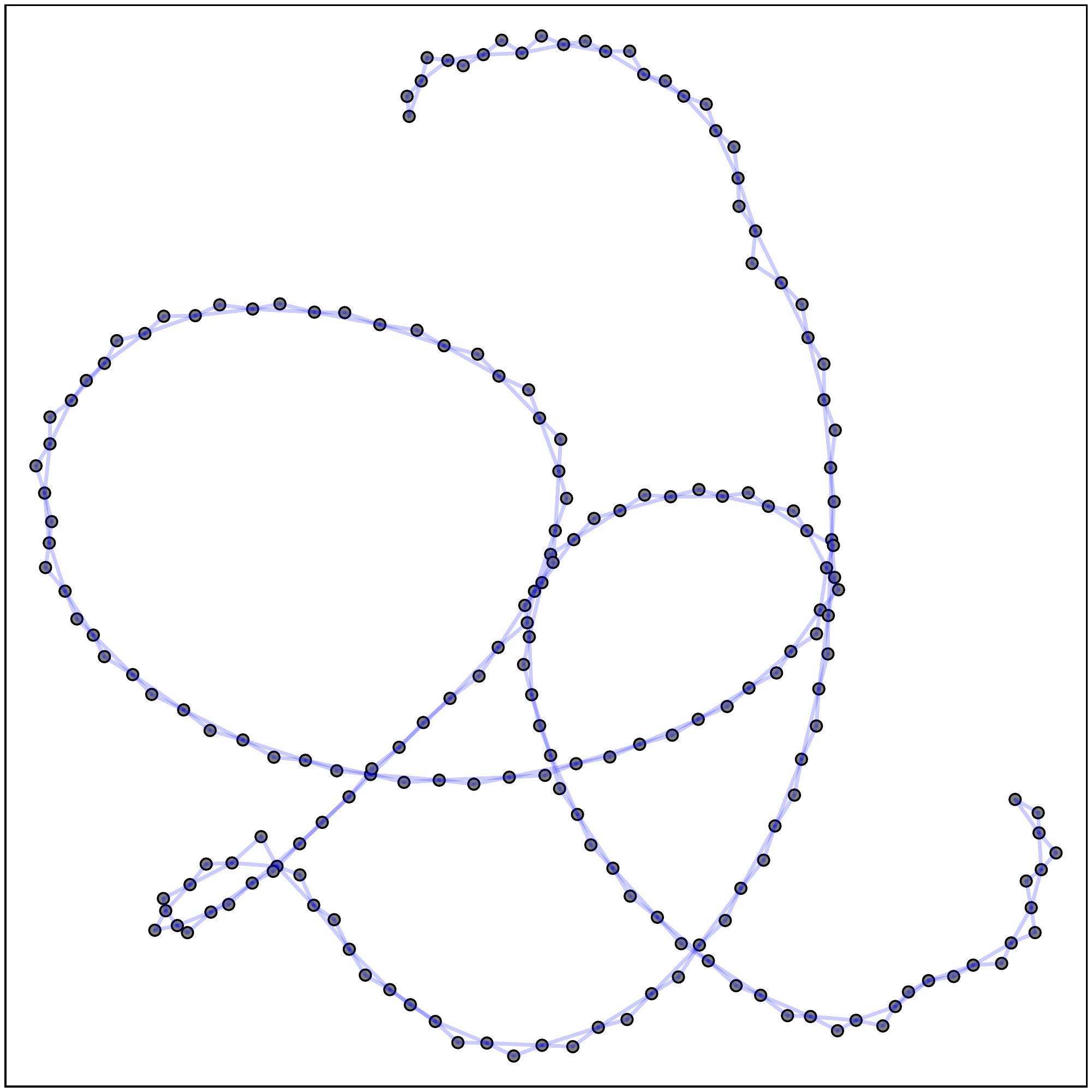}
\end{minipage}%
\begin{minipage}{0.5\textwidth}
\hspace{0.2cm}
\includegraphics[width=7.3cm,keepaspectratio]{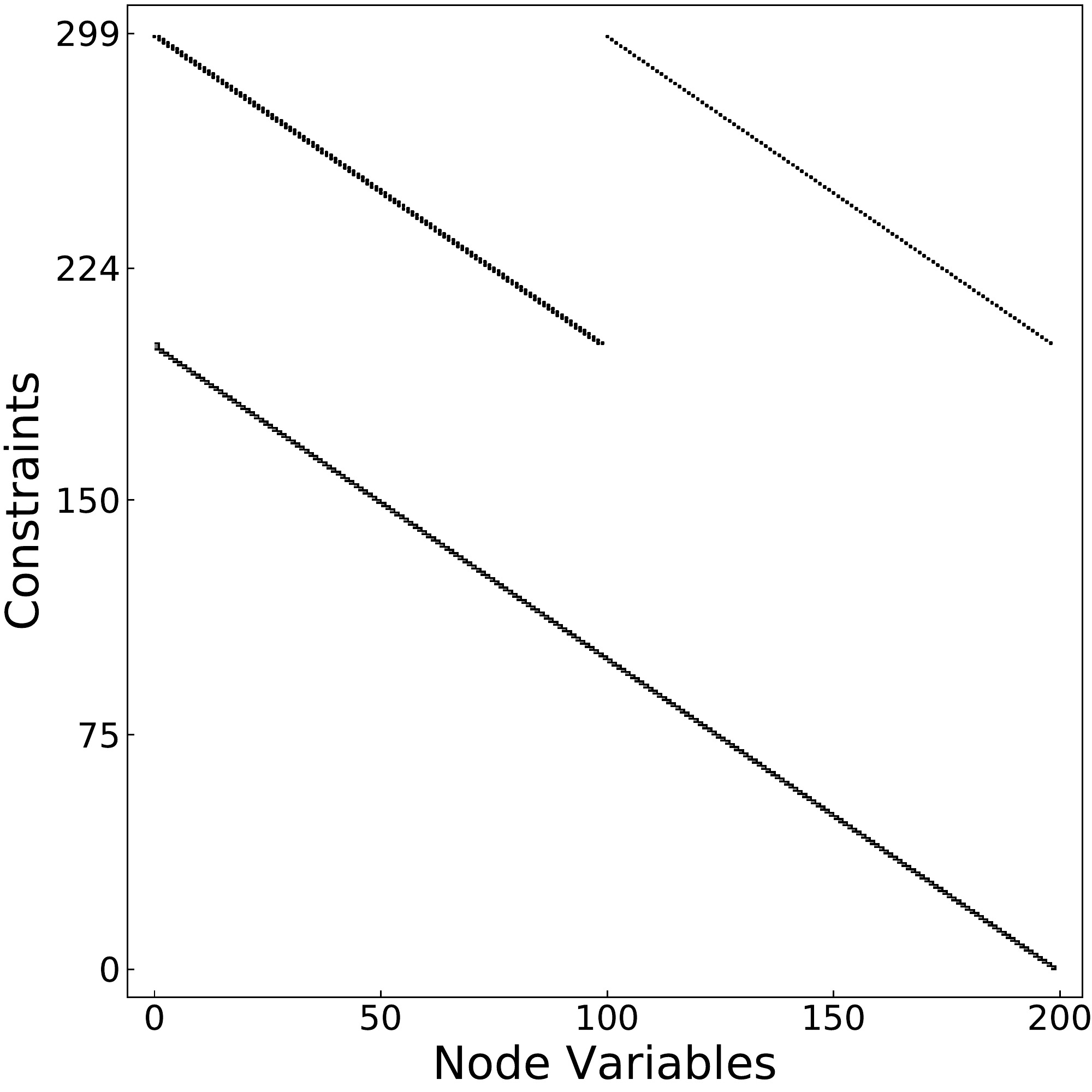}
\end{minipage}
\captionof{figure}{Output visuals for Code Snippet \ref{code:ex4_construction} showing graph structure of dynamic \\ optimization problem.}
\label{fig:ex4_plots_1}
\end{minipage}
\end{figure}

To demonstrate partitioning and manipulation capabilities, we use the simple dynamic optimization problem \cite{Shin2018}:

\begin{subequations}\label{eq:simple_example}
\begin{align}
    \min_{\{ x,u \}} & \sum_{t = 1}^T x_t^2 + u_t^2  &  \label{eq:example_objective}\\
    \textrm{s.t.} \quad & x_{t+1} = x_t + u_t + d_t, \quad t \in \{1,...,T-1\}  & \label{eq:example_temporal}\\
    & x_{1} = 0  &\label{eq:example_initial_condition}\\
    & x_t \ge 0, \quad t \in \{1,...,T\}\\
    & u_t \ge -1000, \quad t \in \{1,...,T-1\}
\end{align}
\end{subequations}
Here, $x$ is a vector of states and $u$ is a vector of control actions which are both
indexed over the set of time indices $t \in \{1,...,T\}$. Equation \eqref{eq:example_objective} minimizes the state trajectory with minimal control effort (energy), \eqref{eq:example_temporal} describes the
state dynamics, and \eqref{eq:example_initial_condition} is the initial condition.  This problem can be formulated using an {\tt OptiGraph} as shown in
Code Snippet \ref{code:ex4_construction} in much the same way as the examples in Section \ref{sec:modeling_with_subgraphs}.
We define the number of time
periods  $T = 100$ and create a disturbance vector {\tt d} (data) on Lines \ref{line:ex4_load_packages_start} through \ref{line:ex4_load_packages_end}.
In our implementation we create separate sets of nodes for the states and controls on Lines \ref{line:ex4_state} and \ref{line:ex4_control},
but it is also possible to define nodes for each individual time interval and add state and control variables to the resulting nodes.
Having control over this granularity is convenient for expressing what can be partitioned (i.e. features defined in an {\tt OptiNode} will remain
in the same partition).  Next we setup the state and control {\tt OptiNodes} on Lines \ref{line:ex4_variables_start} through \ref{line:ex4_variables_end}, we
use a linking constraint to capture dynamic coupling on Line \ref{line:ex4_dynamics} and we show how to solve the problem
with {\tt Ipopt} \cite{Wachter2006} on Line \ref{line:ex4_solve}.  We visualize the graph topology and matrix in Figure \ref{fig:ex4_plots_1}. The layouts depict a
linear graph but we note that the matrix has no obvious structure.

\begin{figure}
\begin{minipage}{\textwidth}
\begin{minipage}[]{\textwidth}
\begin{scriptsize}
\lstset{language=Julia,breaklines = true}
\begin{lstlisting}[caption = Creating subgraphs using hypergraph partitioning with {\tt KaHyPar},label = {code:ex4_kahypar}]
#Import KaHyPar interface
using KaHyPar  |\label{line:kaHyPar}|

#Create a hypergraph representation of the OptiGraph
hypergraph,reference_map = gethypergraph(graph)  |\label{line:hypergraph}|

#Perform hypergraph partitioning using KaHyPar
node_vector = KaHyPar.partition(hypergraph,8,imbalance = 0.1) |\label{line:kwayhyper}|

#Create a Partition object
partition = Partition(graph,node_vector,reference_map)    |\label{line:modelpartition}|

#Create subgraphs using the partition object and reference_map
make_subgraphs!(graph,partition) |\label{line:makesubgraphs}|

plt_graph5 = plot(graph,         |\label{line:ex4_plt2}|
layout_options = Dict(:tol => 0.01,:iterations => 500),
linealpha = 0.2,markersize = 6,subgraph_colors = true);
plt_matrix5 = spy(graph,subgraph_colors = true); |\label{line:ex4_spy2}|
\end{lstlisting}
\end{scriptsize}
\end{minipage}\\
\\
\begin{minipage}{0.5\textwidth}
\includegraphics[width=7.3cm,keepaspectratio]{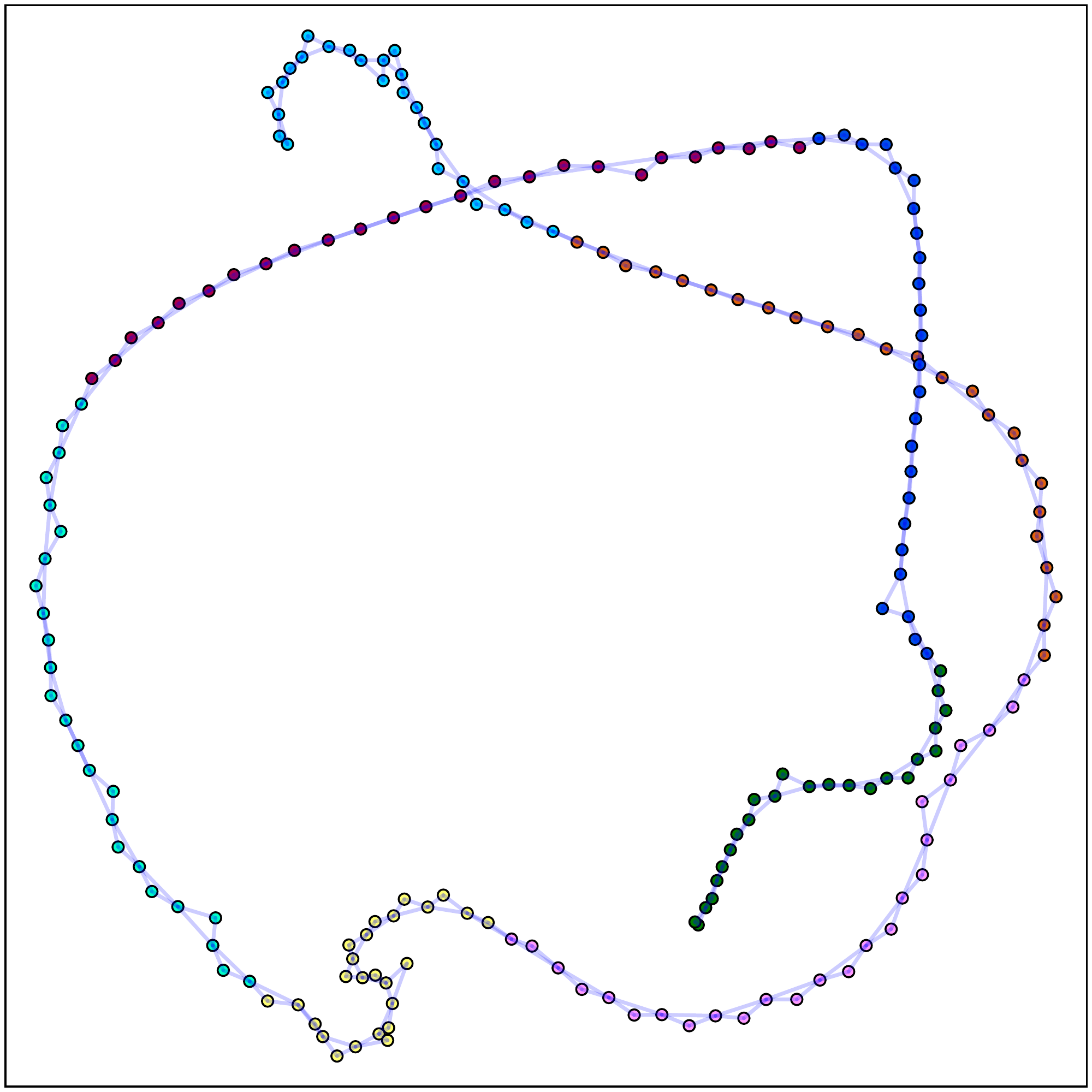}
\end{minipage}%
\begin{minipage}{0.5\textwidth}
\hspace{0.2cm}
\includegraphics[width=7.3cm,keepaspectratio]{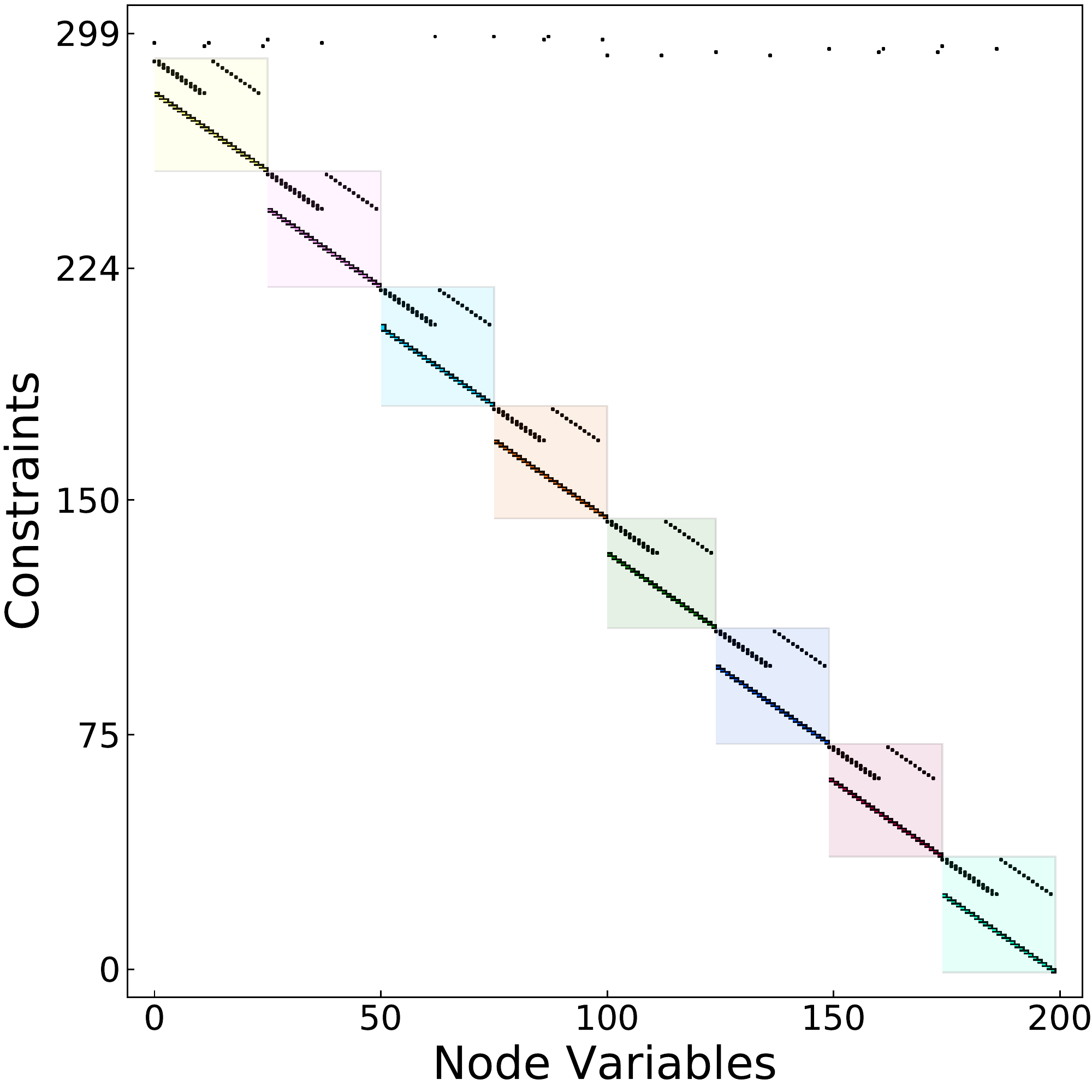}
\end{minipage}
\captionof{figure}{Output visuals for Code Snippet \ref{code:ex4_kahypar} showing partitions and reordering of dynamic \\ optimization problem.}
\label{fig:storage_problem_plots2}
\end{minipage}
\end{figure}

We partition the graph using {\tt KaHyPar}, as shown in Code Snippet \ref{code:ex4_kahypar};
here, Line \ref{line:kaHyPar} imports the {\tt KaHyPar} interface and Line \ref{line:hypergraph} creates the hypergraph representation corresponding to the
{\tt OptiGraph} using {\tt gethypergraph}.  We also return a {\tt reference\_map} which maps the hypergraph elements back to the {\tt OptiGraph}.
Line \ref{line:kwayhyper} performs hypergraph partitioning using {\tt KaHyPar} with a maximum imbalance ($\epsilon_{max}$) of 10\% and
Line \ref{line:modelpartition} creates a {\tt Partition} object using the resulting partition vector and the {\tt reference\_map}.
Line \ref{line:makesubgraphs} creates subgraphs in the graph {\tt graph} using the {\tt Partition} object
and the {\tt make\_subgraphs!} function.  We visualize the topology and matrix of the partitioned problem on
Lines \ref{line:ex4_plt2} and \ref{line:ex4_spy2} and these are shown in Figure \ref{fig:storage_problem_plots2}. This reveals eight distinct partitions and
their corresponding coupling. We note that the partitions are well-balanced and  note that the matrix is now rearranged into a banded structure that is typical of dynamic
optimization problems (partitioning automatically induces reordering).  {\tt Plasmo.jl} can use other graph representations to perform partitioning.  For instance,
we can create a traditional graph representation, as shown in Code Snippet \ref{code:ex4_metis}, and partition it with {\tt Metis}. We then use the {\tt reference\_map} to
obtain the original {\tt OptiGraph} elements to create a graph {\tt Partition}. We could also partition less intuitive representations (e.g., bipartite) in this way; we only require a mapping
from the partition back to the {\tt OptiGraph} elements. The partitioning procedure shown here can also be repeated to create an arbitrary number of subgraph levels.

\begin{figure}
\begin{minipage}{\textwidth}
\begin{scriptsize}
\lstset{language=Julia,breaklines = true}
\begin{lstlisting}[caption = Creating a subgraph partition using graph partitioning with {\tt Metis},label = {code:ex4_metis}]
#Import the Metis interface
using Metis

#Retrieve underlying hypergraph from dynamic opt problem graph
simple_graph,reference_map = getcliquegraph(graph)

#Run Metis direct k-way partitioning with a 8 maximum of partitions
node_vector = Metis.partition(simple_graph,8,alg = :KWAY)

#Create a Partition object using node vector and reference_map
partition = Partition(graph,node_vector,reference_map)

#Create subgraphs using Partition object
make_subgraphs!(graph,partition)
\end{lstlisting}
\end{scriptsize}
\end{minipage}
\end{figure}

\begin{figure}
\begin{minipage}{\textwidth}
\begin{minipage}[]{\textwidth}
\begin{scriptsize}
\lstset{language=Julia,breaklines = true}
\begin{lstlisting}[caption = Aggregating nodes in an OptiGraph,label = {code:combine_snippet}]
#Combine the subgraphs in graph into OptiNodes in a new OptiGraph
aggregated_graph,ref_map = aggregate(graph,0)                    |\label{line:combine_subgraphs}|

#plot the graph a matrix layouts of the aggregated OptiGraph
plt_graph6 = plot(aggregated_graph,                    |\label{line:ex4_plt3}|
layout_options = Dict(:tol => 0.01,:iterations => 10),
node_labels = true,markersize = 30,labelsize = 20,node_colors = true)
plt_matrix6 = spy(aggregated_graph,node_labels = true,node_colors = true)   |\label{line:ex4_spy3}|
\end{lstlisting}
\end{scriptsize} %
\end{minipage}\\
\\
\begin{minipage}{0.5\textwidth}
\includegraphics[width=7.3cm,keepaspectratio]{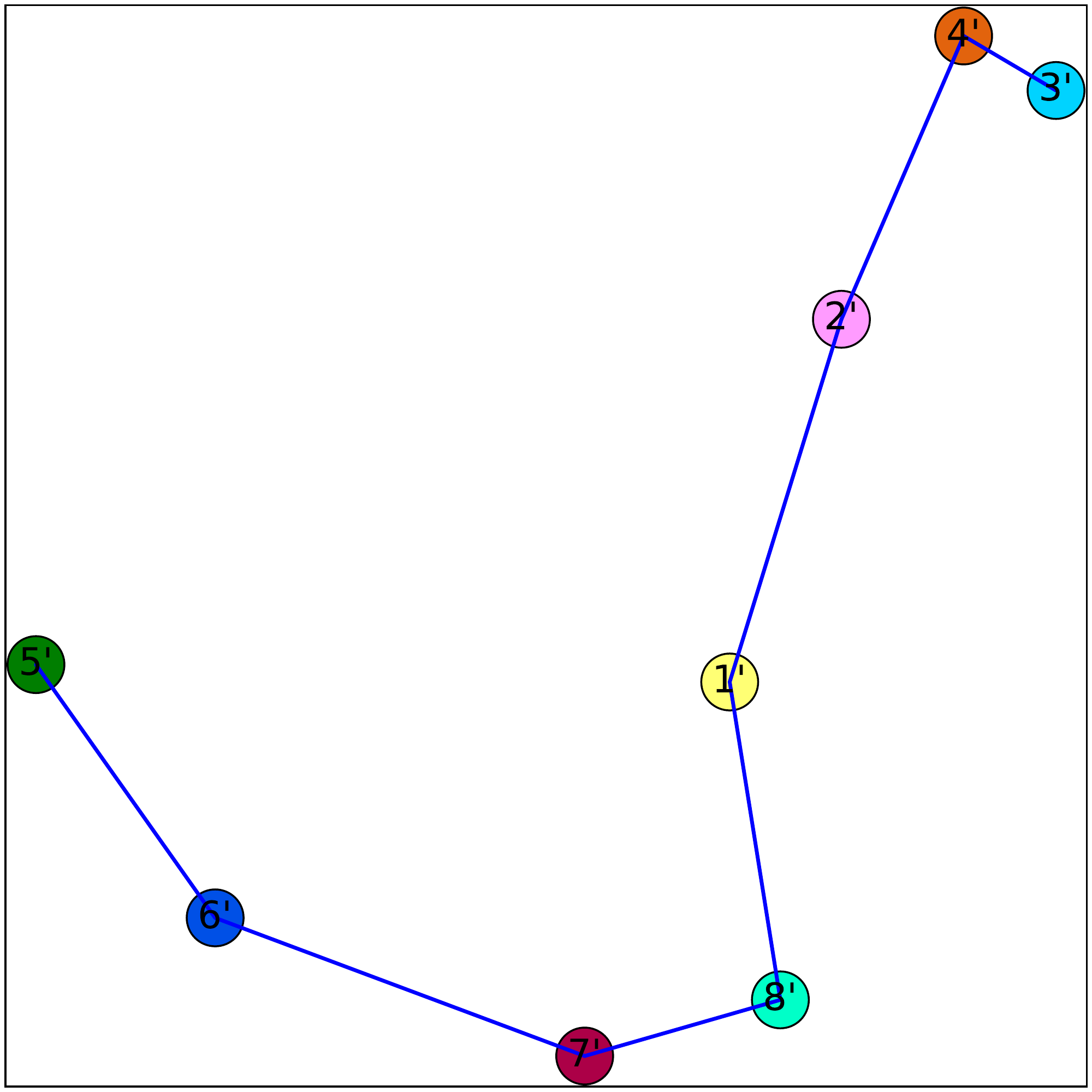}
\end{minipage}%
\begin{minipage}{0.5\textwidth}
\hspace{0.2cm}
\includegraphics[width=7.3cm,keepaspectratio]{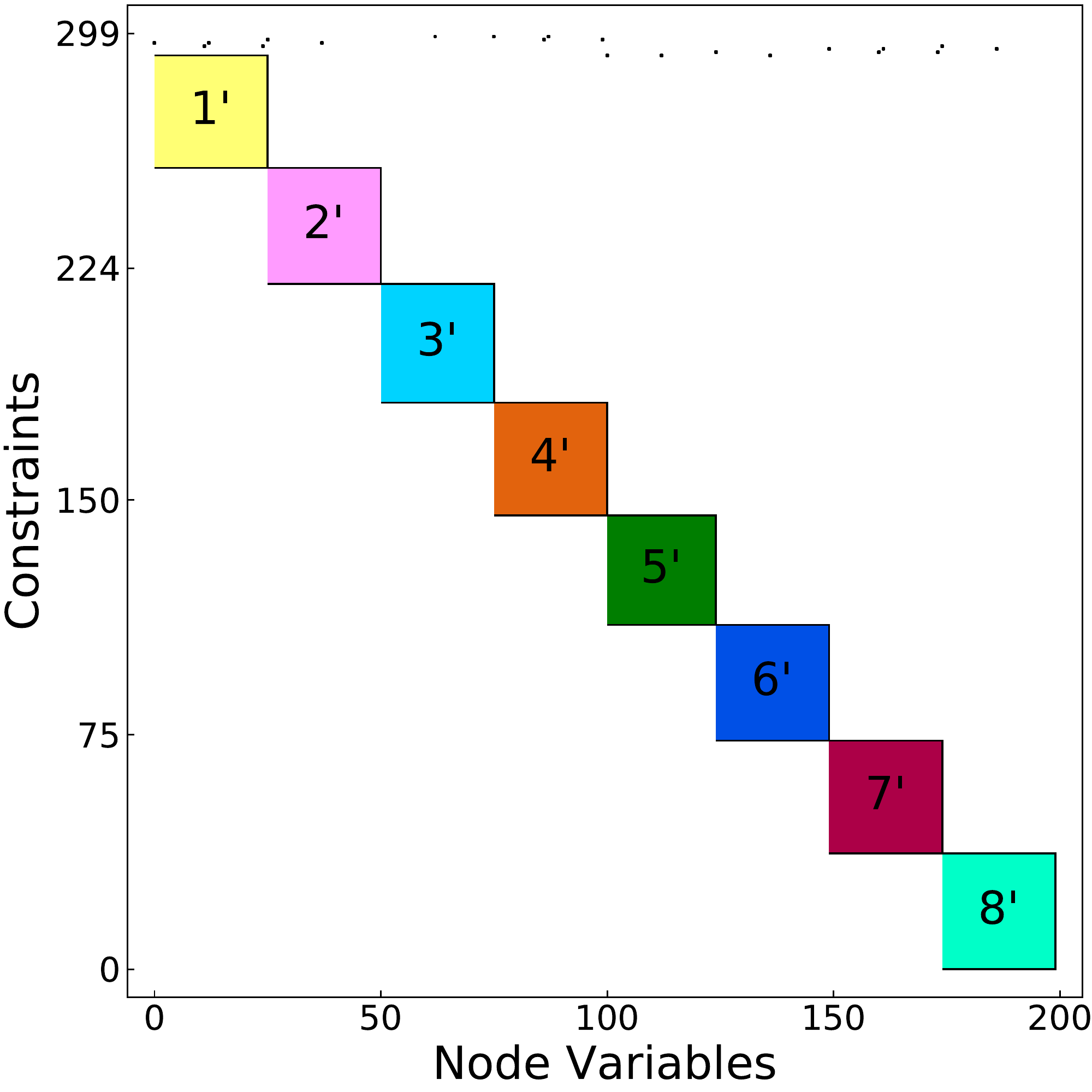}
\end{minipage}
\captionof{figure}{Output visuals for Code Snippet \ref{code:combine_snippet} showing aggregated graph of dynamic \\ optimization problem.}
\label{fig:ex4_plots3}
\end{minipage}
\end{figure}

The {\tt OptiNodes} in a graph can be \emph{aggregated} to form larger nodes, as shown in Figure \ref{fig:modelgraph_combine}.
This aggregation functionality is key to communicate subproblems to decomposition algorithms that operate at different levels of granularity. Aggregation can also be performed to
collapse an entire graph into a single node, producing a standard optimization problem that can be solved with off-the-shelf solvers. Code Snippet \ref{code:combine_snippet} shows how to
aggregate the graph of the dynamic optimization problem into a new aggregated graph with eight nodes.  We create the new ({\tt aggregated\_graph}) on Line \ref{line:combine_subgraphs}
by using {\tt aggregate} on the {\tt graph} from Snippet \ref{code:ex4_kahypar}.  We provide the integer {\tt 0} to the function to
specify that we want zero subgraph levels which converts the eight subgraphs into nodes. For hierarchical graphs with many levels,
we can define how many subgraph levels we wish to retain. We plot the graph and matrix layouts for the aggregated {\tt OptiGraph} on
Lines \ref{line:ex4_plt3} and \ref{line:ex4_spy3} and these are shown in Figure \ref{fig:ex4_plots3}.

\subsubsection{Example 5: Using Graph Topology Functions}
\begin{figure}
\begin{minipage}{\textwidth}
\begin{minipage}[]{\textwidth}
\begin{scriptsize}
\lstset{language=Julia,breaklines = true}
\begin{lstlisting}[caption = Using Graph Topology Functions,label = {code:expand_snippet}]
#Get the current subgraphs in graph
subgraphs = getsubgraphs(graph)             |\label{line:ex5_getsubgraphs}|

#Query the first subgraph
sg1 = subgraphs[1]

#Query the edges inccident to the nodes in sg1
incident_edges = incident_edges(graph,all_nodes(sg1))       |\label{line:ex5_incident_edges}|

distance = 2
#Query the nodes within distance 2 of sg1
nodes = neighborborhood(graph,all_nodes(sg1),distance)        |\label{line:ex5_neighborhood}|

#Broadcast expand function to each subgraph.
#Obtain vector of expanded subgraphs
expanded_subgraphs = expand.(Ref(graph),subgraphs,Ref(distance))     |\label{line:ex5_expand}|

#Plot the expanded subgraphs
plt_graph7 = plot(graph,expanded_subgraphs,           |\label{line:ex5_plt_expand}|
layout_options = Dict(:tol => 0.01,:iterations => 1000),
markersize = 6,linealpha = 0.2)
plt_matrix7 = spy(graph,expanded_subgraphs)    |\label{line:ex5_spy_expand}|
\end{lstlisting}
\end{scriptsize} %
\end{minipage}\\
\\
\begin{minipage}{0.5\textwidth}
\includegraphics[width=7.3cm,keepaspectratio]{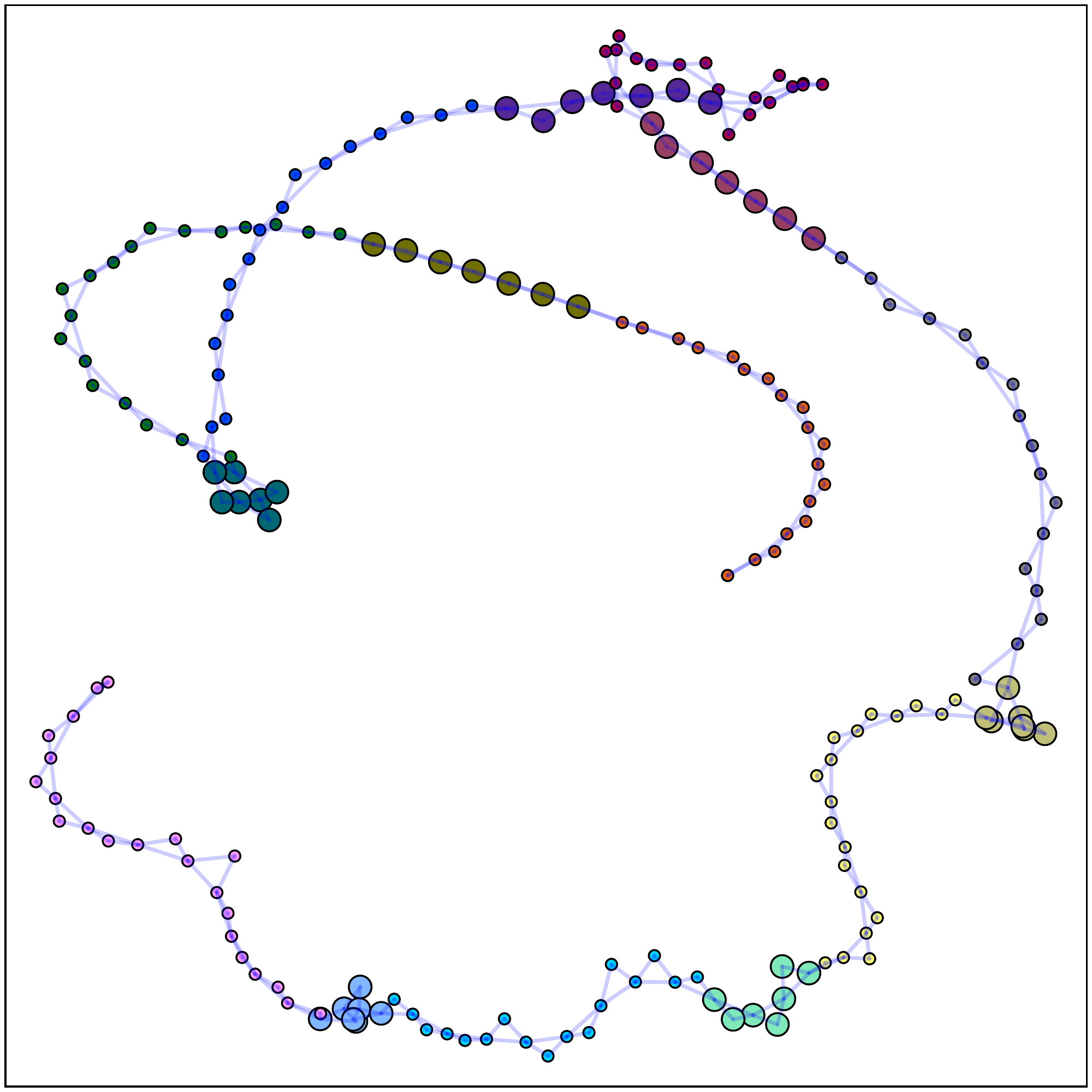}
\end{minipage}%
\begin{minipage}{0.5\textwidth}
\hspace{0.2cm}
\includegraphics[width=7.3cm,keepaspectratio]{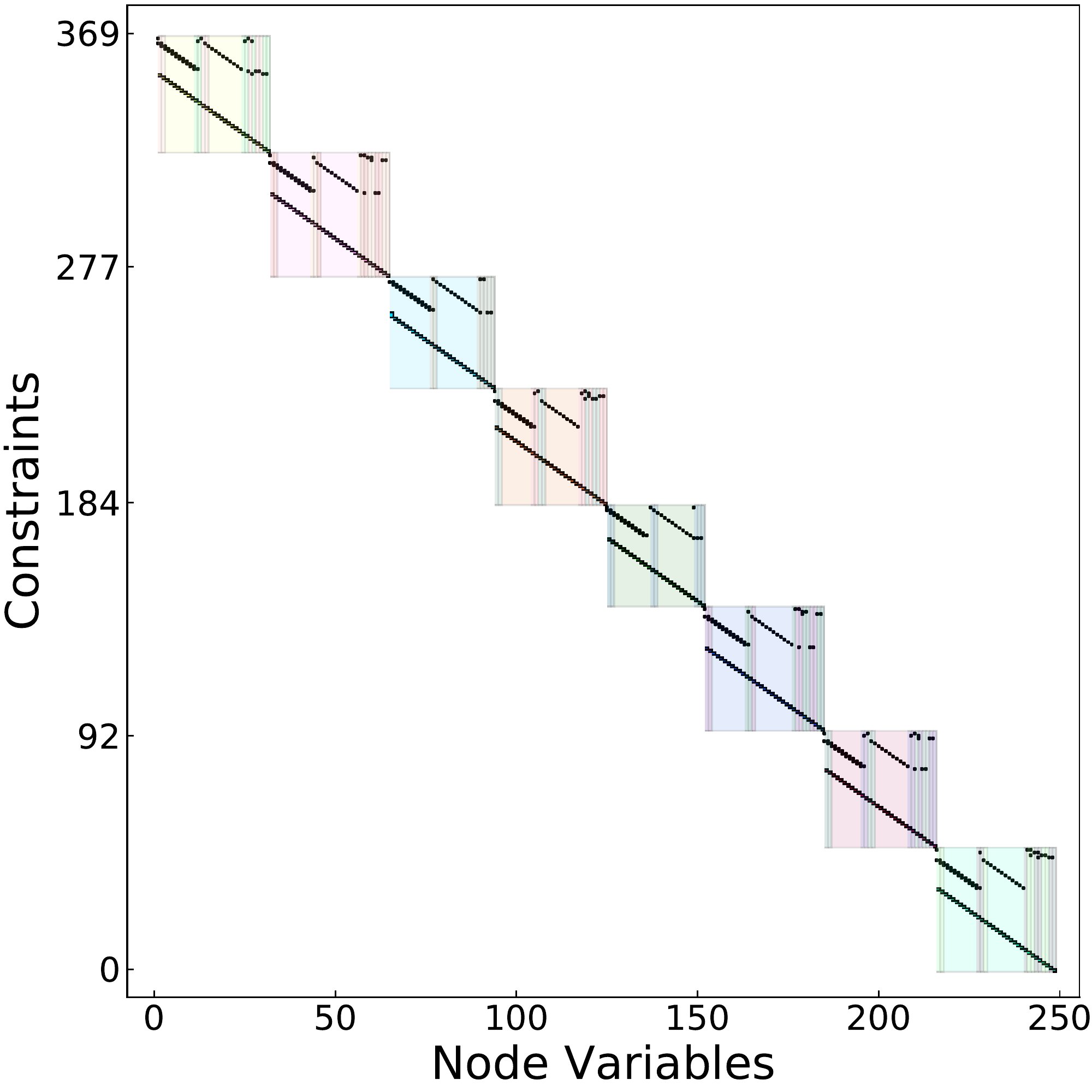}
\end{minipage}
\captionof{figure}{Output visuals for Code Snippet \ref{code:expand_snippet} showing overlapping subgraphs.}
\label{fig:ex5_plots}
\end{minipage}
\end{figure}

We now show how to use graph topology functions by computing {\em overlapping} partitions for the dynamic optimization example. This is shown in
Code Snippet \ref{code:expand_snippet}; here, Line \ref{line:ex5_getsubgraphs} obtains subgraphs from the {\tt OptiGraph} {\tt graph}
created in Snippet \ref{code:ex4_kahypar}, Line
\ref{line:ex5_incident_edges} determines and returns the {\tt OptiEdges} incident to the {\tt OptiNodes} in the first subgraph {\tt sg1}, and Line
\ref{line:ex5_neighborhood} returns the complete neighborhood around the same {\tt OptiNodes} within a distance of two.  On Line \ref{line:ex5_expand}, we broadcast the
{\tt expand} function (using dot syntax {\tt expand.()} and {\tt Ref} as typical in {\tt Julia}) to create a new set of subgraphs, each expanded by a distance of two.
Line \ref{line:ex5_plt_expand} plots the layout of {\tt graph} with the expanded subgraphs as is shown in Figure \ref{fig:ex5_plots}.  Here, we see eight
distinct partitions (each with a unique color) where the larger markers represent nodes that are part of multiple subgraphs (they are also the average color of their
containing subgraphs). Figure \ref{fig:ex5_plots} shows the corresponding matrix layout where highlighted columns indicate that the
node also appears in other subgraph blocks. Overlapping partitions are useful in certain algorithms such as Schwarz decomposition.

\FloatBarrier

\section{Algorithms}\label{sec:algorithm_interfaces}
The {\tt OptiGraph} provides a flexible abstraction that facilitates communication of problem structure to a wide range of decomposition strategies.
The structure can be exploited at a problem level, in which the decomposition strategy treats {\tt OptiNodes} as optimization subproblems whose solutions are coordinated to find
a solution of the entire {\tt OptiGraph}. This strategy is used, for example, in Benders, Lagrangian dual, and ADMM decomposition. The structure can also be exploited at the linear
algebra level, in which a decomposition strategy is used within a general algorithm to compute a search step. Here, the strategy treats {\tt OptiNodes} as linear systems that result from
the optimality conditions of the associated subproblems and coordinates their solutions to find a solution of the linear system associated with the optimality conditions of the
entire {\tt OptiGraph}.  In this section we show how to use the proposed abstraction and {\tt Plasmo.jl} functionality to explain structures at the linear algebra and problem level.

\begin{figure}[]
\centering
\begin{subfigure}[t]{0.48\textwidth}
    \centering
    \includegraphics[width=4cm, height=4cm]{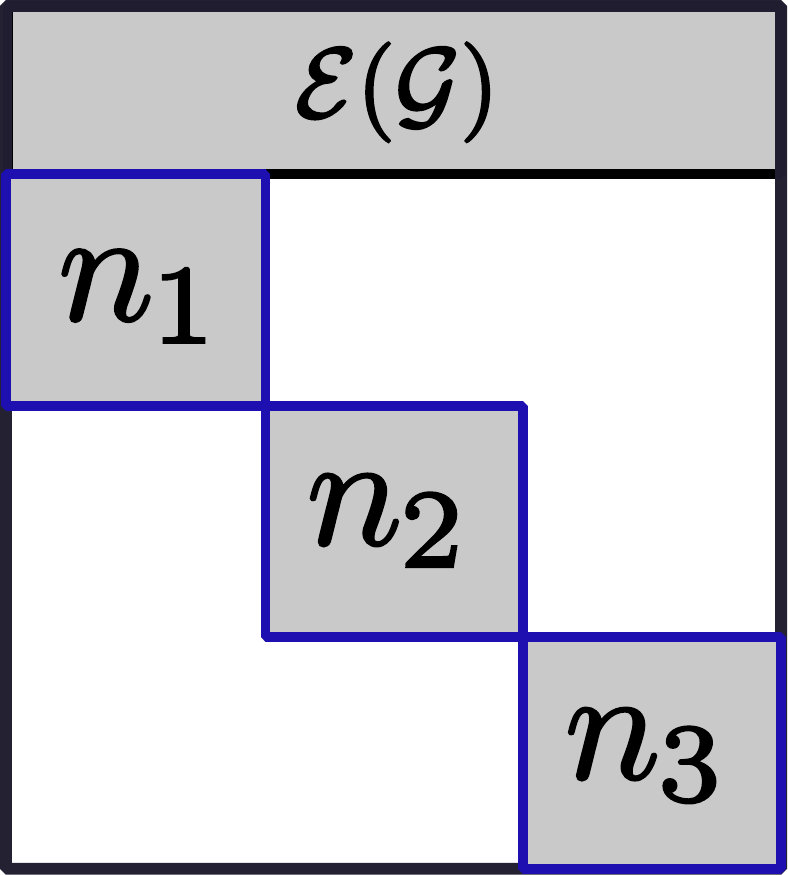}
    \label{fig:single_graph}
\end{subfigure}
\begin{subfigure}[t]{0.48\textwidth}
    \centering
    \includegraphics[width=5cm, height=5cm]{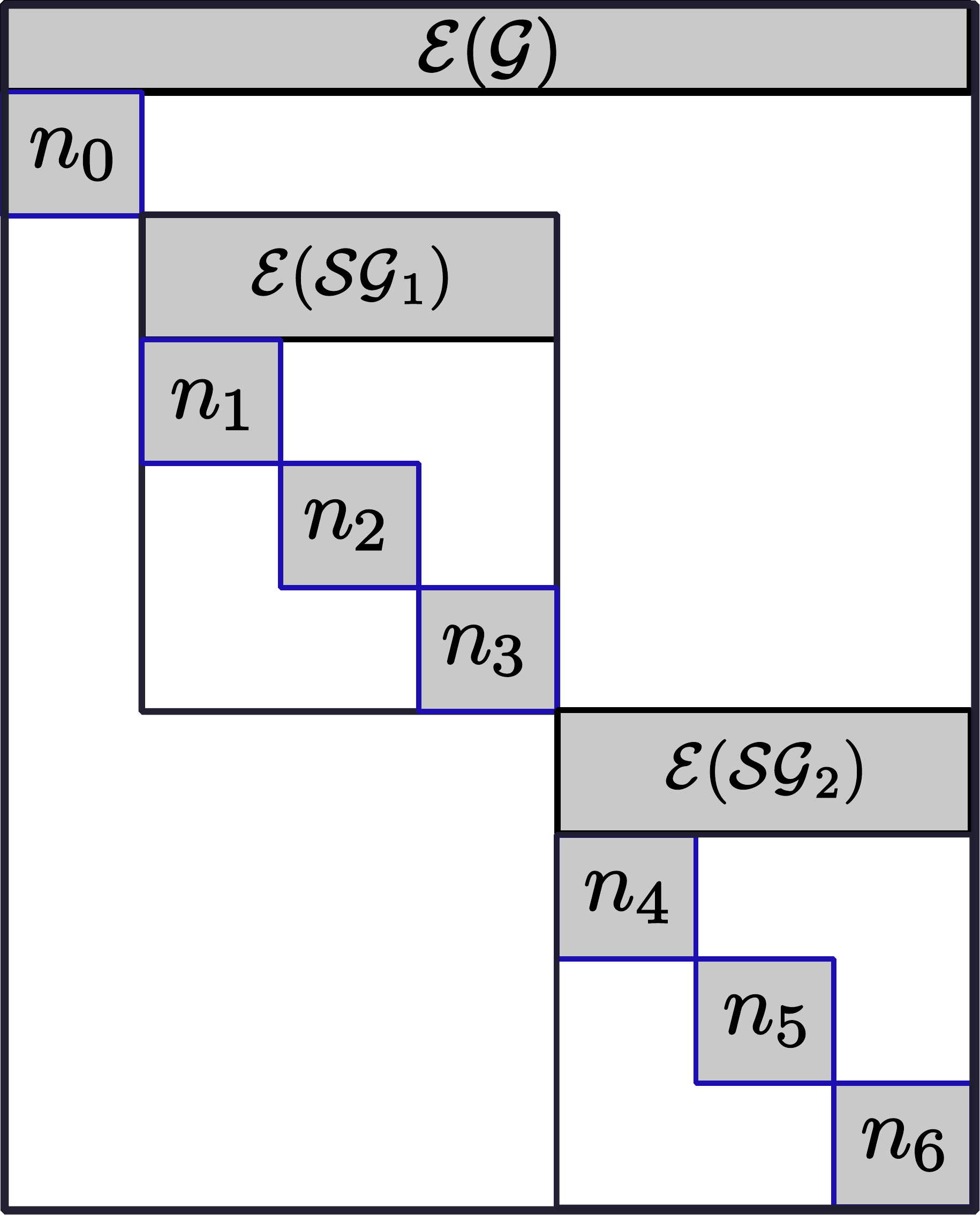}
    \label{fig:nested_graphs}
\end{subfigure}
\caption{Block structure (left) and nested block structure (right) of an {\tt OptiGraph}.}\label{fig:block_structured_problems}
\end{figure}

\subsection{Linear Algebra Decomposition}
It is well-known that block structures can be exploited by linear algebra operations within interior-point
algorithms. For instance,
continuous problems with the partially-separable structure described by \eqref{eq:model-graph-compact}
induce structured linear algebra kernels that can be solved using tailored techniques such as Schur decomposition.
We express the continuous variant of the graph formulation \eqref{eq:model-graph-compact} with feasible sets of the form $\mathcal{X}_n:=\{x\,|\,c_n(x)=0\}$ which gives rise to the
Karush-Kuhn-Tucker (KKT) system given by:

\begin{subequations}\label{eq:graph-KKT}
    \begin{align}
        & \sum_{n \in \mathcal{N}(\mathcal{G})} \Big( \nabla_{x_n} f_n(x_n) + \nabla_{x_n} c_n(x_n) \lambda_n \Big ) +
        \sum_{e \in \mathcal{E}(\mathcal{G})} \nabla_{\{x_n\}_{n \in \mathcal{N}(e)}} g_e(\{x_n\}_{n \in \mathcal{N}(e)}) \lambda_e = 0,\\
        & c_n(x_n) = 0, \quad  n \in \mathcal{N}(\mathcal{G}),\\
        & g_e(\{x_n\}_{n \in \mathcal{N}(e)}) = 0,  \quad e \in \mathcal{E}(\mathcal{G}).
    \end{align}
\end{subequations}
Here, we recall that $\mathcal{N}(e)$ denotes the nodes that
support edge $e$. The terms associated with barrier functions exist in \eqref{eq:graph-KKT} in the presence of inequality constraints.
Upon linearization of \eqref{eq:graph-KKT}, the system of algebraic equations gives rise to the block-bordered KKT system given by:

\begin{align}\label{eq:bordered_KKT}
\left[\begin{array}{cccc|c}
K_{n_1}& &&& B_{n_1}\\
&K_{n_2}& && B_{n_2}\\
&&\ddots&&\vdots\\
&&&K_{n_{|\mathcal{N}|}}&B_{n_{|\mathcal{N}|}}\\\hline
B_{n_1}^T&B_{n_2}^T&\hdots&B_{n_{|\mathcal{N}|}}^T & \\
\end{array}\right]
\left[\begin{array}{c}\Delta w_{n_1}\\ \Delta w_{n_2}\\ \vdots \\ \Delta w_{n_{|\mathcal{N}|}}\\ \hline \Delta \lambda_{\mathcal{E}(\mathcal{G})} \end{array}\right]=
-\left[\begin{array}{c} r_{n_1}\\ r_{n_2}\\  \vdots \\ r_{n_{|\mathcal{N}|}}\\ \hline \ r_{\mathcal{E}(\mathcal{G})}\ \end{array}\right].
\end{align}
Here, $\Delta w_n := (\Delta x_n, \Delta \lambda_n)$ is the primal-dual step for the variables and constraints on node $n$ and
$\Delta \lambda_{\mathcal{E}(\mathcal{G})}$ is a vector of steps corresponding to the dual variables on the {\tt OptiEdges} in {\tt OptiGraph} $\mathcal{G}$,
where $\Delta \lambda_{\mathcal{E}(\mathcal{G})} := \{\Delta \lambda_e\}_{e \in \mathcal{E}(\mathcal{G})}$. We also define
\begin{align*}
K_n := \left[\begin{array}{cc}
W_n& J_n^T\\
J_n & 0
\end{array}\right],
\end{align*}
as the block matrix corresponding to {\tt OptiNode} $n$ where $W_n = \nabla_{x_n,x_n}\mathcal{L}$ is the Hessian of the Lagrange
function of \eqref{eq:model-graph-compact} and $J_n := \nabla_{x_n}c_n(x_n)$ is the constraint Jacobian.  We define coupling blocks $B_n$ as
\begin{align*}
B_n := \left[\begin{array}{cc}
Q_n & 0\\
0 & 0
\end{array}\right],
\end{align*}
where $Q_n := \nabla_{x_n} \{g_e(\{x_{n^\prime}\}_{n^\prime \in \mathcal{E}(n)})\}_{e \in \mathcal{E}(\mathcal{G})}$.
If all of the linking constraints in the graph are \emph{linear}, $Q_n$ reduces to $\Pi^T_n$ where $\Pi_n$ is the matrix
of coefficients corresponding to the {\tt OptiEdges} incident to node $n$.

\subsubsection{Schur Decomposition}
The block-bordered structure described in \eqref{eq:graph-KKT} can be exploited using Schur decomposition or
block preconditioning strategies. The typical Schur decomposition algorithm is described by \eqref{eq:schur_decomposition} in terms of the {\tt OptiGraph} abstraction,
where $S$ is the Schur complement matrix.

\begin{subequations}\label{eq:schur_decomposition}
    \begin{align}
        & S = -\sum_{n \in \mathcal{N}(\mathcal{G})} B_n^T K_n^{-1}B_n \label{eq:form_schur_complement}\\
        & S \Delta \lambda_{\mathcal{E}(\mathcal{G})} = \sum_{n \in \mathcal{N}(\mathcal{G})} B_n^T K_n^{-1}r_n -  r_{\mathcal{E}(\mathcal{G})}  \label{eq:schur_first_stage}\\
        & K_n \Delta w_n = B_n \Delta \lambda_{\mathcal{E}(\mathcal{G})} - r_n \label{eq:schur_second_stage}
    \end{align}
\end{subequations}

Step \eqref{eq:form_schur_complement} requires factorizing the linear system associated to each {\tt OptiNode} ($K_n$)
and computes the Schur complement contribution $B_n^T K_n^{-1} B_n$ on each node (possibly in parallel).
After each contribution is computed (each requires performing a sparse factorization), the total Schur complement matrix $S$ is
created and factorized to solve the linear system \eqref{eq:schur_first_stage} (in serial) and take a step in the {\tt OptiEdge} dual variables
($\lambda_{\mathcal{E}(\mathcal{G})}$).  Step \eqref{eq:schur_second_stage} solves for
the {\tt OptiNode} primal-dual step $\Delta w_n$ given the {\tt OptiEdge} dual step (also possibly in parallel).

The general Schur decomposition exhibits a couple of major {\em computational bottlenecks}.  Forming the contributions $B^T_n K^{-1}_n B_n$
is expensive when there are many columns in $B_n$ (the number of columns corresponds to the number linking constraints incident to node $n$).  Moreover, if
the node blocks have different sizes, the factorization of the blocks $K_n$ can create load imbalance and memory issues.
Second, factorizing the Schur matrix $S$ is challenging when there are many linking constraints because this matrix is dense or is composed of dense sub-blocks.
Consequently, one needs to control the amount of coupling between the blocks. The {\tt OptiGraph} topology directly corresponds with the size of the block matrices that appear in the
Schur decomposition and thus can be manipulated to facilitate computational efficiency. Specifically, we wish to manipulate the partitioning of an {\tt OptiGraph} to accelerate Schur decomposition.
We use the capabilities of {\tt Plasmo.jl} to experiment with different partitioning strategies and with this analyze trade-offs between coupling, imbalance,
and memory use. Section \ref{sec:case_study_pips} details
these performance trade-offs with a large example.

\subsubsection{PIPS-NLP Interface}

{\tt Plasmo.jl} includes an interface to the structure-exploiting interior-point solver {\tt PIPS-NLP} that we call \\ {\tt PipsSolver.jl}.
{\tt PIPS-NLP} provides a general Schur decomposition strategy that performs parallel computations via {\tt MPI}.  {\tt PIPS-NLP} was developed to solve large-scale stochastic programming problems
that adhere to a two-level tree structure (i.e. a single first stage problem coupled to second stage subproblems).
The {\tt OptiGraph} interface to {\tt PIPS-NLP} converts general graph structures into this format (e.g., via aggregation and partitioning).  For instance, using
our example problem \eqref{eq:simple_example}, we can communicate the problem structure to {\tt PIPS-NLP}
and solve the problem in parallel. This is shown in Snippet \ref{code:pipsnlp_problem_snippet}.

\begin{minipage}[]{0.9\linewidth}
\begin{scriptsize}
\lstset{language=Julia,breaklines = true}
\begin{lstlisting}[caption = Solving an {\tt OptiGraph} model in parallel with PIPS-NLP , label={code:pipsnlp_problem_snippet}]
using Distributed             | \label{line:use_distributed} |
using MPIClusterManagers      | \label{line:use_MPIClusterManagers} |

# specify, number of mpi workers
manager=MPIManager(np=2)       | \label{line:mpi_manager} |
# start mpi workers and add them as julia workers too.
addprocs(manager)               | \label{line:addprocs} |

#Setup the worker environments
@everywhere using Plasmo       | \label{line:using_modelgraphs_distributed} |
@everywhere using PipsSolver  #Solver interface to PIPS-NLP      | \label{line:usingMGpipssolver} |

#get the julia ids of the mpi workers
julia_workers = collect(values(manager.mpi2j))  | \label{line:juliaworkers} |

#Use the pips-nlp interface to distribute the OptiGraph among workers
#Here, we create the variable `pips_graph` on each worker
remote_references = PipsSolver.distribute(graph,julia_workers,    | \label{line:distributegraph} |
remote_name = :pips_graph)

#Solve with PIPS-NLP
@mpi_do manager begin     | \label{line:solve_pips} |
    using MPI
    PipsSolver.pipsnlp_solve(pips_graph)
end
\end{lstlisting}
\end{scriptsize}
\end{minipage}

Snippet \ref{code:pipsnlp_problem_snippet} presents a standard template for setting up distributed computing environments to use {\tt PipsSolver.jl}
which is worth discussing. Line \ref{line:use_distributed} imports the {\tt Distributed} module, which is a {\tt Julia} package for performing distributed computing.
On Line \ref{line:use_MPIClusterManagers} we import the {\tt MPIClusterManagers} package which allows us to interface {\tt MPI} ranks
(used by {\tt PIPS-NLP}) with {\tt Julia} worker CPUs.
We create a {\tt manager} object and specify that we want to use 2 workers on Line \ref{line:mpi_manager} and add the Julia workers on Line \ref{line:addprocs}.
Next we setup the model and solver environments for the added workers on Lines \ref{line:using_modelgraphs_distributed} and \ref{line:usingMGpipssolver} and create a reference to the
julia workers by querying the manager on Line \ref{line:juliaworkers}.  We \emph{distribute} the graph among the workers in Line \ref{line:distributegraph}
using the function provided by {\tt PipsSolver}. This function sets up the relevant graph nodes on each worker and creates the graph named {\tt pips\_graph} on each worker
(internally this function inspects the {\tt OptiGraph} and allocates {\tt OptiNodes} to worker CPUs).
Finally, we use the {\tt \@mpi\_do} function from {\tt MPIClusterManagers} to execute {\tt MPI} on each worker and solve the graph.
Each worker executes the {\tt pipnlp\_solve} function
and communicates using {\tt MPI} routines within {\tt PIPS-NLP}.

\subsection{Problem Level Decomposition}

Overlapping Schwarz is a flexible graph-based decomposition strategy that can be used at the linear algebra or problem level \cite{Fromeer2002,Shin2020}. This approach has been used to solve dynamic optimization \cite{Shin2019} and network optimization problems \cite{Shin2020b}. At the problem level, the Schwarz algorithm decomposes the problem graph over {\em overlapping} partitions. The algorithm solves subproblems over the overlapping partitions and coordinates solutions by exchanging primal-dual information over the overlapping regions. The presence of overlap is key to promote convergence; in particular, it has been proven that the convergence rate improves exponentially with the size of the overlap \cite{Shin2020}. The Schwarz scheme is also flexible in that the overlap can be adjusted to trade-off subproblem complexity (time and memory) with convergence speed. When the overlap is the entire graph, the algorithm solves the entire graph once (converges in one iteration).  When the overlap is zero, the algorithm operates as a standard Gauss-Seidel scheme and will exhibit slow convergence (or no convergence at all). We thus have that Schwarz has fully centralized and fully decentralized schemes as extremes.  Schwarz algorithms can also be implemented under synchronous and asynchronous settings to handle load imbalance issues \cite{Shin2020}.

\subsubsection{Development of Schwarz Algorithm}
The Schwarz algorithm iteratively solves subproblems associated with overlapping subgraphs. In particular, one first
constructs {\it expanded subgraphs} $\{\mathcal{SG}^\prime_i\}_{i=1}^N$ from the subgraphs $\{\mathcal{SG}_i\}_{i=1}^N$ obtained from {\tt OptiGraph} partitioning. This procedure can be performed by expanding the subgraphs and adding {\tt OptiNodes} within a prescribed distance $\omega\geq 0$. The optimization subproblems for the expanded subgraph $\mathcal{SG}_i'$ can be formulated as:
\begin{subequations}\label{eq:Schwarz_subproblem}
    \begin{align}
        \min_{{\{x_n}\}_{n \in \mathcal{N}(\mathcal{SG}_i^\prime)}} & \quad \sum_{n \in \mathcal{N}(\mathcal{SG}_i^\prime)} f_n(x_n) -
        \sum_{e \in \mathcal{I}_1(\mathcal{SG}_i^\prime)}(\lambda_{e}^k)^T g_{e}(\{x_{n}\}_{n \in \mathcal{N}(e)\cap\mathcal{N}(\mathcal{SG}^\prime_i)},
        \{x_{n}^k\}_{n \in \mathcal{N}(e)\setminus \mathcal{N}(\mathcal{SG}'_i)}) \label{eq:Schwarz-objective} \\
        \textrm{s.t.} & \quad x_n \in \mathcal{X}_n,  \quad n \in \mathcal{N}(\mathcal{SG}_i^\prime), \label{eq:Schwarz-constraints} \\
        & \quad g_e(\{x_n\}_{n \in \mathcal{N}(e)}) = 0, \quad e \in \mathcal{E}(\mathcal{SG}^\prime_i), \label{eq:Schwarz-local-linkconstraints}\\
        & \quad g_e(\{x_n\}_{n \in \mathcal{N}(e)\cap\mathcal{N}(\mathcal{SG}^\prime_i)},\{x_n^k\}_{n \in \mathcal{N}(e)\setminus \mathcal{N}(\mathcal{SG}^\prime_i)}) = 0,
        \quad e \in \mathcal{I}_2(\mathcal{SG}^\prime_i).     \label{eq:primal-received}
    \end{align}
\end{subequations}
Here, the dual variables for \eqref{eq:Schwarz-local-linkconstraints} and \eqref{eq:primal-received} are denoted by $\lambda_e$, $\mathcal{N}(\mathcal{SG}'_i)$ is
the set of nodes in subgraph $\mathcal{SG}_i'$, and the superscript $(\cdot)^k$ denotes the itertion counter.
For representation we denote $\mathcal{I}_1$ and $\mathcal{I}_2$ as two separate sets of incident edges where $\mathcal{I} := \mathcal{I}_1 \cup \mathcal{I}_2$.
More specifically, $\mathcal{I}(\mathcal{SG}_i^\prime)$ is the set of edges incident to $\mathcal{SG}_i'$ (i.e. edges that contain a node in another subgraph) and $\mathcal{I}_1$ and $\mathcal{I}_2$
denote how the incident linking constraints are formulated within the subproblem using either primal or dual coupling.

In our formulation, \eqref{eq:Schwarz-objective} is the subproblem objective which is the sum of node objective functions contained in the expanded
subgraph $\mathcal{SG}_i^\prime$ where we have added the dual penalties from the incident dual linking constraints on edges $e \in \mathcal{I}_1(\mathcal{SG}_i^\prime)$, and
\eqref{eq:primal-received} represents primal constraints we have added from the edges $e \in \mathcal{I}_2(\mathcal{SG}_i^\prime)$. Note that
incident linking constraints can be directly incorporated into the subproblem as local constraints, or included in the objective function as a
dual penalty (this assigning procedure can be important to the algorithm performance). In our implementation, one may provide preference on whether a linking constraint is treated in primal or dual form (we show in Section \ref{sec:dcopf} how to do this).
The primal-dual solution of the other subproblems obtained from the previous iteration,
in particular $\{\{x_{n}^k\}_{n \in \mathcal{N}(e)\setminus \mathcal{N}(\mathcal{SG}'_i)}\}_{e\in\mathcal{I}(\mathcal{SG}_i')}$ and
$\{\lambda_{e}^k\}_{e \in \mathcal{I}_1(\mathcal{SG}_i^\prime)}$,
also enter into the subproblem formulation. The Schwarz algorithm achieves convergence with this exchange of information.

An important step in the Schwarz algorithm is the {\it restriction} of the subproblem solution.
One obtains the primal-dual solution ${\{x^*_n}\}_{n \in \mathcal{N}(\mathcal{SG}_i^\prime)}$ and
${\{\{\lambda^*_{e}}\}_{e \in \mathcal{E}(\mathcal{SG}_i^\prime) \cup \mathcal{I}_2(\mathcal{SG}_i^\prime)}\}$
by solving \eqref{eq:Schwarz_subproblem}.
We observe that the solutions from different subproblems overlap (the solution for overlapping nodes may appear more than once). We thus use a restriction step to eliminate such multiplicity; in particular, we discard the part of the solution associated with the nodes that are acquired by expansion. The restriction step can be expressed as:
\begin{align*}
  \forall i \in \{1,2,\cdots,N\}, \quad {\{x^k_n}\}_{n \in \mathcal{N}(\mathcal{SG}_i)},\;{\{\lambda^k_{e}}\}_{e\in
  \mathcal{E}(\mathcal{SG}_i)} \leftarrow \text{solution of \eqref{eq:Schwarz_subproblem}}.
\end{align*}
The solution of subproblems $i=1,2,\cdots,N$ can be performed in parallel. Next, the primal-dual residual at any Schwarz iteration $k$ is evaluated according to the residual to the KKT system for \eqref{eq:Schwarz_residual}.
We define $r^{k,Pr}_e$ as the primal residual of the linking constraints on edge $e$ at iteration $k$, and $r^{k,Du}_{e}$ as the dual residual of the linking constraints on edge $e$. Specifically, the primal error evaluates the linking constraints in the higher level graph $e \in \mathcal{E}(\mathcal{G})$, and the dual error evaluates the consensus of the dual values of the expanded subgraphs that contain the edge $e$; for the case where more than two subgraphs are associated with one linking constraint, see \cite[Proposition 1]{Shin2020b} for details on how to evaluate dual infeasiblility. The formal definition of these residuals is given by:
\begin{subequations}\label{eq:Schwarz_residual}
    \begin{align}
        & r^{k,Pr}_e := g_e(\{x^k_n\}_{n \in \mathcal{N}(e)}), \quad e \in \mathcal{E}(\mathcal{G}),\\
        & r^{k,Du}_e := \lambda^k_e(\mathcal{SG}^{\prime}_i) - {\lambda}^k_e(\mathcal{SG}^{\prime}_j), \quad  e \in \mathcal{E}(\mathcal{G}).\label{eq:Schwarz_residual-dual}
    \end{align}
\end{subequations}
The termination criteria can be set as follows:
\begin{align}\label{eq:stop}
  \text{stop if: }\max_{e\in\mathcal{E}(\mathcal{G})} \Vert r^{k,Pr}_e\Vert_\infty,\max_{e\in\mathcal{E}(\mathcal{G})} \Vert r^{k,Du}_e\Vert_\infty \leq \epsilon^{tol},
\end{align}
where $\epsilon^{tol}$ is the prescribed convergence tolerance.

\begin{figure}[]
\centering
\begin{subfigure}[t]{0.4\textwidth}
    \centering
    \includegraphics[width = 6cm,keepaspectratio]{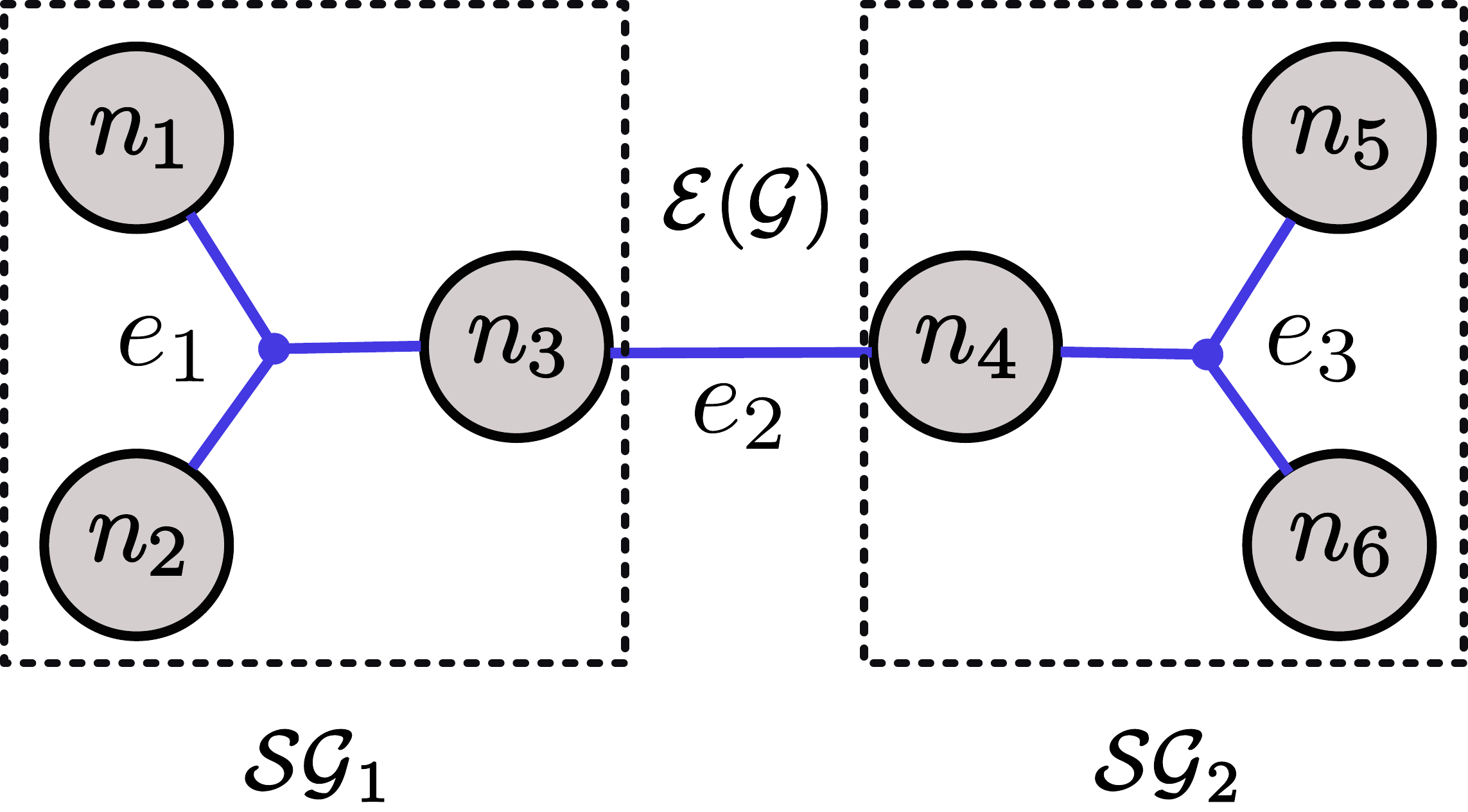}
    \label{fig:Schwarz1}
\end{subfigure}
\begin{subfigure}[t]{0.58\textwidth}
    \centering
    \includegraphics[width = 8cm,keepaspectratio]{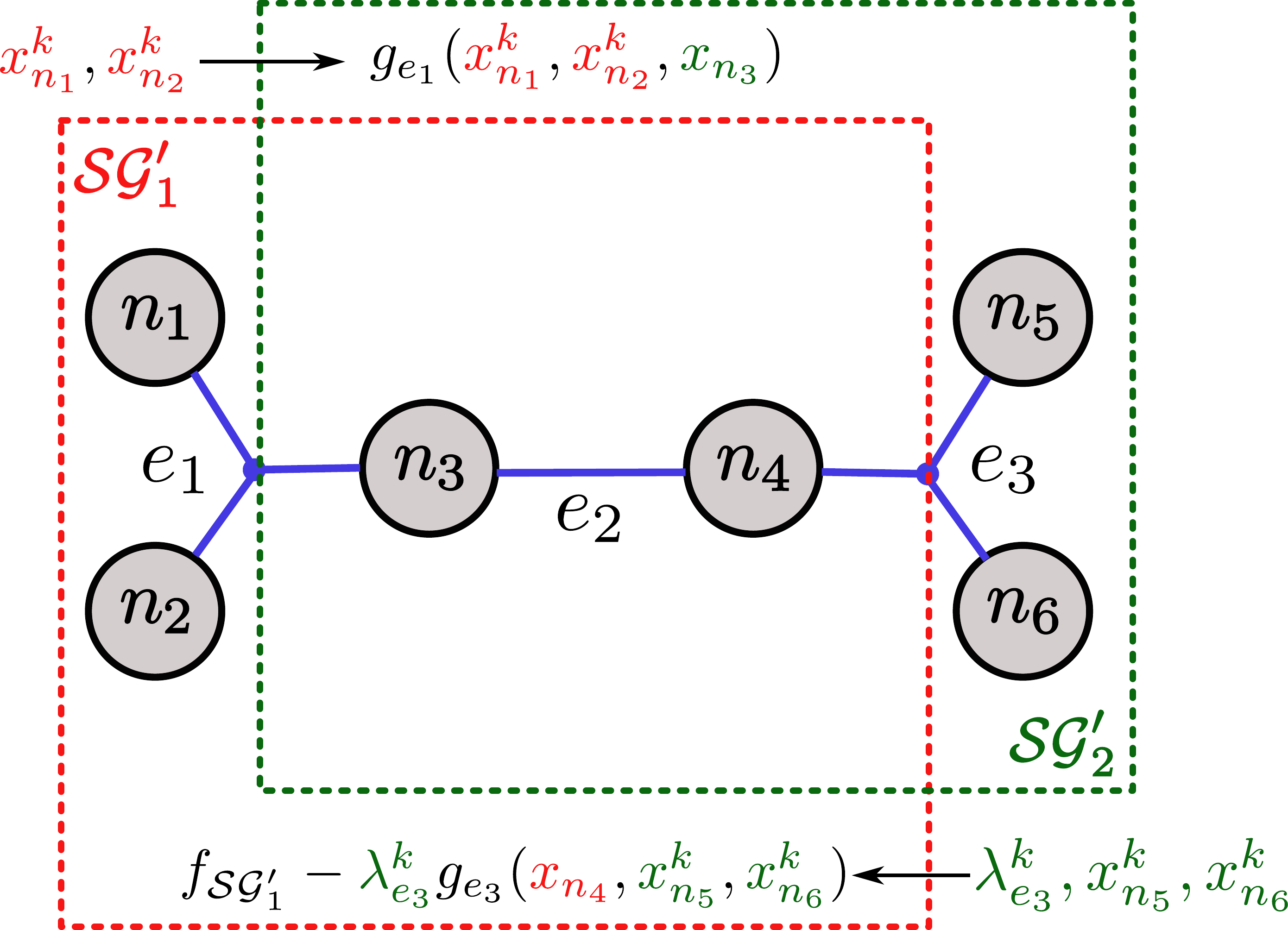}
    \label{fig:Schwarz2}
\end{subfigure}
\caption{Depiction of Schwarz Algorithm.  The original graph containing two subgraphs ($\mathcal{SG}_1$ and $\mathcal{SG}_2$) c
onnected by edge $e_2$ (left), graph with expanded subgraphs ($\mathcal{SG}^\prime_1$ and $\mathcal{SG}^\prime_2$) (right). The expanded subgraphs overlap at nodes $n_3$ and $n_4$.}
\label{fig:Schwarz_depiction}
\end{figure}

The Schwarz algorithm can be expressed using syntax that closely matches that of the \\ {\tt OptiGraph} abstraction, as shown in Algorithm \ref{alg:Schwarz}.
Figure \ref{fig:Schwarz_depiction} depicts how
primal-dual information is exchanged in the overlapping subgraph scheme.
The figure also depicts a simple graph that contains two subgraphs $\mathcal{SG}_1$ and $\mathcal{SG}_2$ and one edge $e_2$
that connects the  subgraphs ($e_2 \in \mathcal{E}(\mathcal{G})$).
The right side of Figure \ref{fig:Schwarz_depiction} shows the expanded subgraphs $\mathcal{SG}_1^\prime$ and
$\mathcal{SG}_2^\prime$.  In this illustration, edge $e_1$ is incident to subgraph $\mathcal{SG}_2^\prime$ and communicates primal
information (i.e. it is in the set $\mathcal{I}_2(\mathcal{SG}_2^\prime)$), and edge $e_3$ is incident to subgraph $\mathcal{SG}_1^\prime$
and communicates dual information to subgraph $\mathcal{SG}_1^\prime$ (i.e. it is in the set $\mathcal{I}_1(\mathcal{SG}_1^\prime)$).

\begin{algorithm}[H]
\caption{Schwarz Algorithm for Solving an {\tt OptiGraph}}\label{alg:Schwarz}
\begin{algorithmic}
\STATE Input graph $\mathcal{G}$, non-expanded subgraphs $\{ \mathcal{SG}_1,...,\mathcal{SG}_N \}$ and expanded subgraphs $\{ \mathcal{SG}^\prime_1,...,\mathcal{SG}^\prime_N \}$.
\STATE Initialize $x^0$, $\lambda^0$
\STATE Formulate subproblems in \eqref{eq:Schwarz_subproblem}
\WHILE{termination criteria not satisfied}
    \FOR{$i = 1:N$ (in parallel)}
        \STATE Retrieve $\{\{x_{n}^k\}_{n \in \mathcal{N}(e)\setminus \mathcal{N}(\mathcal{SG}'_i)}\}_{e\in\mathcal{I}(\mathcal{SG}_i')}$ and
        $\{\lambda_{e}^k\}_{e \in \mathcal{I}_1(\mathcal{SG}_i^\prime)}$ and update subproblems.
        \STATE Solve subproblem \eqref{eq:Schwarz_subproblem} to obtain ${\{x^{k+1}_n}\}_{n \in \mathcal{N}(\mathcal{SG}_i)},\;\{\lambda^{k+1}_{e}\}_{e\in\mathcal{E}(\mathcal{SG}_i)}$
    \ENDFOR
    \STATE Compute $\{r^{k,Pr}_e\}_{e\in \mathcal{E}(\mathcal{G})}$, $\{r^{k,Du}_e\}_{e\in\mathcal{E}(\mathcal{G})}$ and evaluate termination criteria \eqref{eq:stop}.
\ENDWHILE
\end{algorithmic}
\end{algorithm}

\subsubsection{Implementation of Schwarz Algorithm}

We implemented the Schwarz algorithm in {\tt Plasmo.jl} and call this {\tt SchwarzSolver.jl}.
Code Snippet \ref{code:Schwarz_snippet} shows how we can solve the overlapping subgraphs we produced in
Code Snippet \ref{code:expand_snippet} corresponding to the overlap shown in Figure \ref{fig:ex5_plots}. On Line \ref{line:Schwarz_solver} we import the Schwarz solver, on Line \ref{line:Schwarz_ipopt} we create an {\tt Ipopt} optimizer to use as the subproblem solver, and we execute the algorithm on Line \ref{line:Schwarz_solve}. It is also possible to pass an overlap distance instead of the expanded subgraphs and allow the Schwarz solver to perform the subgraph expansion.  The benefit of formulating subgraphs at the user level is that one could formulate custom overlap schemes.

\begin{minipage}[]{0.9\linewidth}
\begin{scriptsize}
\lstset{language=Julia,breaklines = true}
\begin{lstlisting}[caption = Using the Schwarz.jl Solver, label = {code:Schwarz_snippet}]
using SchwarzSolver    |\label{line:Schwarz_solver}|
using Ipopt

#Use Ipopt as the subproblem solver
ipopt = Ipopt.Optimizer           |\label{line:Schwarz_ipopt}|

#Solve with Schwarz Algorithm
schwarz_solve(graph,expanded_subgraphs,                            |\label{line:Schwarz_solve}|
sub_optimizer = ipopt,max_iterations = 100,tolerance = 1e-6)
\end{lstlisting}\label{fig:Schwarz_problem_snippet}
\end{scriptsize}
\end{minipage}

\section{Case Studies}\label{sec:case_study}

We demonstrate modeling and solution capabilities enabled by the {\tt OptiGraph} abstraction and the {\tt Plasmo.jl} implementation using a couple of challenging problems
arising in infrastructure networks. The implementations used for the case studies can be found at \url{https://github.com/zavalab/JuliaBox/tree/master/PlasmoExamples}.
We highlight that the study presented in Section \ref{sec:case_study_pips} uses a dataset that cannot be shared (we only present computational results to demonstrate scalability).  In the GitHub repository we provide a simple network dataset that uses the same features of {\tt Plasmo.jl}.

\subsection{Linear Algebra Decomposition of a Natural Gas Network}\label{sec:case_study_pips}
We solve a large dynamic optimization problem for the natural gas network shown in Figure \ref{fig:gas_network_layout}.
This network includes four gas supplies, 153 time-varying gas demands, 215 pipelines, and 16 compressor stations.
The goal of the optimization problem is to track demands and minimize compressor power over a 24-hour time horizon.  The main complexity of this problem arises from
nonlinear PDEs needed to capture conservation of mass and momentum inside the pipelines. After space-time discretization of the PDE, we obtain a nonconvex NLP with
432,048 variables, 427,512 equality constraints, and 3,887 inequality constraints. Figure \ref{fig:gas_ocp_space_time} depicts
the graph structure induced by the space-time nature of this problem. Our goal is to identify efficient partitions that traverse space-time to efficiently solve the problem
using Schur decomposition in {\tt PIPS-NLP}.

All model details can be found in \cite{Zavala2014} and \cite{Jalving2018} and in the scripts provided with this work.  Here we provide a summary to
showcase modeling features of {\tt Plasmo.jl}.  The gas network contains links and junctions $Net(\mathcal{J},\mathcal{L})$.  The set of junctions
$j \in \mathcal{J}$ connect links and include gas supplies (injections) $\mathcal{S}_j$ and  demands (withdrawals) $\mathcal{D}_j$.
The set of links $\mathcal{L}$ is composed of both pipeline links $\mathcal{L}_p$ and compressor links $\mathcal{L}_c$ such that
$\mathcal{L} := \mathcal{L}_p \cup \mathcal{L}_c$. We also specify the set of time periods as $\mathcal{T} := \{1,...,N_t\}$ with a 24-hour horizon
and $\overline{\mathcal{T}} := \mathcal{T} \setminus N_t$.
For our implementation, each component of the system is modeled as a stand-alone
{\tt OptiGraph} with {\tt OptiNodes} distributed over time and we connect these in a high-level graph to form the complete problem.
This modular construction is one of the benefits of the graph abstraction. Specifically, component models can be developed separately
(e.g., by different people). Moreover, the syntax of the different components can be different because each {\tt OptiGraph} is a self-contained object.
This greatly facilitates model construction and debugging. This differs from standard algebraic modeling languages, in which the syntax in the entire model has to be
compatible (this complicates debugging of large models). The implementations to construct each component model (each {\tt OptiGraph}) in {\tt Plasmo.jl}
can be found in the Appendix.

\paragraph{Junction OptiGraph}\mbox{}\\
The gas junction model is described by \eqref{eq:gas_junctions}, where $\theta_{j,t}$ is the pressure at junction
$j$ and time interval $t$. $\underline{\theta}_j$ is the lower pressure bound for the
junction, $\overline{\theta}_j$ is the upper pressure bound, $f_{j,d,t}^{target}$ is the target demand flow for demand $d$ on junction $j$
and $\overline{f}_{j,s}$ is the available gas generation from supply $s$ on junction $j$.
Code Snippet \ref{code:junction_construction_snippet} in the Appendix shows how to create the junction graph.

\begin{subequations}\label{eq:gas_junctions}
    \begin{align}
        & \underline{\theta}_n \le \theta_{j,t} \le \overline{\theta}_n , \quad j \in \mathcal{J}, \ t \in \mathcal{T}    \label{eq:node_pressure} \\
        &0 \le f_{j,d,t} \le f_{j,d,t}^{target}, \quad d \in \mathcal{D}_j, \ j \in \mathcal{J}, \ t \in \mathcal{T} \label{eq:demand_target} \\
        &0 \le f_{j,s,t} \le \overline{f}_{j,s}, \quad s \in \mathcal{S}_j, \ j \in \mathcal{J}, \ t \in \mathcal{T}. \label{eq:supply_limit}
    \end{align}
\end{subequations}

\paragraph{Pipeline OptiGraph}\mbox{}\\
Each pipeline model is an {\tt OptiGraph} with {\tt OptiNodes} distributed on a space-time grid. Specifically, the nodes of each pipeline graph form a $N_t \times N_x$ mesh
wherein pressure and flow variables are assigned to each node. Flow dynamics within pipelines are expressed with linking constraints that describe the discretized
PDE equations for mass and momentum using finite differences. We assume isothermal flow through each horizontal pipeline segment $\ell \in \mathcal{L}_p$ which are
described by the conservation equations from \cite{Osiadacz1984}.  The pipeline formulation is given by \eqref{eq:pipeline_discretized}
where $p_{\ell,t,k}$ and $f_{\ell,t,k}$ are the pipeline pressure and flow for each pipeline link $\ell$ at each time period $t \in \mathcal{T}$ for each
space point $k \in \mathcal{X}$ where $\mathcal{X} := \{1,...,N_x\}$ is the set of spatial discretization points we use for each pipeline.
We also define $\overline{\mathcal{X}} :=  \mathcal{X} \setminus N_x$ for convenience.
The constants given in this formulation for each pipeline $\ell \in \mathcal{L}_p$ are $c_{1,\ell} = \alpha \frac{c^2}{A_\ell}$, $c_{2,\ell} = \frac{1}{\alpha}A_\ell$,
and $c_{3,\ell} = \alpha \frac{8c^2\lambda_\ell A_{\ell}}{\pi^2D_{\ell}^5}$. Here, the pipeline terms are the cross-sectional
area $A_\ell$, diameter $D_\ell$, friction coefficient $\lambda_\ell$,
the gas speed of sound squared $c^2$, and the scaling coefficient $\alpha$.  We denote expressions for
the flow in and out of each pipeline at each time as $f_{\ell,t}^{in}$ and $f_{\ell,t}^{out}$ as well as the total line pack (inventory of gas)
in each pipeline link at each time as $m_{\ell,t}$. With the defined terms,
\eqref{eq:mass} and \eqref{eq:momentum} are the mass and momentum equations, \eqref{eq:dummyflowin}, \eqref{eq:dummyflowout}, \eqref{eq:dummypressurein}, \eqref{eq:dummypressureout}
define the convenience variables for flow and pressure in and out of each pipeline, \eqref{eq:initial_ss_mass} and \eqref{eq:initial_ss_momentum}
prescribe the system to initially be at steady-state, and \eqref{eq:def_linepack} and \eqref{eq:refill_linepack} define line pack and require
each pipeline to refill its line pack by the end of the time horizon.

Capturing the space-time structure of \eqref{eq:pipeline_discretized} is seemingly complex but
it is straightforward to do so with an {\tt OptiGraph} because each pipeline can be treated independently. Code Snippet \ref{code:pipeline_construction_snippet}
in the Appendix details how each pipeline model is constructed.

\begin{subequations}\label{eq:pipeline_discretized}
    \begin{align}
        & \frac{p_{\ell,t+1,k} - p_{\ell,t,k}}{\Delta t} =  -c_{1,\ell}\frac{f_{\ell,t+1,k+1} - f_{\ell,t+1,k}}{\Delta x_\ell},
        \ell \in \mathcal{L}_p,t \in \overline{\mathcal{T}}, k \in \overline{\mathcal{X}}, \label{eq:mass}\\
        & \frac{f_{\ell,t+1,k} - f_{\ell,t,k}}{\Delta t} = -c_{2,\ell}\frac{p_{\ell,t+1,k+1} - p_{\ell,t+1,k}}{\Delta x_\ell} -
        c_{3,\ell}\frac{f_{\ell,t+1,k}\mods{f_{\ell,t+1,k}}}{p_{\ell,t+1,k}},
        \ell \in \mathcal{L}_p, t \in \overline{\mathcal{T}}, k \in \overline{\mathcal{X}} \label{eq:momentum}\\
        & f_{\ell,t,N_x} = f_{\ell,t}^{out},\quad \ell \in \mathcal{L}_p, t \in \mathcal{T} \label{eq:dummyflowin}\\
        & f_{\ell,t,1} = f_{\ell,t}^{in},\quad \ell \in \mathcal{L}_p, t \in \mathcal{T}    \label{eq:dummyflowout}\\
        & p_{\ell,t,N_x} = p_{\ell,t}^{out},\quad \ell \in \mathcal{L}_p, t \in \mathcal{T} \label{eq:dummypressurein}\\
        & p_{\ell,t,1} = p_{\ell,t}^{in},\quad \ell \in \mathcal{L}_p, t \in \mathcal{T}    \label{eq:dummypressureout}\\
        & \frac{f_{\ell,1,k+1} - f_{\ell,1,k}}{\Delta x_{\ell}} = 0,  \quad \ell \in \mathcal{L}_{p}, k \in \overline{\mathcal{X}} \label{eq:initial_ss_mass} \\
        & c_{2,\ell}\frac{p_{\ell,1,k+1} - p_{\ell,1,k}}{\Delta x_{\ell}} + c_3\frac{f_{\ell,1,k}\mods{f_{\ell,1,k}}}{p_{\ell,1,k}} = 0,
        \quad \ell \in \mathcal{L}_p, k \in \overline{\mathcal{X}} \label{eq:initial_ss_momentum}\\
        & m_{\ell,t} = \frac{A_\ell}{c^2} \sum_{k=1}^{N_x} p_{\ell,t,k} \Delta x_{\ell},\quad \ell \in \mathcal{L}_p, t \in \mathcal{T} \label{eq:def_linepack}\\
        & m_{\ell,N_t} \ge m_{\ell,1}, \quad  \ell \in \mathcal{L}_p \label{eq:refill_linepack}.
    \end{align}
\end{subequations}

\paragraph{Compressor OptiGraph}\mbox{}\\
The compressor models are created analogously to the junction and pipeline models.  We use an ideal isentropic formulation given by
\eqref{eq:ideal_compressor} where $\eta_{\ell,t}$, $p_{\ell,t}^{in}$, $p_{\ell,t}^{out}$, and $P_{\ell,t}$ are the compression ratio,
suction presure, discharge pressure, and horsepower respectively for compressor link $\ell$ at time $t$.  We also define
constant parameters, $c_p$, $T$, and $\gamma$ for heat capacity, temperature and isentropic efficiency. The variables
$f_{\ell,t,in}$ and $f_{\ell,t,out}$ are used for convenience to be consistent with the pipeline link formulation.
The compressor graph construction is straightforward and is shown in Code Snippet \ref{code:compressor_modelgraph} in the Appendix.

\begin{subequations}\label{eq:ideal_compressor}
    \begin{align}
        &p_{\ell,t}^{out} = \eta_{\ell,t} p_{\ell,t}^{in}, \quad \ell \in \mathcal{L}_c, t \in \mathcal{T} \label{eq:compressor_boost}\\
        &P_{\ell,t} = c_p \cdot T \cdot f_{\ell,t} \left(\left(\frac{p_{\ell}^{out}}{p_{\ell}^{in}}\right)^{\frac{\gamma-1}{\gamma}}-1\right),
        \quad \ell \in \mathcal{L}_c, t \in \mathcal{T} \label{eq:compressor_power}\\
        &f_{\ell,t} = f_{\ell,t}^{in} = f_{\ell,t}^{out}, \quad \ell \in \mathcal{L}_c, t \in \mathcal{T}. \label{eq:comprssor_flow}
    \end{align}
\end{subequations}

\paragraph{Network OptiGraph}\mbox{}\\
The network structure of the model is induced using a high-level graph wherein we use linking constraints
to express mass conservation around each junction and boundary conditions for pipeline and compressor links according
to the formulation given by \eqref{eq:network_links}. Here, \eqref{eq:junction_conservation} represents mass conservation
at each junction $j$ and \eqref{eq:bc1} and \eqref{eq:bc2} are the pipeline and compressor link boundary conditions. We also
define $\theta_{rec(\ell),t}$ and $\theta_{snd(\ell),t}$ as the receiving and sending junction for each link $\ell \in \mathcal{L}$ at time $t$
and we define $\mathcal{L}_{rec}(j)$ and $\mathcal{L}_{snd}(j)$ as the set of receiving and sending links to each junction $j$ respectively.
The final optimal control problem is given by \eqref{eq:network_formulation} which seeks to minimize the total operating cost over
the 24-hour period where $\alpha_{\ell}$ is the compressor cost for compressor link $\ell$ and $\alpha_{j,d}$ is the gas delivery price for demand $d$ on junction $j$.
The complete formulation includes all of the junction, pipeline, compressor, and network equations.
Code Snippet \ref{code:gas_network_construction} in the Appendix assembles the complete natural gas model.

\begin{subequations}\label{eq:network_links}
    \begin{align}
        & \sum_{\ell\in\mathcal{L}_{rec}(j)} f^{out}_{\ell,t} - \sum_{\ell \in\mathcal{L}_{snd}(j)} f^{in}_{\ell,t} +
        \sum_{s\in\mathcal{S}_j}f_{j,s,t} - \sum_{d\in \mathcal{D}_j}f_{j,d,t} = 0, \quad j \in\mathcal{J}, t \in \mathcal{T} \label{eq:junction_conservation}\\
        & p_{\ell,t}^{in} = \theta_{rec(\ell),t}, \quad \ell \in \mathcal{L},  t \in \mathcal{T} \label{eq:bc1}\\
        & p_{\ell,t}^{out} = \theta_{snd(\ell),t}, \quad \ell \in \mathcal{L}, t \in \mathcal{T} \label{eq:bc2}
    \end{align}
\end{subequations}

\begin{subequations}\label{eq:network_formulation}
    \begin{align}
        \min_{ \substack{ \{ \eta_{\ell,t},f_{j,d,t} \} \\ \ell \in \mathcal{L}_c, d \in \mathcal{D}_j, j \in \mathcal{J}, t \in \mathcal{T}}} \quad &
        \sum_{\substack{\ell \in \mathcal{L}_c \\ t \in \mathcal{T}}} \alpha_{\ell} P_{\ell,t} -
        \sum_{\substack{d \in \mathcal{D}_j, j \in \mathcal{J}, \\  t \in \mathcal{T}}} \alpha_{j,d} f_{j,d,t} &\\
         s.t. \quad & \text{Junction Limits} & \eqref{eq:gas_junctions} \\
         & \text{Pipeline Dynamics}  & \eqref{eq:pipeline_discretized} \\
         & \text{Compressor Equations} & \eqref{eq:ideal_compressor} \\
         & \text{Network Link Equations} & \eqref{eq:network_links}
    \end{align}
\end{subequations}

\paragraph{Partitioning and Solving the Natural Gas Problem}\mbox{}\\

\begin{figure}[]
\centering
\begin{subfigure}[t]{0.48\textwidth}
    \centering
    \includegraphics[width = 1.0\textwidth,keepaspectratio]{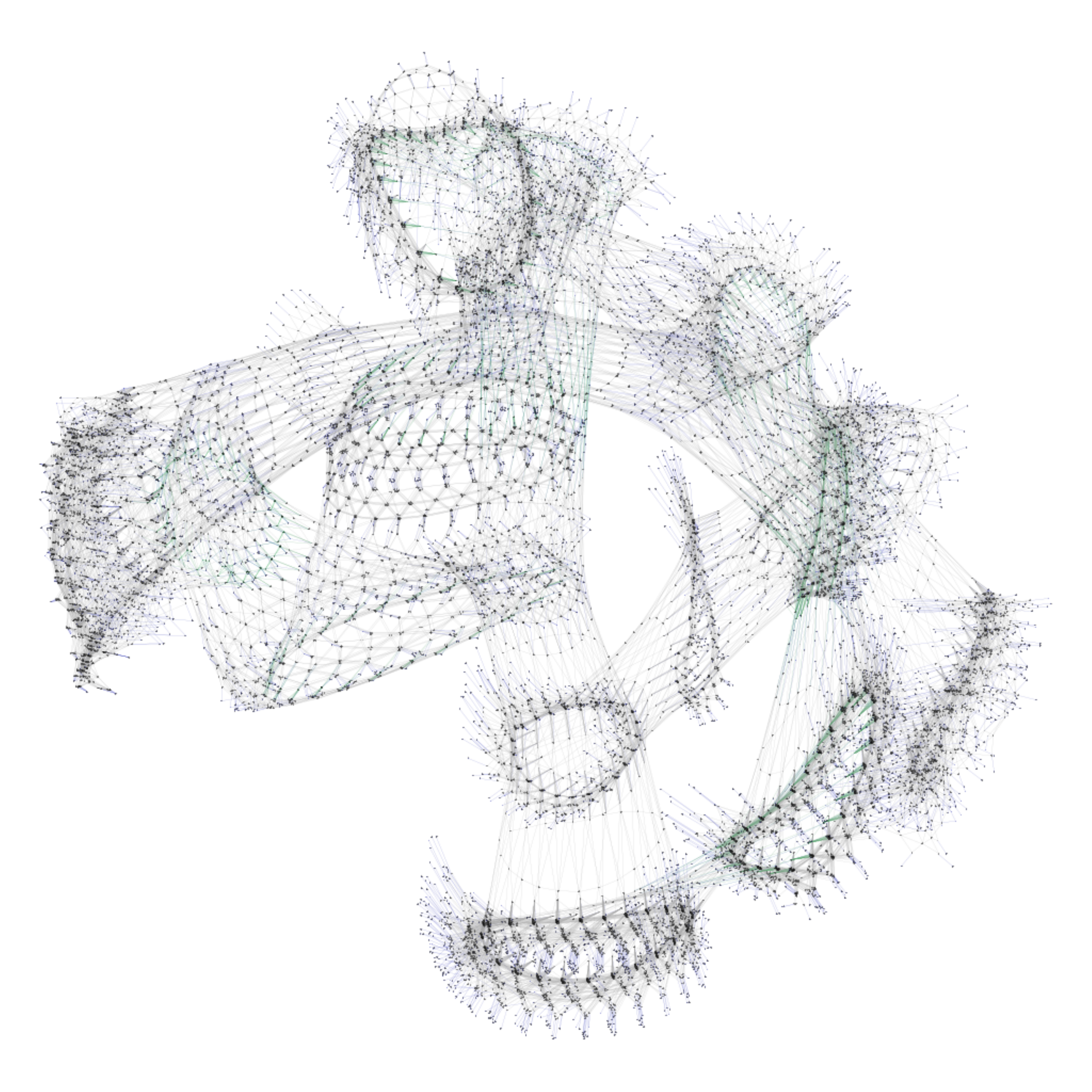}
    \caption{Modeled Components}
    \label{fig:network_components}
\end{subfigure}
\begin{subfigure}[t]{0.45\textwidth}
    \centering
    \includegraphics[width = 1.0\textwidth,keepaspectratio]{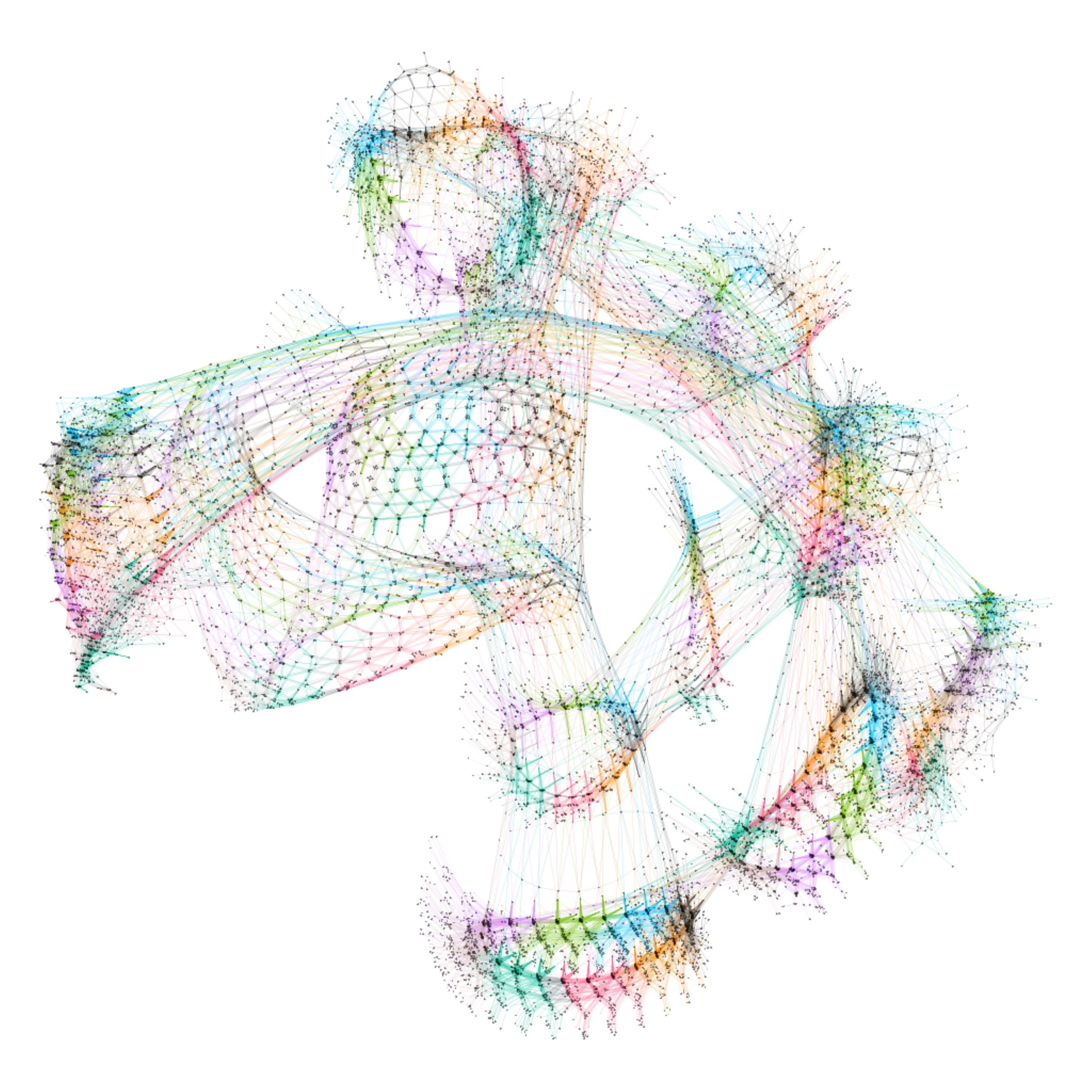}
    \caption{Time Partitions}
    \label{fig:network_time_partition}
\end{subfigure}
\begin{subfigure}[t]{0.48\textwidth}
    \centering
    \includegraphics[width = 1.0\textwidth,keepaspectratio]{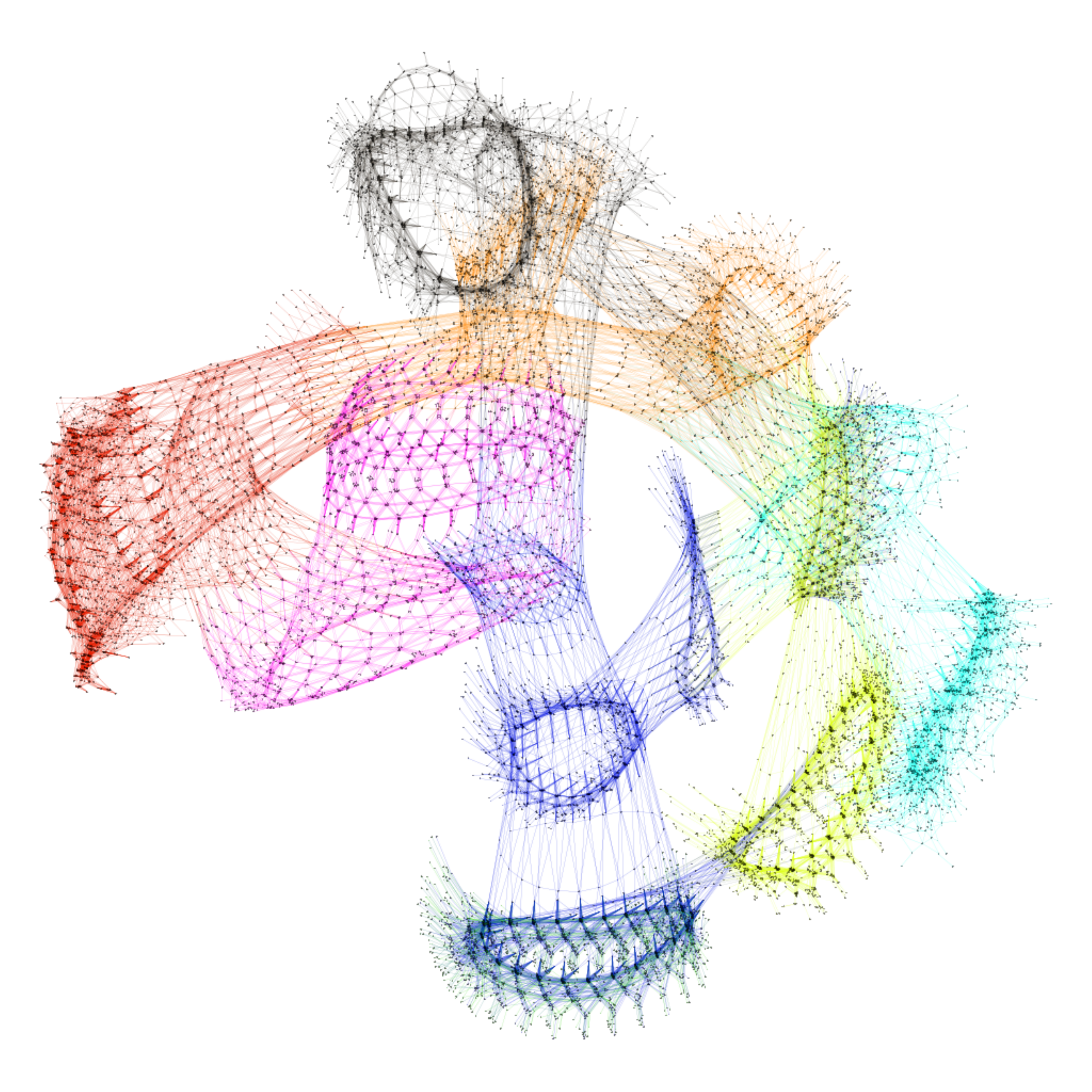}
    \caption{Network Partitions}
    \label{fig:network_network_partition}
\end{subfigure}
\begin{subfigure}[t]{0.45\textwidth}
    \centering
    \includegraphics[width = 1.0\textwidth,keepaspectratio]{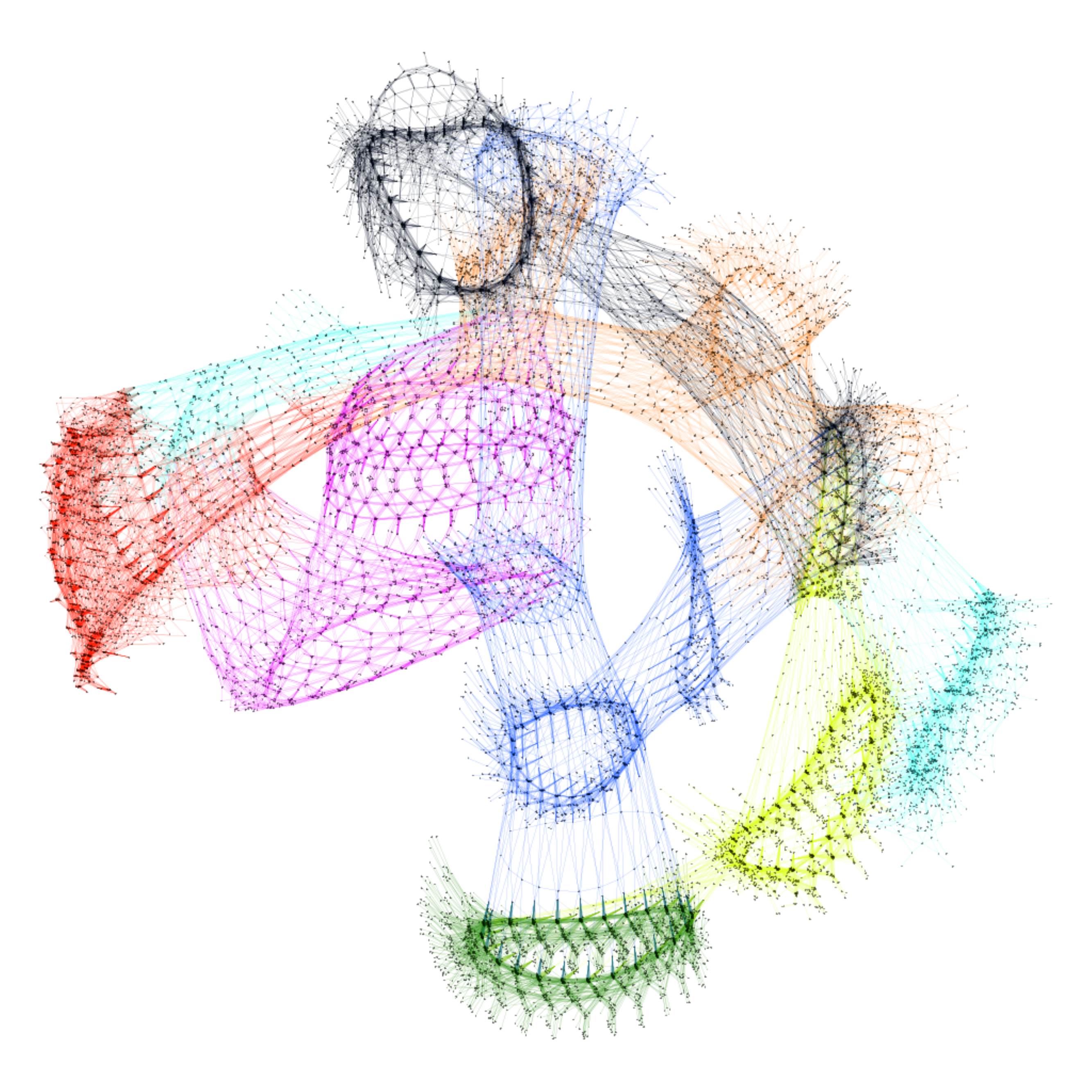}
    \caption{Space-Time Partitions}
    \label{fig:network_space_time_partition}
\end{subfigure}
\caption{Graph layouts of the natural gas network problem.  The graph is colored
by the physical network components (top left), by 8 time partitions (top right), by 8 network partitions (bottom left),
and by 8 space-time partitions (bottom right).}
\label{fig:gas_network_partitions}
\end{figure}

Once we have constructed a graph representation of the problem, we can use
the {\tt Plasmo.jl} tools to explore different partitioning strategies. Figure \ref{fig:gas_network_partitions} visualizes the graph structure of the problem and
shows different partitions.  Figure \ref{fig:network_components} shows the graph components
we assembled in the above snippets with pipeline nodes depicted in grey, compressor nodes in green, and junctions in blue.
Figure \ref{fig:network_time_partition} depicts a pure time partition of the problem with 8 partitions (each with a different color).
Partitioning in time is an intuitive approach but this results in 87,505 linking constraints (making the Schur matrix intractable).
This problem can also be partitioned purely as a network wherein we consider the
partition of the network components themselves (as opposed to the structure of the problem) which is
shown in Figure \ref{fig:network_network_partition}.  The network partition is physically intuitive but it does not
capture the true mathematical structure of the problem nor consider load balancing. We highlight that both time and network partitioning approaches can be performed using
the partitioning framework (by manually defining a partition vector), but the value of having the graph is that
we can efficiently obtain space-time partitions to produce the partition shown in Figure \ref{fig:network_space_time_partition}.
While visually similar to the network partition, the space-time
partition produces considerably less coupling (748 linking constraints against 1,800 linking constraints),
thus leading to a small Schur matrix. We also highlight that the space-time partition requires traversing space-time.

Code Snippet \ref{code:network_partition_snippet} shows how to partition the gas network problem and produce the problem we
communicate to {\tt PIPS-NLP}.  Line \ref{line:gas_network_kahypar} imports the {\tt KaHyPar} hypergraph partitioner,
Line \ref{line:gas_network_hypergraph} formulates the hypergraph representation of our problem, and Lines \ref{line:gas_network_weights_start}
through \ref{line:gas_network_weights_end} setup weight vectors we use for the nodes and edges (we weight nodes by their number of variables and edges
by their number of linking constraints). We partition the {\tt hypergraph} on Line \ref{line:gas_network_part}, create a {\tt Partition} object on
Line \ref{line:gas_network_part_object}, produce new subgraphs on Line \ref{line:gas_network_make_subgraphs} and
finally aggregate the subgraphs on Line \ref{line:gas_network_combine}.

\begin{minipage}[]{0.9\linewidth}
\begin{scriptsize}
\lstset{language=Julia,breaklines = true}
\begin{lstlisting}[caption = Partitioning the gas network control problem,label = {code:network_partition_snippet}]
#Import the KaHyPar interface
using KaHyPar         |\label{line:gas_network_kahypar}|

#Get they hypergraph representation of the gas network
hypergraph,ref_map = gethypergraph(gas_network)      |\label{line:gas_network_hypergraph}|

#Setup node and edge weights
n_vertices = length(vertices(hypergraph))                |\label{line:gas_network_weights_start}|
node_weights = [num_variables(node) for node in all_nodes(gas_network)]
edge_weights = [num_link_constraints(edge) for edge in all_edges(gas_network)]   |\label{line:gas_network_weights_end}|

#Use KaHyPar to partition the hypergraph
node_vector = KaHyPar.partition(hypergraph,8,configuration = :edge_cut,
imbalance = 0.01, node_weights = node_weights,edge_weights = edge_weights)      |\label{line:gas_network_part}|

#Create a model partition
partition = Partition(gas_network,node_vector,ref_map)           |\label{line:gas_network_part_object}|

#Setup subgraphs based on partition
make_subgraphs!(gas_network,partition)               |\label{line:gas_network_make_subgraphs}|

#Aggregate the subgraphs into OptiNodes
new_graph , aggregate_map  = aggregate(gas_network,0)        |\label{line:gas_network_combine}|
\end{lstlisting}
\end{scriptsize}
\end{minipage}

\begin{figure}[]
\centering
\begin{subfigure}[t]{0.48\textwidth}
    \centering
    \includegraphics[width=7cm,keepaspectratio]{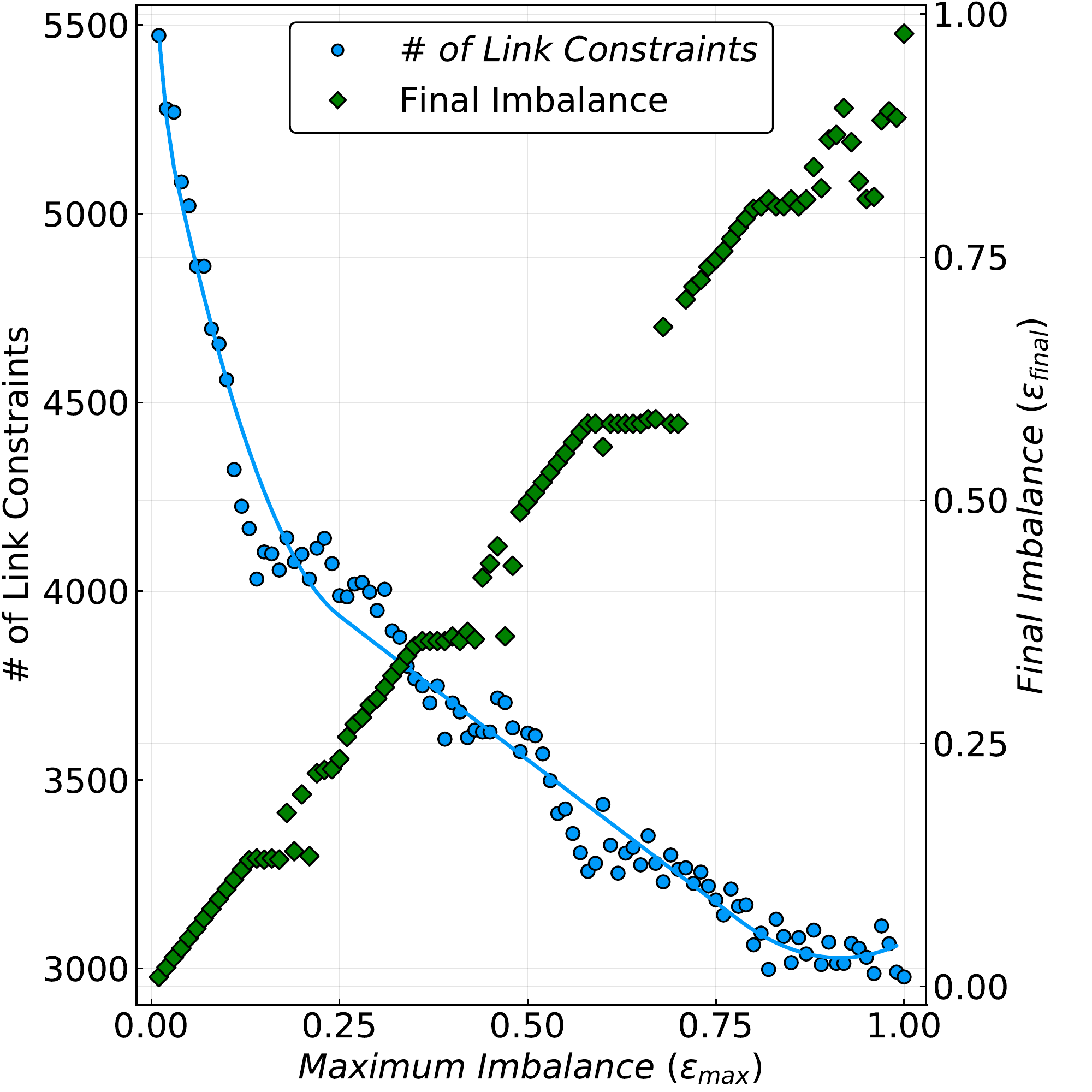}
    \caption*{}
    \label{fig:imbalance_full}
\end{subfigure}
\begin{subfigure}[t]{0.48\textwidth}
    \centering
    \includegraphics[width=7cm,keepaspectratio]{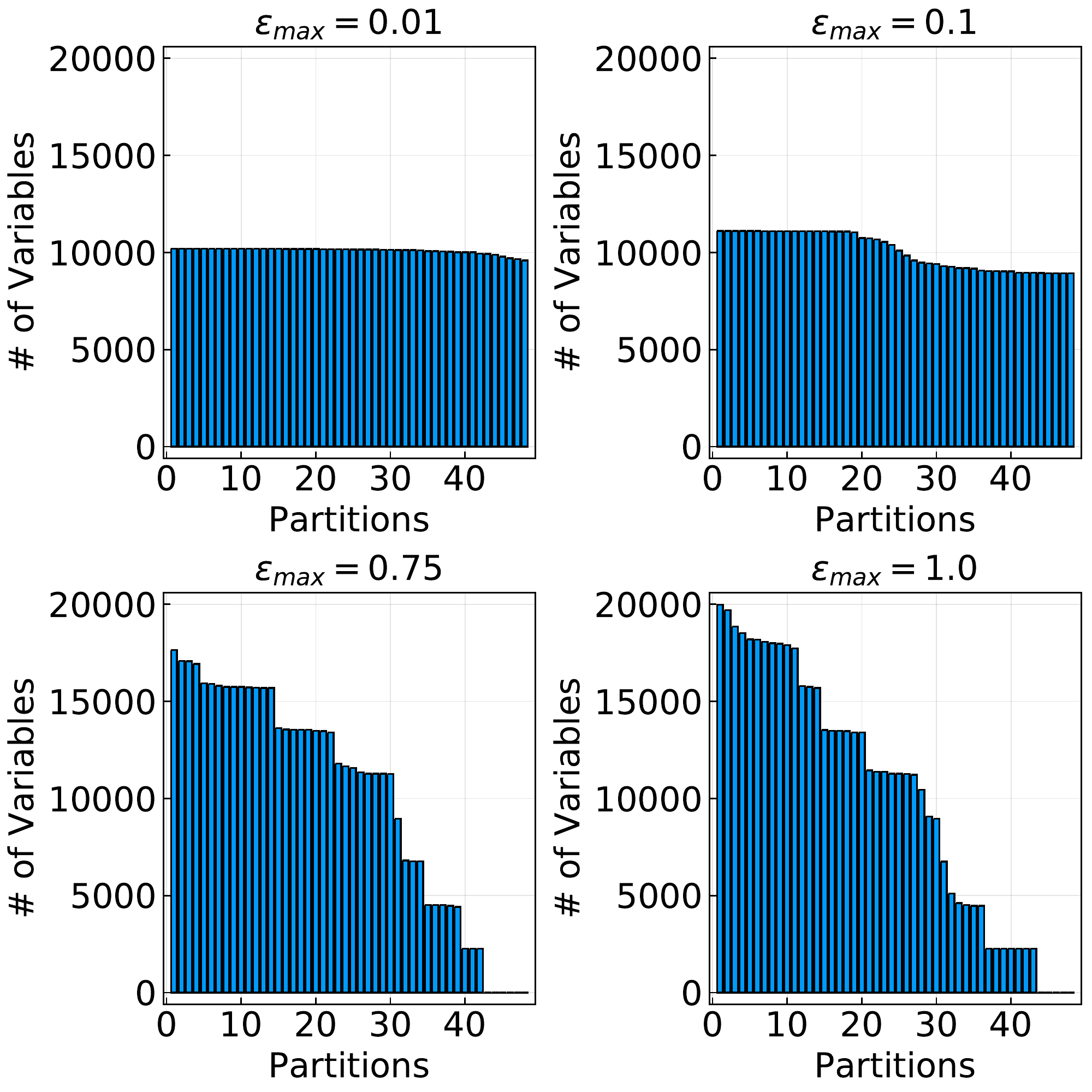}
    \caption*{}
    \label{fig:imbalance}
\end{subfigure}
\caption{Partitioning results with 48 partitions.  Imbalance versus the number of linking constraints and final imbalance (left) and
the distribution of subproblem sizes for select imbalance parameters (right).}
\label{fig:gas_network_imbalance}
\end{figure}

The partitioning in Snippet \ref{code:network_partition_snippet} produces a graph with 8 nodes which we can distribute and solve with {\tt PIPS-NLP}.
We follow the exact same setup used in Snippet \ref{code:pipsnlp_problem_snippet} using {\tt PipsSolver.jl} and distribute the graph between workers to solve in parallel
on a machine with 32 Intel  Xeon  CPUs  E5-2698  v3  running at 2.30  GHz.
We experiment with different numbers of partitions and imbalance values and explore how the {\tt PIPS-NLP} algorithm performs.
Table \ref{table:pips_results} details the results obtained where we vary the number of partitions $|\mathcal{P}|$ (8, 24, and 48) and the maximum
imbalance value $\epsilon_{max}$  (between 0.01 and 1.0). Figure \ref{fig:gas_network_imbalance} shows the effect of increasing the imbalance for the
48 partition case on the true final imbalance the partitioner produced ($\epsilon_{final}$) and the total number of linking constraints.
For the most part, the maximum and final imbalances $\epsilon_{final}$ display a one-to-one relationship but there are
distinct intervals wherein the final imbalance flattens out.
We also see that nominal values of maximum imbalance (less than 25\%) reduce the number of linking constraints considerably but greater values produce diminishing returns.
Here, we also show the distribution of subproblem sizes for a few select maximum imbalance values on the right side of Figure \ref{fig:gas_network_imbalance} for reference.
We can see that, under naive partitioning, the size of the partitions can vary by an order of magnitude.

For each run, we note the true imbalance {\tt KaHyPar} produced $\epsilon_{final}$, the sum of the edge weights $sum(w(\mathcal{E}))$ (which
corresponds to the number of linking constraints), the maximum node size\\
$max(\{s(n)\}_{n \in \mathcal{N}(\mathcal{G})})$ (which corresponds to the node
with the most variables), as well as the average and maximum
number of node connections (i.e. the number of linking constraints incident to a node).  For timing results we observe the time spent building the KKT system $t_{build}$ (i.e.
the time to build \eqref{eq:schur_first_stage}),
the time spent factorizing the
Schur complement $t_{fac}$ (i.e. time to solve \eqref{eq:schur_first_stage}) and the total time spent inside {\tt PIPS-NLP} $t_{pips}$.
The first rows of Table \ref{table:pips_results} show results with 8 partitions and also include cases for time and network partitions corresponding
to Figures \ref{fig:network_time_partition} and \ref{fig:network_network_partition}. As expected, the high degree of coupling in these partitions is
computationally prohibitive; the time partition cannot be solved with Schur complement decomposition and the network partition requires over 2 hours.
Significant improvement is achieved by using the {\tt KaHyPar} partitioner; we note that even a 1\% imbalance is twice as fast as the network partition.
Increasing the imbalance parameter to 50\% results in much better performance (the average and maximum number of columns in $B_n$ decreases), but increasing it to 60\% results in
diminished speed, despite producing a similar partition.  Interestingly, allowing too much imbalance can result in
ill-conditioned subproblems (e.g., ill-conditioned matrices $K_n$). This highlights why having flexibility in partitioning is important. The second and third sets of rows present the
results using 24 and 48 partitions with 24 CPUs.  Increasing the maximum imbalance
for the 24 partition case results in diminished performance (due to the bottleneck of building the KKT system).  This is because either the node connectivity
increases (for the 10\% imbalance),
the maximum subproblem size increases (for the 100\% imbalance), or there is some ill-conditioning of the subproblems.  In contrast, increasing the imbalance for the 48 partition case results in
computational improvement.  This is because more linking constraints (more coupling) shifts the bottleneck step to
factorizing the Schur complement and drastic speedups result from reducing the degree of coupling. This highlights how flexible partitions can help fit the solver to the computing architecture of interest.

\begin{table}
\scriptsize
\begin{center}
\begin{tabular}{|c|c|c|c|c|c|c|c|c|c|}
\hline
Partition & $\epsilon_{max}$ & $\epsilon_{final}$ & {\begin{tabular}{@{}c@{}} $\sum\limits_{e \in \mathcal{E}(\mathcal{G})}(w(e))$ \\ (\# \ Link Constraints)\end{tabular}} &
{\begin{tabular}{@{}c@{}}$max(\{s(n)\}_{n \in \mathcal{N}(\mathcal{G})})$  \\ (Largest node) \end{tabular}} &
{\begin{tabular}{@{}c@{}} {$\{\sum\limits_{e \in \mathcal{E}(n)}(w(e))\}_{n \in \mathcal{N}(\mathcal{G})}$} \\ (\# Incident Node Links) \\  mean , max \end{tabular}} &
{\begin{tabular}{@{}c@{}} $t_{build}$  \\ (sec) \end{tabular}} &
{\begin{tabular}{@{}c@{}} $t_{fac}$   \\ (sec) \end{tabular}} &
{\begin{tabular}{@{}c@{}}  $t_{pips}$   \\ (sec) \end{tabular}} \\
\hline
\multirow{5}{*}{\begin{tabular}{@{}l@{}}$|\mathcal{P}|$ = 8 \\ \# Proc = 8\end{tabular}} & time & 0.0 & 87505 & 60567 & 21877 , 24940 & - & -  & -\\
& network & 0.41 & 1800 & 85248 & 459 , 1368 & 7236 & 56  & 7505\\
& 0.01 & 0.0073  & 748  & 61008 & 205 , 316 & 2944 & 6.3 & 3115\\
& 0.5  & 0.37 & 434 & 83202 & 121 , 218 & 2269 & 1.7 & 2475\\
& 0.6  & 0.37 & 458 & 83202 & 124 , 218 & 2862 & 1.9 & 3069\\
\hline
\multirow{3}{*}{\begin{tabular}{@{}l@{}}$|\mathcal{P}|$ = 24 \\ \# Proc = 24\end{tabular}} & 0.01 & 0.01 & 2889 & 20390 & 259 , 431 & 1746 & 273 & 2117\\
& 0.1 & 0.1 & 2748 & 22200 & 256 , 529 & 1934 & 236 & 2270\\
& 1.0 & 0.99 & 1572 & 40080 & 127 , 362 & 2279 & 54 & 2447\\
\hline
\multirow{4}{*}{\begin{tabular}{@{}l@{}}$|\mathcal{P}|$ = 48 \\ \# Proc = 24\end{tabular}} & 0.01 & 0.01 & 5472 & 10194 & 245 , 482 & 2029 & 1769 & 3954\\
& 0.1 & 0.1 & 4560 & 11104 & 210 , 440 & 1871 & 1031 & 3054\\
& 0.75 & 0.75 & 3182 & 17642 & 151 , 387 & 2126 & 368 & 2670\\
& 1.0 & 0.98 & 2978 & 19985 & 143 , 480 & 1826 & 298 & 2247\\
\hline
\end{tabular}
\caption{PIPS-NLP results for different graph partitions.}
\label{table:pips_results}
\end{center}
\end{table}

\subsection{Overlapping Schwarz Decomposition for DC OPF}\label{sec:dcopf}
This case study demonstrates how we can use graph partitioning and topology manipulation to pose a complex power network problem to the Schwarz overlapping scheme.
We consider a 9,241 bus test case obtained from pglib-opf (v19.05) \cite{pglib}, which we obtain through {\tt Power}-{\tt Models.jl} \cite{powermodels}.
We denote the power grid as a network $Net(\mathcal{V},\mathcal{L}_{grid})$ containing a set of
electric buses $\mathcal{V}$ that connect power lines $\ell \in \mathcal{L}_{grid}$.  Each electric bus $i \in \mathcal{V}$ can include
generators $q \in \Omega_i$ and serves a total power load $P_i^L$. The total set of generators
is given by $\Omega$ such that $\bigcup_{i\in\mathcal{V}}\Omega_i = \Omega$.  The network is described by the DC power flow
equations \cite{Sun2018} given by \eqref{eq:dcopf} where \eqref{eq:dcopf_objective} seeks to minimize the total
generation cost and voltage angle difference where $\beta$ is a regularization parameter. \eqref{eq:power_balance} enforces energy conservation, and \eqref{eq:power_limits} and \eqref{eq:angle_limits} denote limits
on power generation and voltage angles.  In this formulation, $v_i$ is the bus voltage angle for each bus $i \in \mathcal{V}$, $P_q$ is the power
generation from generator $q$ with cost coefficients $c_{q,1}$ and $c_{q,2}$ and lower and upper limits $\underline{P_q}$ and $\overline{P_q}$,
and $v_{\ell,j}$ and $v_{\ell,k}$ are the inlet and outlet
voltage angles on $\ell$ with ramp limit $\overline{v_{\ell}}$. We also define $\mathcal{L}_{rec}(i)$ as the set of power links
received by bus $i$ and $\mathcal{L}_{snd}(i)$ as the set of links sent from bus $i$.
The power flow on line $\ell$ is defined by \eqref{eq:power_flow},
where $Y_\ell$ is the branch admittance for line $\ell$, $v_{\ell,j}$ is the source bus voltage angle
and $v_{\ell,k}$ is the destination bus voltage angle. $src(\ell)$ denotes the source bus of line $\ell$ and $dst(\ell)$ is
destination bus for line $\ell$. We denote the voltage angles on the set of reference buses in \eqref{eq:reference_angles}
where $\mathcal{V}^{ref}$ is the set of reference buses.

\begin{subequations}\label{eq:dcopf}
    \begin{align}
        \min_{\substack{ \{v_i\}_{i \in \mathcal{V}} \\ \{P_q\}_{q \in \Omega}}} & \sum_{q \in \Omega}
        \left( c_{q,1} P_q + c_{q,2} P_q^2  \right) + \frac{\beta}{2}\sum_{\ell \in \mathcal{L}_{grid}} (v_{\ell,j} - v_{\ell,k})^2 \label{eq:dcopf_objective}\\
        s.t. & \sum_{q \in \Omega_i} P_q + \sum_{\ell \in \mathcal{L}_{rec}(i)} P_\ell - \sum_{\ell \in \mathcal{L}_{snd}(i)} P_\ell = P_i^L, \quad i \in \mathcal{V} \label{eq:power_balance}\\
        & v_{\ell,j} = v_{src(\ell)} , v_{\ell,k} = v_{dst(\ell)}, \quad \ell \in \mathcal{L}_{grid}\\
        & P_\ell = Y_\ell(v_{\ell,j} - v_{\ell,k}), \quad \ell \in \mathcal{L}_{grid} \label{eq:power_flow}\\
        & \underline{P_q} \le P_q \le \overline{P_q} \label{eq:power_limits}, \quad q \in \Omega\\
        & -\overline{v_{\ell}} \le v_{\ell,j} - v_{\ell,k} \le \overline{v_{\ell}}, \quad \ell \in \mathcal{L}_{grid} \label{eq:angle_limits} \\
        & v_i = v_i^{ref}, \quad i \in \mathcal{V}^{ref} \label{eq:reference_angles}
    \end{align}
\end{subequations}

We construct the DC OPF model with Code Snippet \ref{code:dcopf_snippet}, which produces an optimization problem with over 100,000 variables and constraints.
The DC OPF model is partitioned using {\tt KaHyPar}; once we obtain partitions,
we create subgraphs, expand them, and solve the problem using the {\tt SchwarzSolver.jl}. Code Snippet \ref{code:dcopf_overlap_solve} demonstrates
how we carry this out using a maximum imbalance of 10\% and an overlap size of $\omega = 10$. Line \ref{line:dcopf_kahypar} imports {\tt KaHyPar} and
Line \ref{line:dcopf_hypergraph} creates a {\tt hypergraph} and {\tt ref\_map} which we use for partitioning.  Lines \ref{line:dcopf_weights_start} through
\ref{line:dcopf_weights_end} query the graph for edge weights and node sizes, Line \ref{line:dcopf_part} partitions the DC OPF problem,
and Lines \ref{line:dcopf_part_object} and \ref{line:dcopf_make_subgraphs} create a {\tt Partition} and use it to define
subgraphs for the problem.  Once we have subgraphs, we perform a subgraph expansion on Line \ref{line:dcopf_expand} and execute
the Schwarz solver on Line \ref{line:dcopf_solve}. We also show how to tell the solver how to treat
linking constraints.  We provide the keyword argument {\tt primal\_links} with the vector of power flow linking constraints and we provide the
keyword {\tt dual\_links} with a vector of voltage angle linking constraints to denote how subproblems are formulated;
the constraints that are specified as {\tt primal\_links} are treated as direct constraints while the constraints in {\tt dual\_links} are incorporated as dual penalty
as described in \eqref{eq:Schwarz_subproblem}.

\begin{minipage}[]{0.9\linewidth}
\begin{scriptsize}
\lstset{language=Julia,breaklines = true}
\begin{lstlisting}[caption = Partitioning and Formulating Overlapping Subproblems,label = {code:dcopf_overlap_solve}]
using KaHyPar                                           |\label{line:dcopf_kahypar}|
using Ipopt
using SchwarzSolver

#Get the hypergraph representation of the gas network
hypergraph,ref_map = gethypergraph(dcopf)        |\label{line:dcopf_hypergraph}|

#Setup node and edge weights
n_vertices = length(vertices(hypergraph))                |\label{line:dcopf_weights_start}|
node_weights = [num_variables(node) for node in all_nodes(dcopf)]
edge_weights = [num_link_constraints(edge) for edge in all_edges(dcopf)]   |\label{line:dcopf_weights_end}|

#Use KaHyPar to partition the hypergraph
node_vector = KaHyPar.partition(hypergraph,4,configuration = :edge_cut,
imbalance = 0.1, node_weights = node_weights,edge_weights = edge_weights)      |\label{line:dcopf_part}|

#Create a partition object
partition = Partition(dcopf,node_vector,ref_map)           |\label{line:dcopf_part_object}|

#Setup subgraphs based on partition
make_subgraphs!(dcopf,partition)               |\label{line:dcopf_make_subgraphs}|

distance = 10
subgraphs = getsubgraphs(dcopf)

#Expand the subgraphs
sub_expand = expand.(Ref(dcopf),subgraphs,Ref(distance))      |\label{line:dcopf_expand}|

ipopt = Ipopt.Optimizer

schwarz_solve(dcopf,sub_expand,primal_links = power_links,  |\label{line:dcopf_solve}|
dual_links = [angle_i;angle_j]],
sub_optimizer = ipopt,max_iterations = 100,tolerance = 1e-3)
\end{lstlisting}
\end{scriptsize}
\end{minipage}

\begin{figure}[]
\centering
\begin{subfigure}[t]{0.24\textwidth}
    \centering
    \includegraphics[width = 4cm]{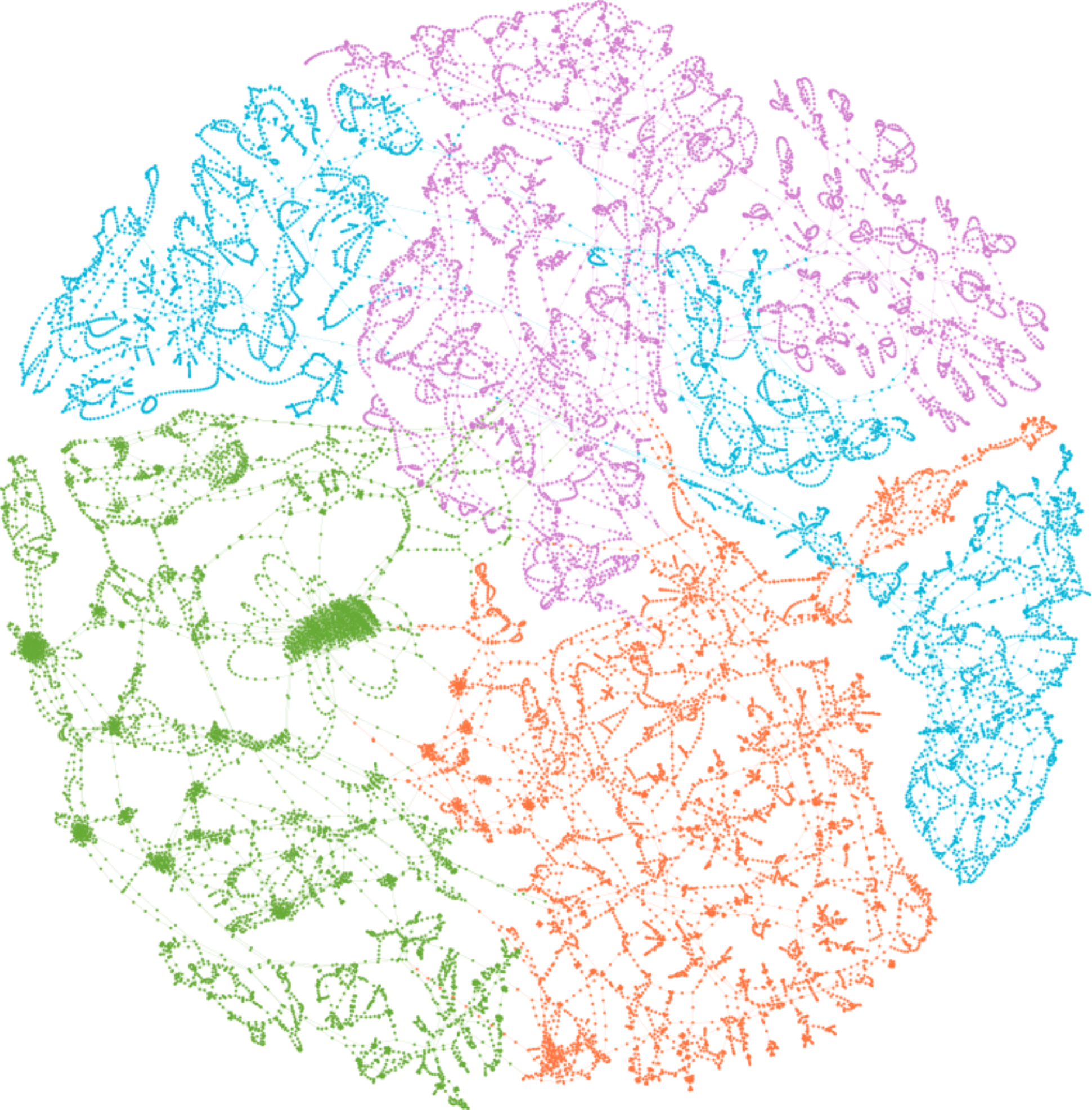}
    \caption*{DC OPF Partitions}
\end{subfigure}
\begin{subfigure}[t]{0.74\textwidth}
    \centering
    \begin{subfigure}[t]{0.4\textwidth}
        \includegraphics[width = 4cm]{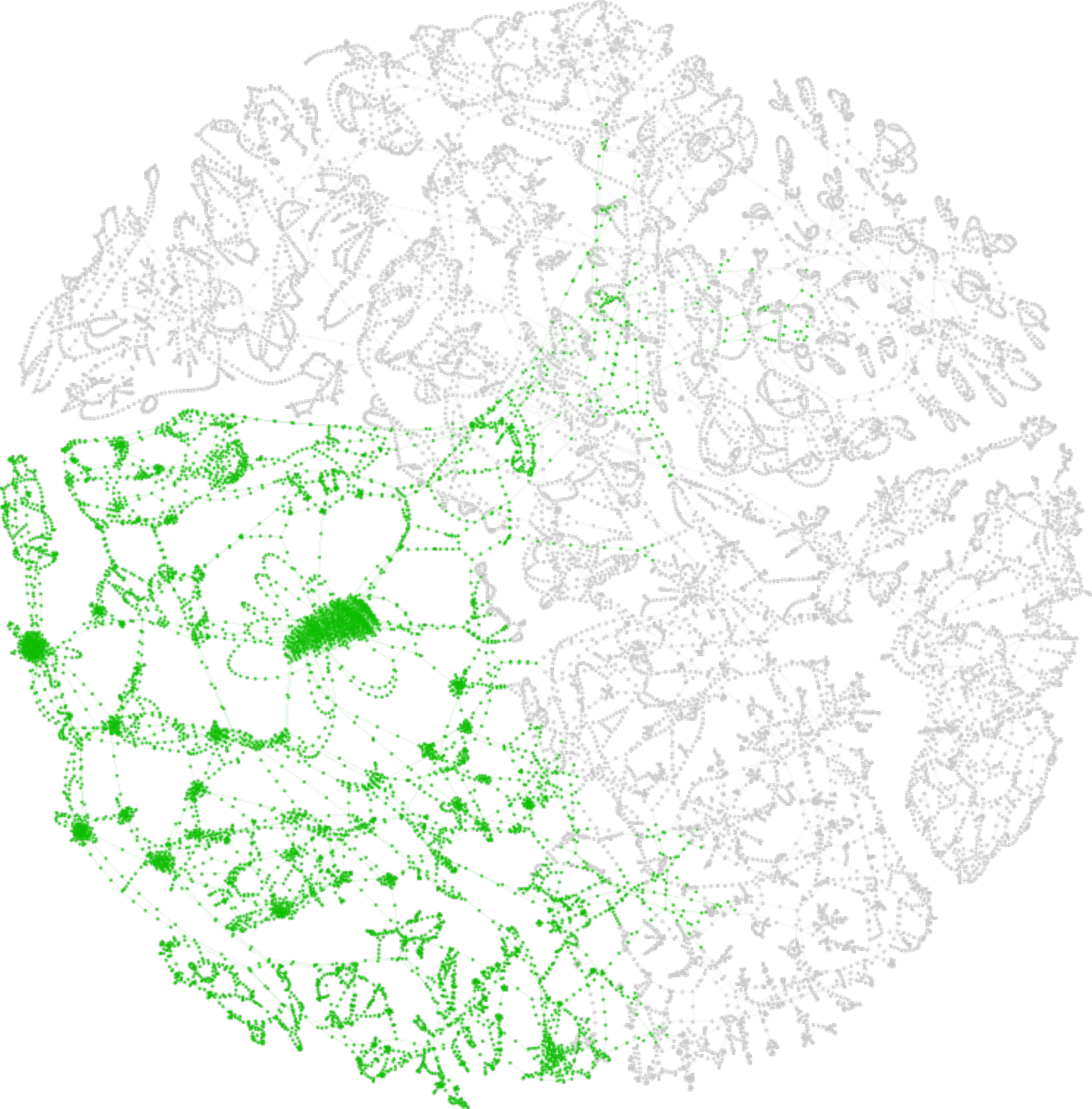}
        \caption*{Partition 1 Overlap}
    \end{subfigure}
    \begin{subfigure}[t]{0.4\textwidth}
        \centering
        \includegraphics[width = 4cm]{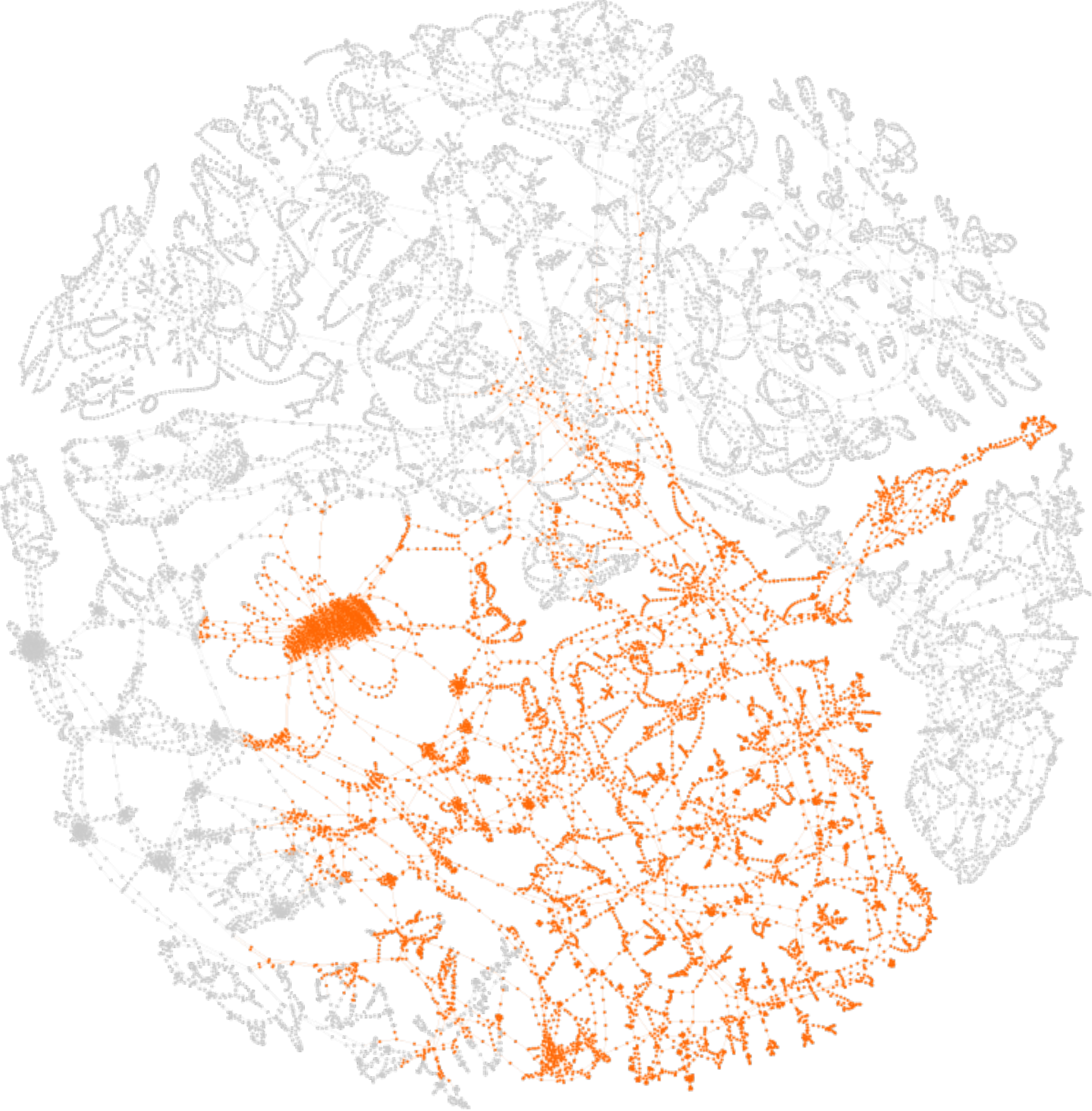}
        \caption*{Partition 2 Overlap}
    \end{subfigure}
    \\
    \begin{subfigure}[t]{0.4\textwidth}
        \centering
        \includegraphics[width = 4cm]{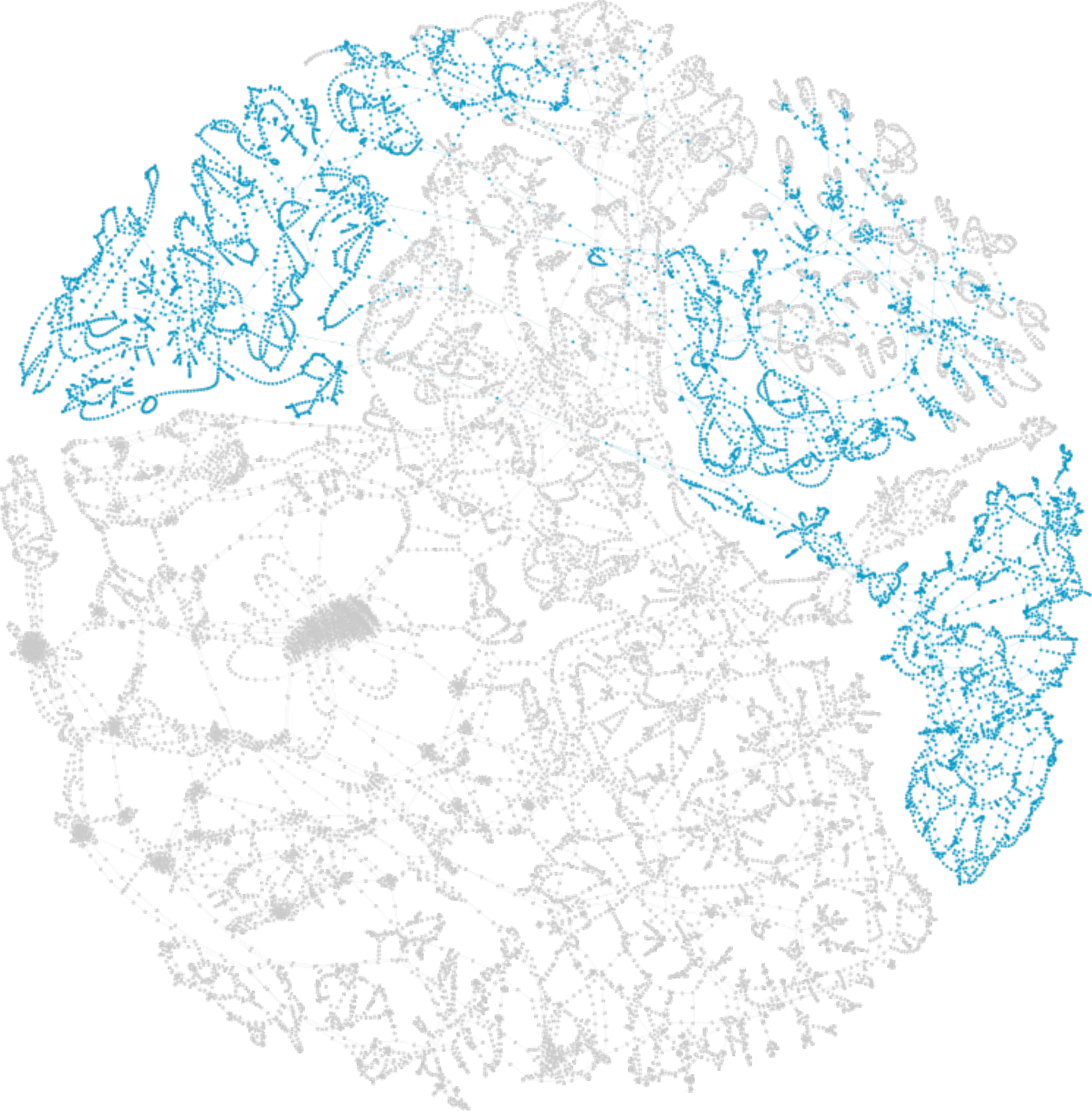}
        \caption*{Partition 3 Overlap}
    \end{subfigure}
    \begin{subfigure}[t]{0.4\textwidth}
        \centering
        \includegraphics[width = 4cm]{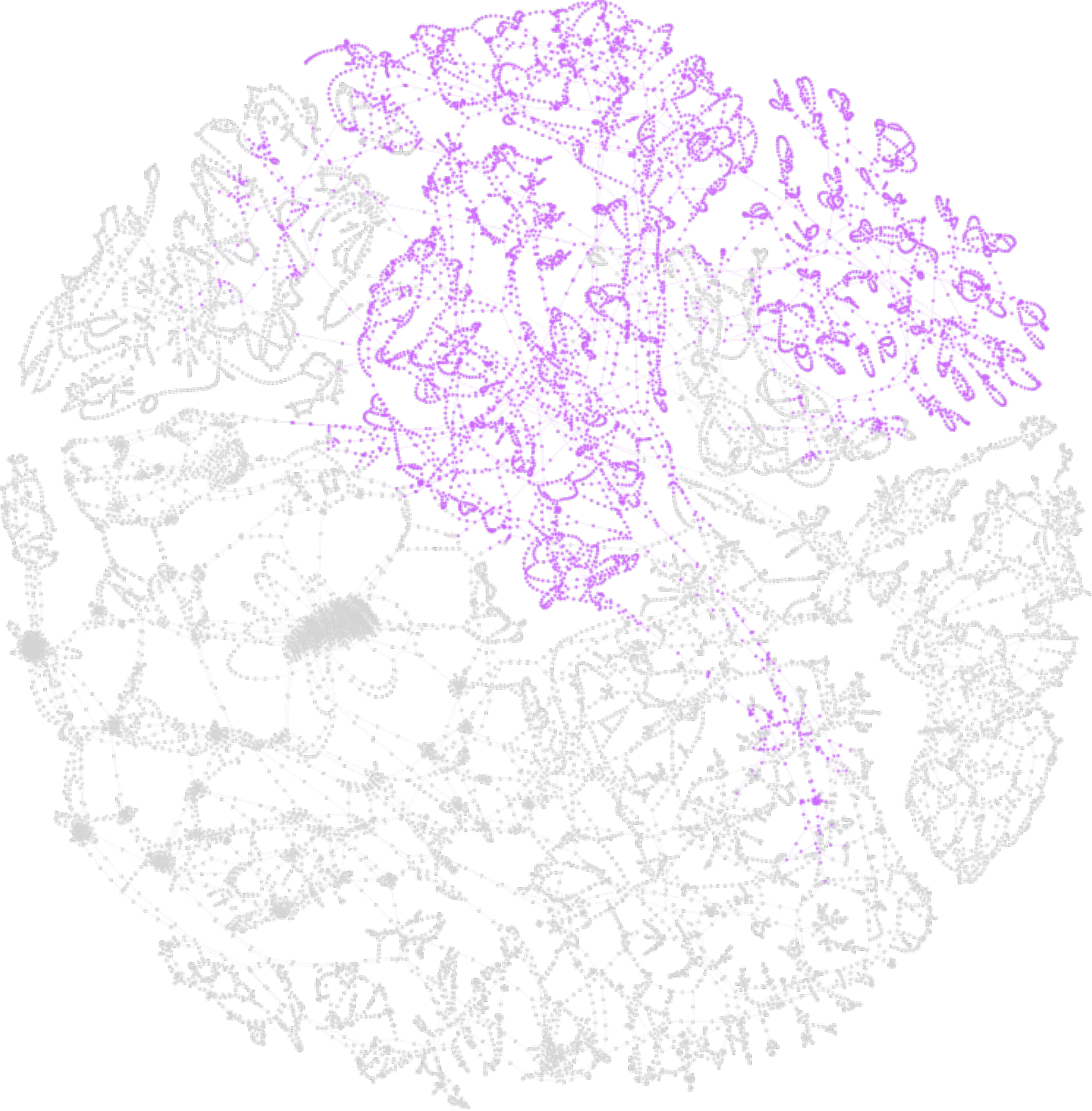}
        \caption*{Partition 4 Overlap}
    \end{subfigure}
\end{subfigure}
\caption{Depiction of DC OPF problem with four partitions. The original calculated partitions with $\epsilon_{max} = 0.1$ (left) and
the corresponding overlap partitions with $\omega = 10$ (right).}
\label{fig:dcopf_parts}
\end{figure}

\begin{figure}[]
\centering
\begin{subfigure}[t]{0.32\textwidth}
    \centering
    \includegraphics[scale = 0.25]{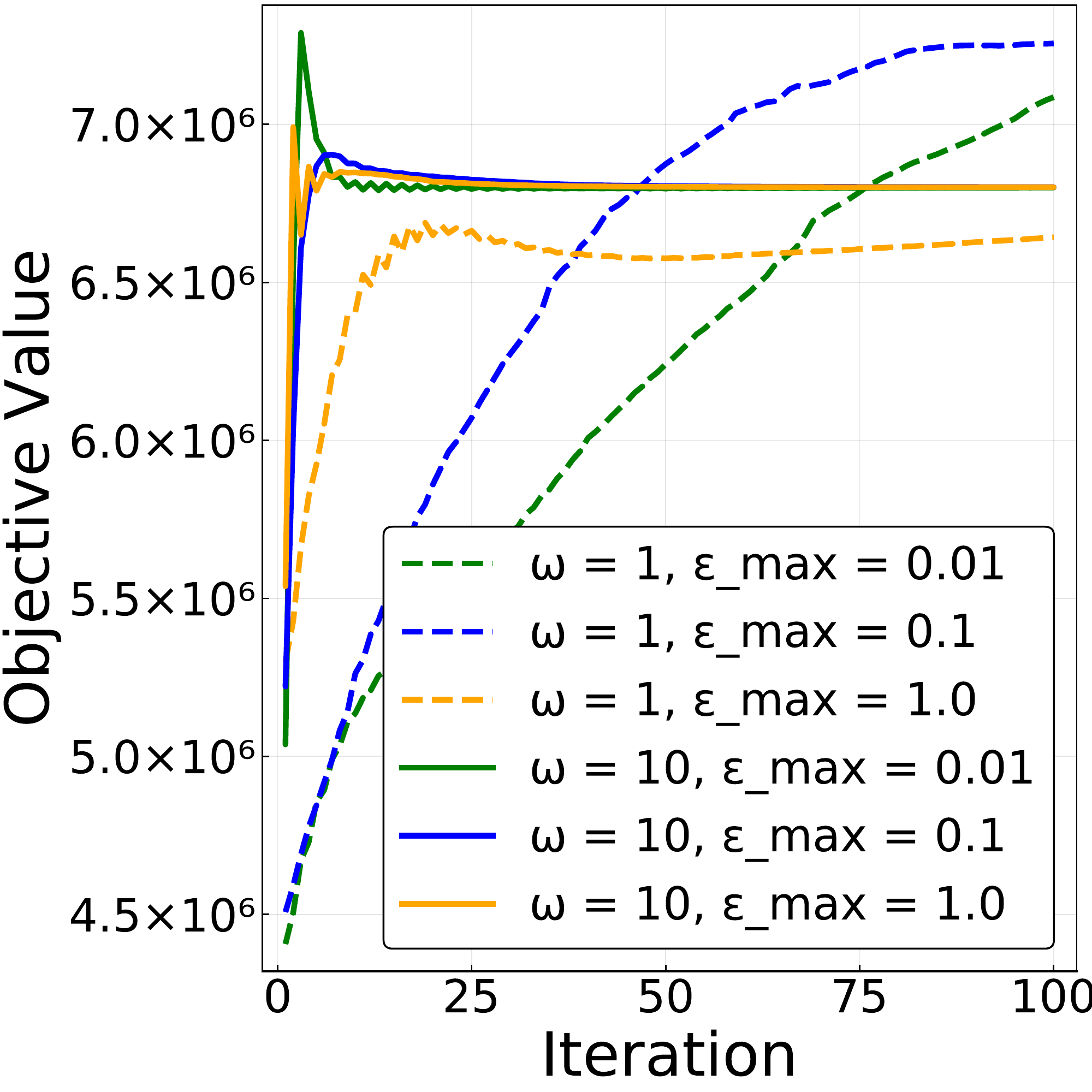}
    \caption*{}
\end{subfigure}
\begin{subfigure}[t]{0.32\textwidth}
    \centering
    \includegraphics[scale = 0.25]{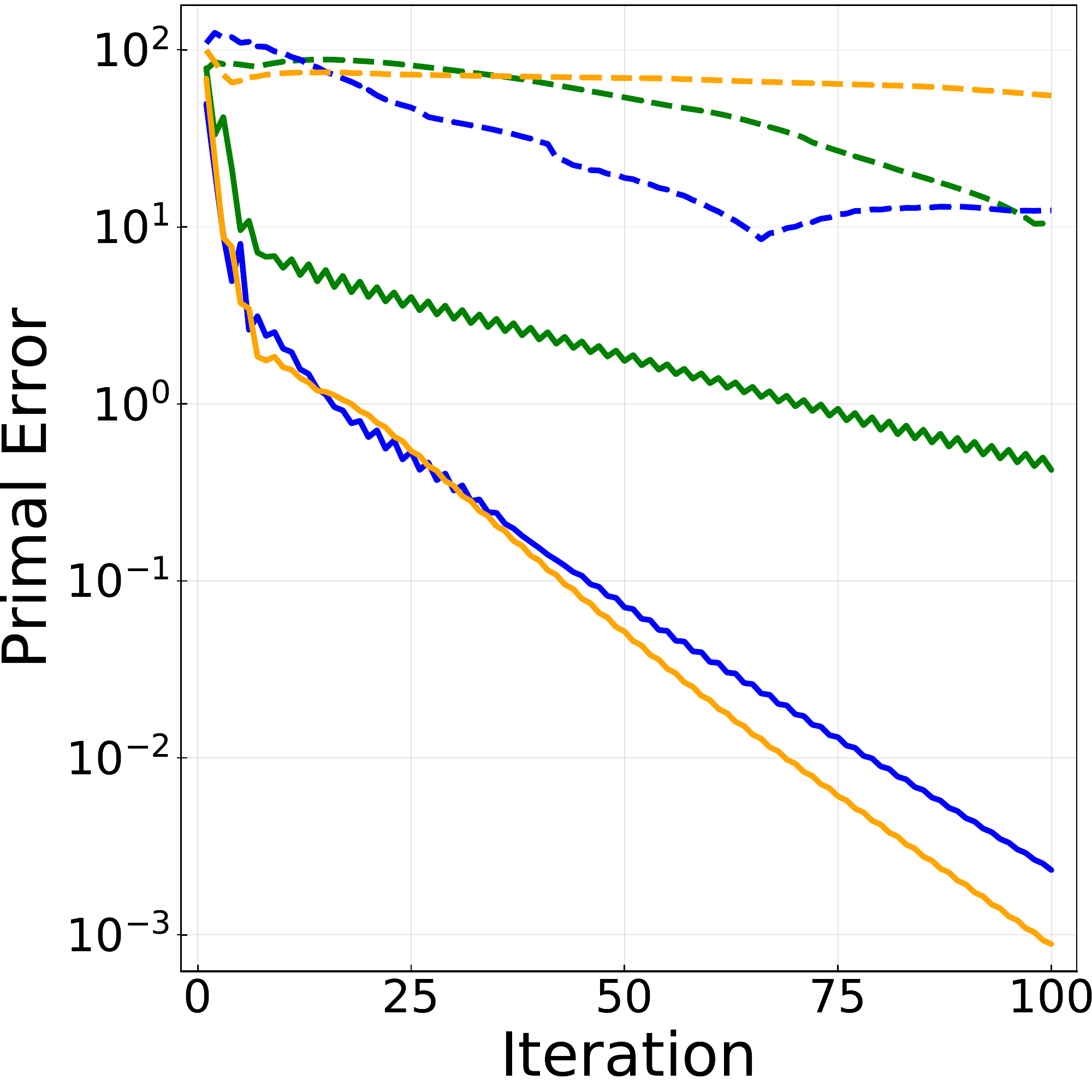}
    \caption*{}
\end{subfigure}
\begin{subfigure}[t]{0.32\textwidth}
    \centering
    \includegraphics[scale = 0.25]{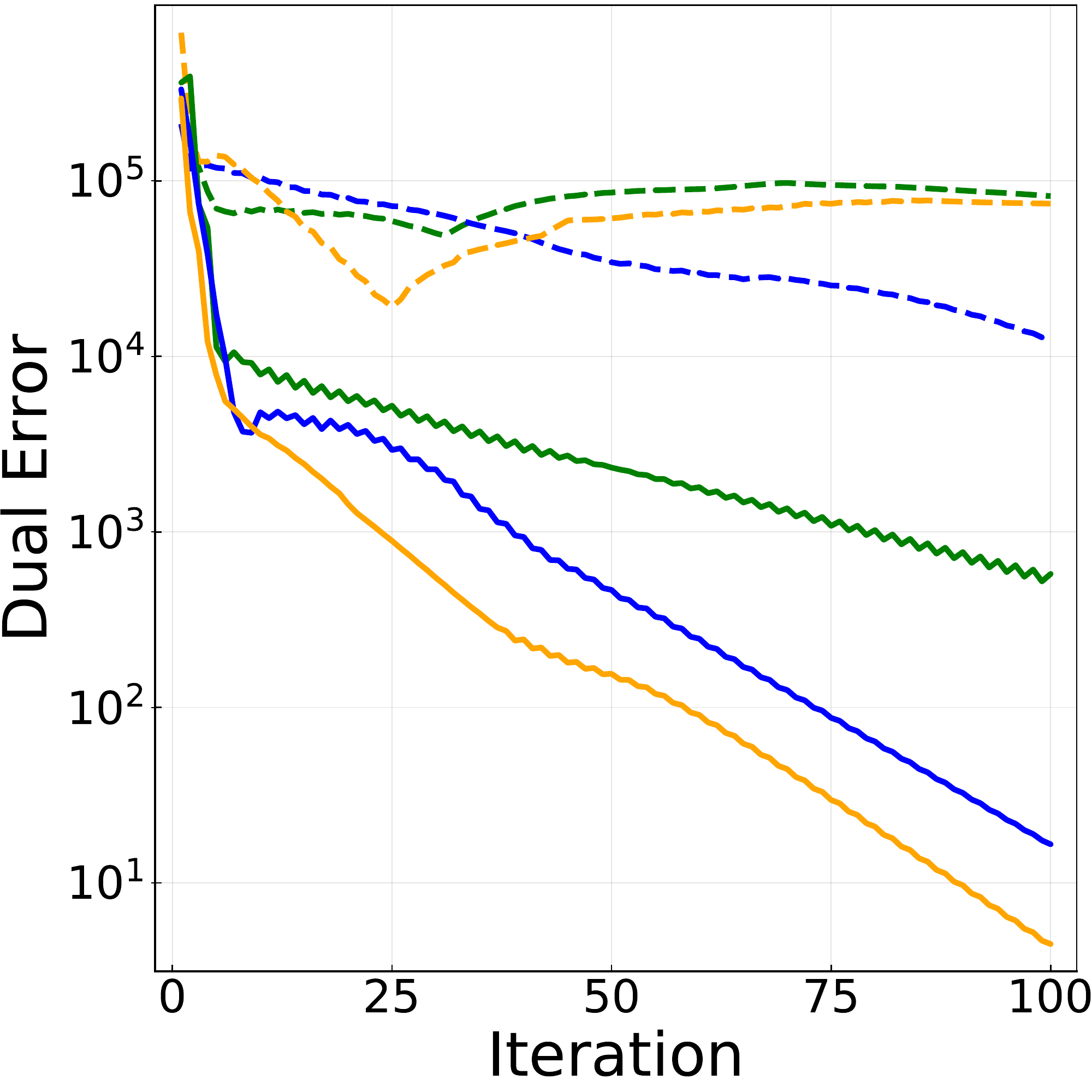}
    \caption*{}
\end{subfigure}
\caption{Comparison of Schwarz algorithm for different values of overlap $\omega$ and maximum partition imbalance $\epsilon_{max}$}
\label{fig:dcopf_results}
\end{figure}

Figure \ref{fig:dcopf_parts} depicts the original and overlapping  partitions obtained from Snippet \ref{code:dcopf_overlap_solve}. These are visualized using
{\tt Gephi}. We experiment with different values for maximum imbalance and overlap and obtain the results
in Figure \ref{fig:dcopf_results}. We see that the Schwarz algorithm fails to converge with an overlap value of one (for any imbalance value) which is consistent with the convergence analysis in
\cite{Shin2020b}, but a sufficient overlap of 10 produces smooth convergence (for each imbalance value). We also found that larger partition
imbalance results in faster convergence (with sufficient overlap) which is likely due to the smaller edge cut and fewer linking constraints that need to be satisfied. We thus see that the trade-offs of imbalance and coupling are complex and differ under different problems and algorithms.

\FloatBarrier

\section{Outlook and Extensions for {\tt Plasmo.jl}}\label{sec:conclusion}
We presented a graph-based modeling abstraction for optimization that we call an {\tt OptiGraph}.  We showcased its implementation in {\tt Plasmo.jl} and demonstrated how this
provides flexibility to use different graph analysis and visualization tools. We also showed how these capabilities facilitate the use of decomposition strategies at both the linear
algebra and problem levels. There are broad opportunities to refine and expand the capabilities of the proposed abstraction and of {\tt Plasmo.jl}. The most evident next step is to
develop interfaces to a wider range of decomposition solvers. As such, we are currently developing interfaces
to {\tt DSP} \cite{Kim2016} (a Benders and Lagrangian dual decomposition solver for
stochastic optimization),  and {\tt SNGO} \cite{Cao2019} (a Julia-based global solver for
nonlinear stochastic programs). We are also interested in new parallel interior point solvers that use recent Schur decomposition strategies
to solve massive energy system models \cite{Onheit2019}. In this regard, the graph could be a natural interface for {\tt BELTISTOS} \cite{Kourounis2019},
a nonlinear solver that exploits specialized time decomposition structures that arise in Schur decomposition. These capabilities will enable solution of a
wide range of application problems.

As ongoing work we are also developing \emph{parallel modeling} capabilities to create graph objects in parallel. Parallel modeling becomes important in instances
when the model itself is too large to fit in memory, such as with stochastic programs with many scenarios.
Scenario-based optimization problems have been modeled in parallel using {\tt StructJuMP} and {\tt MPI}, but more complex structures (such as those with many linking constraints) are
more challenging to implement in a parallel modeling framework.  In this context, the Parallel Structured Model Generator ({\tt PSMG}) \cite{Qiang2014} is
the only framework we know of that performs parallel model generation which works with structured conveying modeling language ({\tt SML}) \cite{Colombo2009}.
{\tt PSMG} focuses on efficient distributed
memory design to achieve parallel generation, but it requires solvers to perform subproblem distribution and load balancing based on a well-defined model structure.
In contrast, we hope to facilitate the modeling, partitioning, and parallel generation tasks in a more systematic interactive framework using the {\tt OptiGraph}.

Finally, we demonstrated how to partition graphs using powerful tools such as {\tt Metis} and {\tt KaHyPar}.  Exploiting
such graph partitioning tools opens many possibilities for using decomposition - based solvers, but it is limited to expressing problem characteristics
strictly in the form of edge weights and node sizes. We are interested in developing more customized graph partitioning algorithms
that work directly on the graph attributes (such as accounting for integer variables). We are also interested in developing algorithms that
calculate overlap directly as opposed to performing subgraph expansion.
Such modeling advances can provide new capabilities that open diverse algorithmic development and offer interesting directions
for future research.

\section*{Acknowledgements}
This material is based on work supported by the U.S. Department of Energy (DOE), Office of Science, under Contract No. DE-AC02-06CH11357 as well as the
DOE Office of Electricity Delivery and Energy Reliability’s Advanced Grid Research and Development program (AGR\&D).
This work was also partially supported by the U.S. Department of Energy grant DE-SC0014114.
We also acknowledge partial support from the National Science Foundation under award NSF-EECS-1609183. The authors acknowledge Yankai Cao for
his work developing the original {\tt PIPS}-{\tt NLP} interface for {\tt Plasmo.jl}.  The authors also acknowledge Michel Schanen from Argonne National Laboratory for his
assistance using and debugging {\tt PIPS-NLP} as well as
Kibaek Kim from Argonne National Laboratory for his help in the conceptualization of {\tt Plasmo.jl}.

\appendix

\section*{Appendix: Code Snippets}

\paragraph{Junction OptiGraph Implementation}\mbox{}\\
The junction {\tt OptiGraph} is implemented in Code Snippet \ref{code:junction_construction_snippet}.
We define the function {\tt create\_junction\_model} on Line \ref{line:junction_function} which accepts junction specific {\tt data} and the number of time
periods {\tt nt}.  We create the OptiGraph {\tt graph} on Lines \ref{line:junction_create_start} through \ref{line:junction_create_end} where
we add a node for each time interval and then create the variables and constraints for each node in a loop.  We also use the {\tt JuMP} specific
{\tt @expression} macro to refer to expressions for total gas supplied, total gas delivered, and total cost for convenience.  The junction model
is returned from the function on Line \ref{line:junction_return}.

\begin{minipage}[]{0.9\linewidth}
\begin{scriptsize}
\lstset{language=Julia,breaklines = true}
\begin{lstlisting}[caption = {Creating a gas junction graph},label = {code:junction_construction_snippet}]
#Define function to create junction OptiGraph
function create_junction_model(data,nt)            |\label{line:junction_function}|
   graph = OptiGraph()                           |\label{line:junction_create_start}|

   #Add OptiNode for each time interval
   @optinode(graph,nodes[1:nt])

   #query number of supply and demands on the junction
   n_demands = length(data[:demand_values])
   n_supplies = length(data[:supplies])

   #Loop and create variables, constraints, and objective for each OptiNode
   for (i,node) in enumerate(nodes)
       @variable(node, data[:pmin] <= pressure <= data[:pmax], start = 60)
       @variable(node, 0 <= fgen[1:n_supplies] <= 200, start = 10)
       @variable(node, fdeliver[1:n_demands] >= 0)
       @variable(node, fdemand[1:n_demands] >= 0)

       @constraint(node,[d = 1:n_demands],fdeliver[d] <= fdemand[d])

       @expression(node, total_supplied, sum(fgen[s] for s = 1:n_supplies))
       @expression(node, total_delivered,sum(fdeliver[d] for d = 1:n_demands))
       @expression(node, total_delivercost,sum(1000*fdeliver[d] for d = 1:n_demands))

       @objective(node,Min,total_delivercost)
   end                                              |\label{line:junction_create_end}|

   #Return the graph
   return graph                                    |\label{line:junction_return}|
end
\end{lstlisting}
\end{scriptsize} %
\end{minipage}

\paragraph{Pipeline OptiGraph Implementation}\mbox{}\\
The pipeline {\tt OptiGraph} is given by Code Snippet \ref{code:pipeline_construction_snippet}.
We again use a function ({\tt create\_pipeline\_model}) to create the graph
and we unpack the input data and create nodes, variables, and constraints on Lines \ref{line:pipe_start} through \ref{line:pipe_end}.
We define the set of {\tt OptiNodes} on Line \ref{line:pipe_grid} as {\tt grid}, which we use to refer to {\tt OptiNodes} using time and space indices.
To capture the space-time coupling of each pipeline we use linking constraints on Lines \ref{line:pipe_linkconstraints_start} through
\ref{line:pipe_linkconstraints_end} to define the finite difference scheme, the steady-state condition, and the line pack constraint.

\begin{minipage}[]{0.9\linewidth}
\begin{scriptsize}
\lstset{language=Julia,breaklines = true}
\begin{lstlisting}[caption = Creating a pipeline graph, label = code:pipeline_construction_snippet]
function create_pipeline_model(data,nt,nx)
    #Unpack data
    c1 = data[:c1]; c2 = data[:c2]; c3 = data[:c3]                      |\label{line:pipe_start}|
    dx = data[:pipe_length] / (nx - 1)

    #Create pipeline OptiGraph
    graph = OptiGraph()

    #Create grid of OptiNodes
    @optinode(graph,grid[1:nt,1:nx])                                           |\label{line:pipe_grid}|

    #Create variables on each node in the grid
    for node in grid
        @variable(node, 1 <= px <= 100)
        @variable(node, 0 <= fx <= 100)
        @variable(node,slack >= 0)
        @NLnodeconstraint(node, slack*px - c3*fx*fx == 0)
    end

    #Setup dummy variable referenes
    @expression(graph,fin[t=1:nt],grid[:,1][t][:fx])
    @expression(graph,fout[t=1:nt],grid[:,end][t][:fx])
    @expression(graph,pin[t=1:nt],grid[:,1][t][:px])
    @expression(graph,pout[t=1:nt],grid[:,end][t][:px])
    @expression(graph,linepack[t=1:nt],c2/A*sum(grid[t,x][:px]*dx for x in 1:nx-1))           |\label{line:pipe_end}|

    #Finite differencing.  Backward difference in time from t, Forward difference in space from x.
    @linkconstraint(graph, press[t=2:nt,x=1:nx-1],                                             |\label{line:pipe_linkconstraints_start}|
    (grid[t,x][:px]-grid[t-1,x][:px])/dt +
    c1*(grid[t,x+1][:fx] - grid[t,x][:fx])/dx == 0)

    @linkconstraint(graph, flow[t=2:nt,x=1:nx-1],(grid[t,x][:fx] -
    grid[t-1,x][:fx])/dt == -c2*(grid[t,x+1][:px] -
    grid[t,x][:px])/dx - grid[t,x][:slack])

    #Initial steady state
    @linkconstraint(graph,ssflow[x=1:nx-1],grid[1,x+1][:fx] - grid[1,x][:fx] == 0)
    @linkconstraint(graph,sspress[x = 1:nx-1], -c2*(grid[1,x+1][:px] -
    grid[1,x][:px])/dx - grid[1,x][:slack] == 0)

    #Refill pipeline line pack
    @linkconstraint(graph,linepack[end] >= linepack[1])                        |\label{line:pipe_linkconstraints_end}|
    return graph
end
\end{lstlisting}
\end{scriptsize}
\end{minipage}

\paragraph{Compressor OptiGraph Implementation}\mbox{}\\
The compressor {\tt OptiGraph} construction is fairly straight-forward and given by
\ref{code:compressor_modelgraph}. Line \ref{line:compressor_function}
defines the function \\ {\tt create\_compressor\_model} which takes {\tt data} and the number of time periods {\tt nt}.
Line \ref{line:compressor_optigraph} creates the {\tt OptiGraph} for the compressor model, Line \ref{line:compressor_nodes}
creates an {\tt OptiNode} for each time point, and Lines \ref{line:compressor_start} through \ref{line:compressor_stop} define the model variables
and constraints.  Lines \ref{line:compressor_exp1} and \ref{line:compressor_exp2} define helpful expressions for the flow in and out of the compressor and Line \ref{line:compressor_return}
returns the {\tt OptiGraph}.

\begin{minipage}[]{0.9\linewidth}
\begin{scriptsize}
\lstset{language=Julia,breaklines = true}
\begin{lstlisting}[caption = Creating a compressor graph,label = {code:compressor_modelgraph}]
function create_compressor_model(data,nt)    |\label{line:compressor_function}|
    #Create compressor OptiGraph
    graph = OptiGraph()                     |\label{line:compressor_optigraph}|
    @optinode(graph,nodes[1:nt])                         |\label{line:compressor_nodes}|

    #Setup variables, constraints, and objective
    for node in nodes                                        |\label{line:compressor_start}|
        @variable(node, 1 <= psuction <= 100)
        @variable(node, 1 <= pdischarge <= 100)
        @variable(node, 0 <= power <= 1000)
        @variable(node, flow >= 0)
        @variable(node, 1 <= eta <= 2.5)
        @NLnodeconstraint(node, pdischarge == eta*psuction)
        @NLnodeconstraint(node, power == c4*flow*((pdischarge/psuction)^om-1) )
        @objective(node, Min, cost*power*(dt/3600.0))                |\label{line:compressor_stop}|
    end

    #Create references for flow in and out
    @expression(graph,fin[t=1:nt],nodes[t][:flow])               |\label{line:compressor_exp1}|
    @expression(graph,fout[t=1:nt],nodes[t][:flow])              |\label{line:compressor_exp2}|

    #Return compressor OptiGraph
    return graph                                     |\label{line:compressor_return}|
end
\end{lstlisting}
\end{scriptsize}
\end{minipage}

\paragraph{Network OptiGraph Implementation}\mbox{}\\
The network {\tt OptiGraph} couples the other component graphs at a higher level and is implemented in Code Snippet \ref{code:gas_network_construction}.
We create the graph {\tt gas\_network} on Line \ref{line:create_gas_network} and add the component subgraphs
on Lines \ref{line:add_subgraph_network_start} through \ref{line:add_subgraph_network_end}.
Once we have the multi-level subgraph structure we create linking constraints at the {\tt gas\_network} level on Lines
\ref{line:network_link_start} through \ref{line:network_link_end} which impose the junction conservation and
boundary conditions for the network links (pipelines and compressors).

\begin{minipage}[]{0.9\linewidth}
\begin{scriptsize}
\lstset{language=Julia,breaklines = true}
\begin{lstlisting}[caption = Formulating the complete gas network graph,label = {code:gas_network_construction}]
#Create OptiGraph for entire gas network
gas_network = OptiGraph()                          |\label{line:create_gas_network}|

#Add every device graph to the network
for pipe in pipelines                               |\label{line:add_subgraph_network_start}|
    add_subgraph!(gas_network,pipe)
end
for compressor in compressors
    add_subgraph!(gas_network,compressor)
end
for junction in junctions
    add_subgraph!(gas_network,junction)
end                                                 |\label{line:add_subgraph_network_end}|

#Link pipelines in gas network
for pipe in pipelines                               |\label{line:network_link_start}|
    junction_from,junction_to = pipe_map[pipe]
    @linkconstraint(gas_network,[t = 1:nt],pipe[:pin][t] ==
    junction_from[:pressure][t])
    @linkconstraint(gas_network,[t = 1:nt],pipe[:pout][t] ==
    junction_to[:pressure][t])
end

#Link compressors in gas network
for compressor in compressors
    junction_from,junction_to = compressor_map[compressor]
    @linkconstraint(gas_network,[t = 1:nt],compressor[:pin][t] ==
    junction_from[:pressure][t])
    @linkconstraint(gas_network,[t = 1:nt],compressor[:pout][t]  ==
    junction_to[:pressure][t])
end

#Link junctions in gas network
for junction in junctions
    devices_in = junction_map_in[junction]
    devices_out = junction_map_out[junction]

    flow_in = [sum(device[:fout][t] for device in devices_in) for t = 1:nt]
    flow_out = [sum(device[:fin][t] for device in devices_out) for t = 1:nt]

    total_supplied = [junction[:total_supplied][t] for t = 1:nt]
    total_delivered = [junction[:total_delivered][t] for t = 1:nt]

    @linkconstraint(gas_network,[t = 1:nt], flow_in[t] - flow_out[t] +
    total_supplied[t] - total_delivered[t] == 0)
end                                                  |\label{line:network_link_end}|
\end{lstlisting}\label{fig:network_link_snippet}
\end{scriptsize}
\end{minipage}

\paragraph{DC OPF Model Implementation}\mbox{}\\
The DC optimal power flow model is implemented in Code Snippet \ref{code:dcopf_snippet}.
Line \ref{line:load_powergrid} loads the bus system using {\tt PowerModels.jl}, Line \ref{line:create_dcopf} creates the graph,
Lines \ref{line:create_buses} and \ref{line:create_lines} create nodes
for each bus and line in the network, and Line \ref{line:load_data} assigns relevant data to each bus and line node using a {\tt load\_data!} function.  The
DC OPF model is constructed in the same fashion as described in earlier examples where Lines \ref{line:line_models_start} through \ref{line:voltage_angle_links}
setup models on each node and add linking constraints corresponding to power conservation and voltage angle connections.

\begin{minipage}[]{0.9\linewidth}
\begin{scriptsize}
\lstset{language=Julia,breaklines = true}
\begin{lstlisting}[caption = Creating the DC OPF Problem,label = {code:dcopf_snippet}]
    include("ogplib_data.jl")                               |\label{line:load_powergrid}|

    #Create graph based on network
    grid = OptiGraph()                                     |\label{line:create_dcopf}|

    #Create bus and power line nodes
    @optinode(grid,buses[1:N_buses])                             |\label{line:create_buses}|
    @optinode(grid,lines[1:N_lines])                             |\label{line:create_lines}|

    #Load data and setup bus and line mappings
    load_data!(lines,buses)                                   |\label{line:load_data}|

    #Setup line OptiNodes
    for line in lines                                        |\label{line:line_models_start}|
        bus_from = line_map[line][1]
        bus_to = line_map[line][2]
        delta = line.ext[:angle_rate]
        @variable(line,va_i,start = 0)
        @variable(line,va_j,start = 0)
        @variable(line,flow,start = 0)
        @constraint(line,flow == B[line]*(va_i - va_j))
        @constraint(line,delta <= va_i - va_j <= -delta)
        @objective(line,Min,gamma*(va_i - va_j)^2)             |\label{line:line_models_end}|
    end

    #Setup bus OptiNodes                                     |\label{line:bus_models_start}|
    power_links = []
    for bus in buses
        va_lower = bus.ext[:va_lower]; va_upper = bus.ext[:va_upper]
        gen_lower = bus.ext[:gen_lower]; gen_upper = bus.ext[:gen_upper]

        @variable(bus,va_lower <= va <= va_upper)
        @variable(bus,P[j=1:ngens[bus]])
        @constraint(bus,[j = 1:ngens[bus]]gen_lower[j] <= P[j] <= gen_upper[j])

        lines_in = node_map_in[bus] ; lines_out = node_map_out[bus]

        @variable(bus,power_in[1:length(lines_in)])
        @variable(bus,power_out[1:length(lines_out)])

        @constraint(bus,power_balance, sum(bus[:P][j] for j=1:ngens[bus]) -
        sum(power_in) + sum(power_out) - load_map[bus] == 0)
        @objective(bus,Min,sum(bus.ext[:c1][j]*bus[:P][j] +
        bus.ext[:c2][j]*bus[:P][j]^2 for j = 1:ngens[bus]))

        #Link power flow
        p1 = @linkconstraint(grid,[j = 1:length(lines_in)],          |\label{line:bus_models_links}|
        bus[:power_in][j] == lines_in[j][:flow])
        p2 = @linkconstraint(grid,[j = 1:length(lines_out)],
        bus[:power_out][j] == lines_out[j][:flow])                   |\label{line:bus_models_end}|
        push!(power_links,p1) ; push!(power_links,p2)
    end

    #Link voltage angles
    @linkconstraint(grid,angle_i[line = lines],line[:va_i] ==        |\label{line:voltage_angle_links}|
    line_map[line][1][:va])
    @linkconstraint(grid,angle_j[line = lines],line[:va_j] ==
    line_map[line][2][:va])
\end{lstlisting}
\end{scriptsize}
\end{minipage}

\clearpage

\bibliography{algebraicgraphs}

\end{document}